\documentclass[12pt]{book}
\usepackage{amsmath,amssymb}
\makeindex

\begin{document}

\title{Selberg zeta functions for spaces of higher rank}
\author{Anton Deitmar}

\date{}
\maketitle

\tableofcontents

\setlength{\parskip}{7pt}

\def \1{{\bf 1}}
\def \a{{{\mathfrak a}}}
\def \ad{{\rm ad}}
\def \al{\alpha}
\def \ar{{\alpha_r}}
\def \A{{\mathbb A}}
\def \Ad{{\rm Ad}}
\def \Aut{{\rm Aut}}
\def \b{{{\mathfrak b}}}
\def \bs{\backslash}
\def \B{{\cal B}}
\def \c{{\mathfrak c}}
\def \cent{{\rm cent}}
\def \C{{\mathbb C}}
\def \CA{{\cal A}}
\def \CB{{\cal B}}
\def \CC{{\cal C}}
\def \CE{{\cal E}}
\def \CF{{\cal F}}
\def \CG{{\cal G}}
\def \CH{{\cal H}}
\def \CHC{{\cal HC}}
\def \CL{{\cal L}}
\def \CM{{\cal M}}
\def \CN{{\cal N}}
\def \CP{{\cal P}}
\def \CQ{{\cal Q}}
\def \CO{{\cal O}}
\def \CS{{\cal S}}
\def \CT{{\cal T}}
\def \CV{{\cal V}}
\def \det{{\rm det}}
\def \df{\ \begin{array}{c} _{\rm def}\\ ^{\displaystyle =}\end{array}\ }
\def \diag{{\rm diag}}
\def \dist{{\rm dist}}
\def \End{{\rm End}}
\def \eps{\varepsilon}
\def \Fx{{\mathfrak x}}
\def \FX{{\mathfrak X}}
\def \g{{{\mathfrak g}}}
\def \ga{\gamma}
\def \Ga{\Gamma}
\def \GL{{\rm GL}}
\def \h{{{\mathfrak h}}}
\def \Hom{{\rm Hom}}
\def \im{{\rm im}}
\def \Im{{\rm Im}}
\def \Ind{{\rm Ind}}
\def \k{{{\mathfrak k}}}
\def \K{{\cal K}}
\def \l{{\mathfrak l}}
\def \la{\lambda}
\def \lap{\triangle}
\def \La{\Lambda}
\def \m{{{\mathfrak m}}}
\def \mod{{\rm mod}}
\def \n{{{\mathfrak n}}}
\def \name{\bf}
\def \N{\mathbb N}
\def \o{{\mathfrak o}}
\def \ord{{\rm ord}}
\def \O{{\cal O}}
\def \p{{{\mathfrak p}}}
\def \ph{\varphi}
\def \prf{\noindent{\bf Proof: }}
\def \Per{{\rm Per}}
\def \q{{\mathfrak q}}
\def \qed{\ifmmode\eqno $\square$\else\noproof\vskip 12pt plus 3pt minus 9pt \fi}
 \def\noproof{{\unskip\nobreak\hfill\penalty50\hskip2em\hbox{}%
     \nobreak\hfill $\square$\parfillskip=0pt%
     \finalhyphendemerits=0\par}}
\def \Q{\mathbb Q}
\def \res{{\rm res}}
\def \R{{\mathbb R}}
\def \Re{{\rm Re \hspace{1pt}}}
\def \r{{\mathfrak r}}
\def \ra{\rightarrow}
\def \rank{{\rm rank}}
\def \Rep{{\rm Rep}}
\def \supp{{\rm supp}}
\def \Spin{{\rm Spin}}
\def \t{{{\mathfrak t}}}
\def \T{{\mathbb T}}
\def \tr{{\hspace{1pt}\rm tr\hspace{2pt}}}
\def \vol{{\rm vol}}
\def \z{\zeta}
\def \Z{\mathbb Z}
\def \={\ =\ }
\newcommand{\choice}[4]{\left\{
            \begin{array}{cl}#1&#2\\ #3&#4\end{array}\right.}
\newcommand{\rez}[1]{\frac{1}{#1}}
\newcommand{\der}[1]{\frac{\partial}{\partial #1}}
\newcommand{\norm}[1]{\parallel #1 \parallel}
\renewcommand{\matrix}[4]{\left( \begin{array}{cc}#1 & #2 \\ #3 & #4 \end{array}
            \right)}

\newtheorem{theorem}{Theorem}[section]
\newtheorem{conjecture}[theorem]{Conjecture}
\newtheorem{lemma}[theorem]{Lemma}
\newtheorem{corollary}[theorem]{Corollary}
\newtheorem{proposition}[theorem]{Proposition}

$$ $$
\newpage
\begin{center} {\bf Introduction} \end{center}

In 1956 A. Selberg introduced the zeta function
$$
Z(s) = \prod_c \prod_{N\geq 0} (1-e^{-(s+N)l(c)}),\ \ \ \ \Re(s)>>0,
$$
where the first product is taken over all primitive closed
geodesics in a compact Riemannian surface of genus $\geq
2$, equipped with the hyperbolic metric, and $l(c)$
denotes the length of the geodesic $c$. Selberg proved
that the product converges if the real part of $s$ is
large enough and that $Z(s)$ extends to an entire function
on $\C$, that it satisfies a functional equation as $s$ is
replaced by $1-s$, and that a generalized Riemann
hypothesis holds for the function $Z$.

The following list of references is only a small and
strongly subjective selection from the existing literature:
\cite{BuOl-buch, CaVo, Efrat-det, Efrat-dyn, Els, Hej,
Sarnak, Venkov, Voros}.

Generalizing Selberg's definition R. Gangolli \cite{Gang}
gave a zeta function with similar properties for all
locally symmetric spaces of rank one. These results were
later extended to vector bundles by M. Wakayama
\cite{Wak}. Both, Gangolli and Wakayama kept the spirit of
Selberg's original proof in that they explicitly inverted
the Selberg integral transform (which required more
knowledge on harmonic analysis then in Selberg's case,
however), then handicrafted a test function to be plugged
into the trace formula which would yield a higher
logarithmic derivative of the Selberg zeta function. In
\cite{Fr-at} D. Fried introduced heat kernel methods which
replaced the inversion of the transform conveniently.
Using supersymmetry arguments his approach was transferred
to some cases of higher rank by H. Moscovici and R.
Stanton \cite{MS-tors}.

In the present paper a new approach is given which allows
to extend the theory of the Selberg zeta function to
arbitrary groups. The central idea is to use a geometric
approach to shift down to a Levi subgroup of rank one less
and then plug in Euler-Poincar\'e functions with respect
to that Levi group. The themes of the paper are:
meromorphicity of the Selberg zeta function, the Patterson
conjecture relating the vanishing order of the zeta
function to group cohomology, and the evaluation of the
zeta function at special points to get topological
invariants.

The contents of the chapters are as follows: In the first
chapter notations are set and results from literature are
collected. In chapter 2 the existence of
pseudo-coefficients and Euler-Poincar\'e functions is
proved in a very general setting. Euler-Poincar\'e
functions are characterized by the fact that they give
Euler-Poincar\'e characteristics as traces under
irreducible representations. They are of great use in
various applications of the trace formula, in particular
in the theory of automorphic forms. Section 2 generalizes
existing approaches in that firstly the group need not be
connected and that secondly the representation to which it
is attached need not be spinorial. The trace formula is
presented in chapter 3. We restrict to the cocompact case
but we allow a broad class of test functions to enter the
trace formula.

In chapter 4 we present the Lefschetz formula which is
called so because it plays the same role for the Selberg
zeta function as the Lefschetz trace formula of \'etale
cohomology does for the Hasse-Weil zeta function of a
finite variety. In chapter 5 we extend the range of the
Lefschetz formula eventually getting the meromorphic
continuation of the Selberg zeta function.

The Patterson conjecture, which expresses the divisor of
the zeta function in terms of group cohomology is, in this
setting, proven in chapter 6. In chapter 7 we
generalize rank one results by expressing the holomorphic
torsion of a hermitian locally symmetric space by special
values of Selberg zeta functions. Finally, in chapter 8 we present the whole theory in
the $p$-adic setting where all analytic arguments become a lot easier.

Chapter one is merely a collection of facts from the
literature, but chapter 2 is entirely new and has not been
published elsewhere. In chapter 4 there is a new part in
section 4.4, where a new interpretation of the Lefschetz
formula is given which is coherent with the Deninger
philosophy. Section 5 is new since the theory of
Euler-Poincar\'e functions of chapter 2 allows for the
meromorphic continuation of a much wider class of zeta
functions than was possible hitherto. Chapter 7
is entirely new and so is the Lefschetz formula in chapter 8.

\pagebreak[4] {\it Acknowledgements:} I thank the RIMS at Kyoto,
the Waseda University at Tokyo and the Princeton University for
their hospitality. I thank Ulrich Bunke, Christopher Deninger,
Werner Hoffmann, Georg Illies, Andreas Juhl, Martin Olbrich and
Peter Sarnak for their comments and other help.

\chapter{Prerequisites}

In this chapter we collect some facts from literature which
will be needed in the sequel. Proofs will only be given by
sketches or references.

\section{Notations} \label{notations}

We denote Lie groups by upper case roman letters $G,H,K,\dots$ and
the corresponding real Lie algebras by lower case German letters
with index $0$, that is: $\g_0,\h_0,\k_0,\dots$. The complexified
Lie algebras will be denoted by $\g,\h,\k,\dots$, so, for
example: $\g =\g_0\otimes_\R \C$.

For a Lie group $L$ with Lie algebra $\l_0$ let $\Ad : L\ra
GL(\l_0)$ be the adjoint representation (\cite{bourb} III.3.12).
By definition, $\Ad(g)$ is the differential of the map $x\mapsto
gxg^{-1}$ at the point $x=e$. Then $\Ad(g)$ is a Lie algebra
automorphism of $\g_0$. A Lie group $L$ is said to be {\it of
inner type}\index{of inner type} if $\Ad(L)$ lies in the complex
adjoint group of the Lie algebra $\l$.

A real Lie group $G$ is said to be a {\it real reductive
group}\index{real reductive group} if there is a linear algebraic
group $\CG$ defined over $\R$ which is reductive as an algebraic
group \cite{Borel-lingroups} and a morphism $\alpha : G\ra \CG(\R)$
with finite kernel and cokernel (Since we do not insist that the
image of $\alpha$ is normal the latter condition means that ${\rm
im}(\alpha)$ has finite index in $\CG(\R)$). This implies in
particular that $G$ has only finitely many connected components
(\cite{Borel-lingroups} 24.6.c).

Note that any real reductive group $G$ of inner type is of {\it
Harish-Chandra class},\index{Harish-Chandra class} i.e., $G$ is of
inner type, the Lie algebra $\g_0$ of $G$ is reductive, $G$ has
finitely many connected components and the connected subgroup
$G_{der}^0$\index{$G_{der}^0$} corresponding to the Lie subalgebra
$[\g_0,\g_0]$ has finite center.

The following are of importance:
\begin{itemize}
\item
a connected semisimple Lie group with finite center is a real
reductive group of inner type (\cite{wall-rg1} 2.1.3),
\item
if $G$ is a real reductive group of inner type and $P=MAN$ is the
Langlands decomposition (\cite{wall-rg1} 2.2.7) of a parabolic
subgroup then the groups $M$ and $AM$ are real reductive of inner
type (\cite{wall-rg1} 2.2.8)
\end{itemize}

The usual terminology of algebraic groups carries over to real
reductive groups, for example a {\it torus}\index{torus} (or a
Cartan subgroup) in $G$ is the inverse image of (the real points
of) a torus or a Cartan subgroup in $\CG(\R)$. An element of $G$
will be called {\it semisimple}\index{semisimple element} if it
lies in some torus of $G$. The {\it split component} of $G$ is the
identity component of the greatest split torus in the center of
$G$. Note that for a real reductive group the Cartan subgroups are
precisely the maximal tori.\index{Cartan subgroup}

Let $G$ be a real reductive group then there exists a {\it Cartan
involution}\index{Cartan involution} i.e., an automorphism
$\theta$\index{$\theta$} of $G$ satisfying $\theta^2=Id$ whose
fixed point set is a maximal compact subgroup $K$ and which is the
inverse ($a\mapsto a^{-1}$) on the split component of $G$. All
Cartan involutions are conjugate under automorphisms of $G$.

Fix a Cartan involution $\theta$ with fixed point set the maximal
compact subgroup $K$ and let $\k_0$ be the Lie algebra of $K$.
\index{$\k_0$} The group $K$ acts on $\g_0$ via the adjoint
representation and there is a $K$-stable decomposition
$\g_0=\k_0\oplus\p_0$, \index{$\p_0$} where $\p_0$ is the
eigenspace of (the differential of) $\theta$ to the eigenvalue
$-1$. Write $\g=\k\oplus\p$ for the complexification. This is
called the {\it Cartan decomposition}.\index{Cartan decomposition}

\begin{lemma}\index{$B(X,Y)$} There is a symmetric bilinear form $B : \g_0\times\g_0\ra \R$
such that
\begin{itemize}
\item
$B$ is invariant, that is $B(\Ad(g)X,\Ad(g)Y)=B(X,Y)$ for all
$g\in G$ and all $X,Y\in\g_0$ and
\item
$B$ is negative definite on $\k_0$ and positive definite on its
orthocomplement $\p_0=\k_0^\perp\subset\g_0$.
\end{itemize}
\end{lemma}

\prf For $X\in\g_0$ let $\ad(X):\g_0\ra\g_0$ be the adjoint
defined by $\ad(X)Y=[X,Y]$. Since for $g\in G$ the map $\Ad(g)$ is
a Lie algebra homomorphism we infer that
$\ad(\Ad(g)X)=\Ad(g)\ad(X)\Ad(g)^{-1}$ and therefore the {\it
Killing form}\index{Killing form}
$$
B_K(X,Y)\= \tr (\ad(X)\ad(Y))
$$
is invariant. It is known that if $B_K$ is nondegenerate, i.e.,
$\g_0$ is semisimple, then $B=B_K$ satisfies the claims of the
lemma. In the general case we have
$\g_0=\a_0\oplus\c_0\oplus\g_0'$, where $\a_0\oplus\c_0$ is the
center of $\g_0$ and $\g_0'$ its derived algebra, which is
semisimple and so $B_K|_{\g_0'}$ is nondegenerate, whereas
$B_K|_{\a_0\oplus\c_0}=0$. Further $\c_0$ is the eigenspace of
$\theta$ in the center of $\g_0$ corresponding to the eigenvalue
$1$ whereas $\a_0$ is the eigenspace of $-1$. For $g\in G$ the
adjoint $\Ad(g)$ is easily seen to preserve $\a_0$ and $\c_0$, so
we get a representation $\rho:G\ra GL(\a_0)\times GL(\c_0)$. This
representation is trivial on the connected component $G^0$ of $G$
hence it factors over the finite group $G/G^0$. Therefore there is
a positive definite symmetric bilinear form $B_a$ on $\a_0$ which
is invariant under $\Ad_c$ and similarly a negative definite
symmetric bilinear form $B_c$ on $\c_0$ which is invariant. Let
$$
B\= B_a\oplus B_c\oplus B_K|_{\g_0'}.
$$
Then $B$ satisfies the claims of the lemma. \qed

Let $U(\g_0)$\index{$U(\g_0)$} be the universal enveloping algebra
of $\g_0$. It can be constructed as the quotient of the tensorial
algebra
$$
T(\g_0) \= \R\oplus\g_0\oplus (\g_0\otimes\g_0)\oplus\dots
$$
by the two-sided ideal generated by all elements of the form
$X\otimes Y-Y\otimes X-[X,Y]$, where $X,Y\in\g_0$.

The algebra $U(\g_0)$ can be identified with the $\R$-algebra of
all left invariant differential operators on $G$. Let $\g
=\g_0\otimes\C$\index{$\g$} be the complexification of $\g_0$ and
$U(\g)$ be its universal enveloping algebra which is the same as
the complexification of $U(\g_0)$. Then $\g$ is a subspace of
$U(\g)$ which generates $U(\g)$ as an algebra and any Lie algebra
representation of $\g$ extends uniquely to a representation of
the associative algebra $U(\g)$.

The form $B$ gives an identification of $\g_0$ with its dual space
$\g_0^*$. On the other hand $B$ defines an element in
$\g_0^*\otimes \g_0^*$. Thus we get a natural element in
$\g_0\otimes\g_0\subset T(\g_0)$. The image $C$ of this element
in $U(\g_0)$ is called the {\it Casimir operator} attached to $B$.
\index{Casimir operator} It is a differential operator of order
two and it lies in the center of $U(\g_0)$. In an more concrete
way the Casimir operator can be described as follows: Let $X_1,
\dots, X_m$ be a basis of $\g_0$ and let $Y_1, \dots, Y_m$ be the
dual basis with respect to the form $B$  then the Casimir operator
is given by
$$
C\= X_1 Y_1 +\dots +X_m Y_m.
$$
We have

\begin{lemma}\label{C-central}
For any $g\in G$ the Casimir operator is invariant under
$\Ad(g)$, that is $\Ad(g) C=C$.
\end{lemma}

\prf This follows from the invariance of $B$
\qed

Let $X$ denote the quotient manifold $G/K$. The tangent space at
$eK$ identifies with $\p_0$ and the form $B$ gives a
$K$-invariant positive definite inner product on this space.
Translating this by elements of $G$ defines a $G$-invariant
Riemannian metric on $X$. This makes $X$ the most general
globally symmetric space of the noncompact type \cite{helg}.

Let $\hat{G}$\index{$\hat{G}$} denote the {\it unitary
dual}\index{unitary dual} of $G$, i.e., $\hat{G}$ is the set of
isomorphism classes of irreducible unitary representations of $G$.

Let $(\pi,V_\pi)$ be a continuous representation of $G$ on some
Banach space $V_\pi$. The subspace $V_\pi^\infty$ of {\it smooth
vectors}\index{smooth vectors} is defined to be the subspace of
$V_\pi$ consisting of all $v\in V_\pi$ such that the map
$g\mapsto \pi(g)v$ is smooth. The universal enveloping algebra
$U(\g)$ operates on $V_\pi^\infty$ via
$$
\pi(X): v\mapsto X_g(\pi(g)v)\ |_{g=e}
$$
for $X$ in $\g$.

A {\it $(\g,K)$-module}\index{$(\g,K)$-module} is by definition a
complex vector space $V$ which is a $K$-module such that for each
$v\in V$ the space spanned by the orbit $K.v$ is finite
dimensional. Further $V$ is supposed to be a $\g$-module and the
following compatibility conditions should be satisfied:
\begin{itemize}
\item
for $Y\in\k\subset\g$ and $v\in V$ it holds
$$
Y.v\= \left.\frac{d}{dt}\right|_{t=0} \exp(tY).v,
$$
\item
for $k\in K$, $X\in\g$ and $v\in V$ we have
$$
k.X.v\= \Ad(k)X.k.v.
$$
\end{itemize}
A $(\g,K)$-module $V$ is called {\it irreducible} or simple if it
has no proper submodules and it is called of finite length if
there is a finite filtration
$$
0=V_0\subset\dots\subset V_n=V
$$
of submodules such that each quotient $V_j/V_{j-1}$ is
irreducible. Further $V$ is called {\it admissible} if for each
$\tau\in\hat{K}$ the space $\Hom_K(\tau,V)$ is finite dimensional.
An admissible $(\g,K)$-module of finite length is called a {\it
Harish-Chandra module}\index{Harish-Chandra module}.

Now let $(\pi,V_\pi)$ again be a Banach representation of $G$. Let
$V_{\pi,K}$ be the subspace of $V_\pi$ consisting of all vectors
$v\in V_\pi$ such that the $K$-orbit $\pi(K)v$ spans a finite
dimensional space. Then $V_{\pi,K}$  is called the space of {\it
$K$-finite vectors} in $V_\pi$. \index{$K$-finite vectors} The
space $V_{\pi,K}$ is no longer a $G$-module but remains a
$K$-module. Further the space $V_{\pi,K}^\infty = V_{\pi,K}\cap
V_\pi^\infty$ is dense in $V_{\pi,K}$ and is stable under $K$, so
$V_{\pi,K}^\infty$ is a $(\g,K)$-module. By abuse of notation we
will often write $\pi$ instead of $V_\pi$ and $\pi_K^\infty$
instead of $V_{\pi,K}^\infty$. The representation $\pi$ is called
admissible if $\pi_K^\infty$ is. In that case we have
$\pi_K^\infty =\pi_K$ since a dense subspace of a finite
dimensional space equals the entire space. So then $\pi_K$ is a
$(\g,K)$-module.

\begin{lemma}
Let $(\pi,V_\pi)$ be an irreducible admissible representation of
$G$ then the Casimir operator $C$ acts on $V_\pi^\infty$ by a
scalar denoted $\pi(C)$.
\end{lemma}

\prf By the formula $\pi(g)\pi(C)\pi(g)^{-1}=\pi(\Ad(g)C)$ and
Lemma \ref{C-central} we infer that $\pi(C)$ commutes with
$\pi(g)$ for every $g\in G$. Therefore the claim follows from the
Lemma of Schur (\cite{wall-rg1} 3.3.2).  \qed

Any $f\in L^1(G)$ will define a continuous operator
\index{$\pi(f)$}
$$
\pi(f)\=\int_Gf(x)\pi(x) dx
$$
on $V_\pi$.

Let $N$ be a natural number and let
$L_{2N}^1(G)$\index{$L_{2N}^1(G)$} be the set of all $f\in
C^{2N}(G)$ which satisfy $Df\in L^1(G)$ for any $D\in U(\g)$ with
$\deg(D)\le 2N$.

\begin{lemma}
Let $N$ be an integer $>\frac{\dim G}{2}$. Let $f\in L_{2N}^1(G)$
then for any irreducible unitary representation $\pi$ of $G$ the
operator $\pi(f)$ will be of trace class.
\end{lemma}

\prf Let $C$ denote the Casimir operator of $G$ and let $C_K$ be
the Casimir operator of $K$. Let $\lap = -C+2C_K\in U(\g)$ the
group Laplacian. It is known that for some $a>0$ the operator
$\pi(\lap +a)$ is positive and $\pi(\lap +a)^{-N}$ is of trace
class. Let $g=(\lap +a)^Nf$ then $g\in L^1(G)$, so $\pi(g)$ is
defined and gives a continuous linear operator on $V_\pi$. We
infer that $\pi(f)=\pi(\lap +a)^{-N}\pi(g)$ is of trace class.
\qed

Finally we need some more notation. The form $\langle X,Y\rangle
=-B(X,\theta(Y))$ is positive definite on $\g_0$ and therefore
induces a positive definite left invariant top differential form
$\omega_L$ on any closed subgroup $L$ of $G$. If $L$ is compact we
set
$$
v(L)\= \int_L\omega_L.
$$
Let $H=AB$ be a $\theta$-stable Cartan subgroup where $A$ is the
connected split component of $H$ and $B$ is compact. The double
use of the letter $B$ here will not cause any confusion. Then
$B\subset K$. Let $\Phi$ denote the root system of $(\g,\h)$,
where $\g$ and $\h$ are the complexified Lie algebras of $G$ and
$H$. Let $\g =\h\oplus\bigoplus_{\alpha\in\Phi}\g_\alpha$ be the
root space decomposition. Let $x\ra x^c$ denote the complex
conjugation on $\g$ with respect to the real form $\g_0= Lie(G)$.
A root $\alpha$ is called {\it imaginary}\index{imaginary root} if
$\alpha^c=-\alpha$ and it is called {\it real}\index{real root} if
$\alpha^c=\alpha$. Every root space $\g_\alpha$ is one dimensional
and has a generator $X_\alpha$ satisfying:
$$
[X_\alpha,X_{-\alpha}]=Y_\alpha\ \ \ \ {\rm with}\ \ \ \
\alpha(.)=B(Y_\alpha,.)
$$ $$
B(X_\alpha,X_{-\alpha})=1
$$
and $X_\alpha^c=X_{\alpha^c}$ if $\alpha$ is non-imaginary and
$X_\alpha^c=\pm X_{-\alpha}$ if $\alpha$ is imaginary. An
imaginary root $\alpha$ is called {\it compact}\index{compact
root} if $X_\alpha^c=-X_{-\alpha}$ and {\it noncompact} otherwise.
Let $\Phi_n$ be the set of noncompact imaginary roots and choose a
set $\Phi^+$ of positive roots such that for $\alpha\in\Phi^+$
nonimaginary we have that $\alpha^c\in\Phi^+$. Let $W=W(G,H)$ be
the {\it Weyl group} of $(G,H)$\index{Weyl group of $(G,H)$}, that
is
$$
W\= \frac{{\rm normalizer}(H)}{{\rm centralizer}(H)}.
$$
Let ${\rm rk}_\R(G)$ be the dimension of a maximal $\R$-split
torus in $G$ and let $\nu =\dim G/K -{\rm rk}_\R(G)$. We define
the {\it Harish-Chandra constant} of $G$ by
$$
c_G\= (-1)^{|\Phi_n^+|} (2\pi)^{|\Phi^+|}
2^{\nu/2}\frac{v(T)}{v(K)}|W|.
$$

\section{Normalization of Haar measures}
Although the results will not depend on normalizations we will
need to normalize Haar measures for the computations along the
way. First for any compact subgroup $C\subset G$ we normalize its
Haar measure so that it has total mass one, i.e., $\vol(C)=1$.
Next let $H\subset G$ be a reductive subgroup, and let $\theta_H$
be a Cartan involution on $H$ with fixed point set $K_H$. The
same way as for $G$ itself the form $B$ restricted to the Lie
algebra of $H$ induces a Riemannian metric on the manifold
$X_H=H/K_H$. Let $dx$ denote the volume element of that metric.
We get a Haar measure on $H$ by defining
$$
\int_Hf(h) dh \= \int_{X_H}\int_{K_H} f(xk) dk dx
$$
for any continuous function of compact support $f$ on $H$.

Let $P\subset G$ be a {\it parabolic subgroup}\index{parabolic
subgroup} of $G$ (see \cite{wall-rg1} 2.2). Let $P=MAN$ be the
{\it Langlands decomposition} \index{Langlands decomposition} of
$P$. Then $M$ and $A$ are reductive, so there Haar measures can be
normalized as above. Since $G=PK=MANK$ there is a unique Haar
measure $dn$ on the unipotent radical $N$ such that for any
constant function $f$ of compact support on $G$ it holds:
$$
\int_G f(x)dx \= \int_M\int_A\int_N\int_K f(mank) dkdndadm.
$$
Note that these normalizations coincide for Levi subgroups with
the ones met by Harish-Chandra in (\cite{HC-HA1} sect. 7).

\section{Invariant distributions}
In this chapter we shall throughout assume that $G$ is a real
reductive group of inner type. A distribution $T$ on $G$, i.e., a
continuous linear functional $T: C_c^\infty(G)\ra \C$ is called
{\it invariant}\index{invariant distribution} if for any $f\in
C_c^\infty(G)$ and any $y\in G$ it holds: $T(f^y)=T(f)$, where
$f^y(x)=f(yxy^{-1})$. Examples are:
\begin{itemize}
\item orbital integrals: $f\mapsto \CO_g(f)=\int_{G_g\bs G}f(x^{-1}gx)dx$ and
\item traces: $f\mapsto \tr\pi(f)$ for $\pi\in\hat{G}$.
\end{itemize}
These two examples can each be expressed in terms of the other.
Firstly, Harish-Chandra proved that for any $\pi\in\hat{G}$ there
exists a conjugation invariant locally integrable function
$\Theta_\pi$ on $G$ such that for any $f\in C_c^\infty(G)$
$$
\tr\pi(f)\=\int_Gf(x)\Theta_\pi(x)  dx.
$$
Recall the {\it Weyl integration formula} which says that for any
integrable function $\ph$ on $G$ we have
$$
\int_G\ph(x)dx\=
\sum_{j=1}^r\rez{|W(G,H_j)|}\int_{H_j}\int_{G/H_j}\ph(xhx^{-1})|\det(1-h|\g
/\h_j)|dx dh,
$$
where $H_1,\dots,H_r$ is a maximal set of nonconjugate Cartan
subgroups in $G$ and for each Cartan subgroup $H$ we let
$W(G,H)$\index{$W(G,H)$} denote its Weyl group, i.e., the quotient
of the normalizer of $H$ in $G$ by its centralizer.

An element $x$ of $G$ is called {\it regular}\index{regular} if
its centralizer is a Cartan subgroup. The set of regular elements
$G^{reg}$ is open and dense in $G$ therefore the integral above
can be taken over $G^{reg}$ only. Letting $H_j^{reg}:=H_j\cap
G^{reg}$ we get

\begin{proposition}
Let $N$ be a natural number bigger than $\frac{\dim G}{2}$, then
for any $f\in L_{2N}^1(G)$ and any $\pi\in\hat{G}$ we have
$$
\tr\pi(f)\=\sum_{j=1}^r\rez{|W(G,H_j)|}\int_{H_j^{reg}}
\CO_h(f)\Theta_\pi(h)|\det(1-h|\g /\h_j)| dh.
$$
\end{proposition}

That is, we have expressed the trace distribution in terms of
orbital integrals. In the other direction it is also possible to
express semisimple orbital integrals in terms of traces.

At first let $H$ be a $\theta$-stable Cartan subgroup of $G$. Let
$\h$ be its complex Lie algebra and let $\Phi =\Phi(\g,\h)$ be the
set of roots. Let $x\ra x^c$ denote the complex conjugation on
$\g$ with respect to the real form $\g_0= Lie_\R(G)$. Choose an
ordering $\Phi^+\subset \Phi$ and let $\Phi^+_I$ be the set of
positive imaginary roots. To any root $\alpha\in\Phi$ let
\begin{eqnarray*}
H &\ra & \C^\times\\
 h &\mapsto & h^\alpha
\end{eqnarray*}\index{$h^\alpha$}
be its character, that is, for $X\in\g_\alpha$ the root space to
$\alpha$ and any $h\in H$ we have $\Ad(h)X=h^\alpha X$. Now put
$$
\ '\lap_I(h)\= \prod_{\alpha\in\Phi_I^+}(1-h^{-\alpha}).
$$\index{$\ '\lap_I(h)$}
Let $H=AT$ where $A$ is the connected split component and $T$ is
compact. An element $at\in AT=H$ is called {\it split
regular}\index{split regular} if the centralizer of $a$ in $G$
equals the centralizer of $A$ in $G$. The split regular elements
form a dense open subset containing the regular elements of $H$.
Choose a parabolic $P$ with split component $A$, so $P$ has
Langlands decomposition $P=MAN$. For $at\in AT =H$ let
\begin{eqnarray*}
\lap_+(at) &=& \left| \det((1-\Ad((at)^{-1}))|_{\g /\a\oplus
\m})\right|^{\rez{2}}\\ {}\\
 & =& \left|\det((1-\Ad((at)^{-1}))|_\n)\right| a^{\rho_P}\\ {}\\
 &=&
 \left|\prod_{\alpha\in\Phi^+-\Phi_I^+}(1-(at)^{-\alpha})\right|
 a^{\rho_P},
\end{eqnarray*}\index{$\lap_+(at)$}
where $\rho_P$\index{$\rho_P$} is the half of the sum of the roots
in $\Phi(P,A)$, i.e., $a^{2\rho_P}=\det(a|\n)$. We will also write
$h^{\rho_P}$ instead of $a^{\rho_P}$.

For any $h\in H^{reg} = H\cap G^{reg}$ let
$$
'F_f^H(h)\= 'F_f(h)\= \ '\lap_I(h) \lap_+(h) \int_{G/A}
f(xhx^{-1}) dx.
$$
It then follows directly from the definitions that for $h\in
H^{reg}$ it holds
$$
\CO_h(f) \= \frac{'F_f(h)}
                 {h^{\rho_P}\det(1-h^{-1}|(\g/\h)^+)},
$$
where $(\g/\h)^+$ is the sum of the root spaces attached to
positive roots. There is an extension of this identity to
nonregular elements as follows: For $h\in H$ let $G_h$ denote its
centralizer in $G$. Let $\Phi^+(\g_h,\h)$ be the set of positive
roots of $(\g_h,\h)$. Let
$$
\varpi_h \= \prod_{\alpha\in\Phi^+(\g_h,\h)}Y_\alpha,
$$
then $\varpi_h$ defines a left invariant differential operator on
$G$.

\begin{lemma}
For any $f\in L_{2N}^1(G)$ and any $h\in H$ we have
$$
\CO_h(f) \= \frac{\varpi_h 'F_f(h)}
                 {c_h h^{\rho_P}\det(1-h^{-1}|(\g/\g_h)^+)}.
$$
\end{lemma}

\prf This is proven in section 17 of \cite{HC-DS}. \qed

Our aim is to express orbital integrals in terms of traces of
representations. By the above lemma it is enough to express
$'F_f(h)$ it terms of traces of $f$ when $h\in H^{reg}$. For this
let $H_1=A_1T_1$ be another $\theta$-stable Cartan subgroup of $G$
and let $P_1=M_1A_1N_1$ be a parabolic with split component $A_1$.
Let $K_1=K\cap M_1$. Since $G$ is connected the compact group
$T_1$ is an abelian torus and its unitary dual $\widehat{T_1}$ is
a lattice. The Weyl group $W=W(M_1,T_1)$ acts on $\widehat{T_1}$
and $\widehat{t_1}\in\widehat{T_1}$ is called {\it regular} if its
stabilizer $W(\widehat{t_1})$ in $W$ is trivial. The regular set
$\widehat{T_1}^{reg}$ modulo the action of $W(K_1,T_1)\subset
W(M_1,T_1)$ parameterizes the discrete series representations of
$M_1$ (see \cite{Knapp}). For $\widehat{t_1}\in\widehat{T_1}$
Harish-Chandra \cite{HC-S} defined a distribution
$\Theta_{\widehat{t_1}}$ on $G$ which happens to be the trace of
the discrete series representation $\pi_{\widehat{t_1}}$ attached
to $\widehat{t_1}$ when $\widehat{t_1}$ is regular. When
$\widehat{t_1}$ is not regular the distribution
$\Theta_{\widehat{t_1}}$ can be expressed as a linear combination
of traces as follows. Choose an ordering of the roots of
$(M_1,T_1)$ and let $\Omega$ be the product of all positive roots.
For any $w\in W$ we have $w\Omega = \epsilon(w)\Omega$ for a
homomorphism $\epsilon : W\ra \{ \pm 1\}$. For nonregular
$\widehat{t_1}\in\widehat{T_1}$ we get
$\Theta_{\widehat{t_1}}=\rez{|W(\widehat{t_1})|}\sum_{w\in
W(\widehat{t_1})}\epsilon(w)\Theta'_{w,\widehat{t_1}}$, where
$\Theta'_{w,\widehat{t_1}}$ is the character of an irreducible
representation $\pi_{w,\widehat{t_1}}$ called a limit of discrete
series representation. We will write $\pi_{\widehat{t_1}}$ for the
virtual representation $\rez{|W(\widehat{t_1})|}\sum_{w\in
W_{\widehat{t_1}}}\epsilon(w)\pi_{w,\widehat{t_1}}$.

Let $\nu :a\mapsto a^\nu$ be a unitary character of $A_1$ then
$\widehat{h_1}=(\nu,\widehat{t_1})$ is a character of
$H_1=A_1T_1$. Let $\Theta_{\widehat{h_1}}$ be the character of the
representation $\pi_{\widehat{h_1}}$ induced parabolically from
$(\nu,\pi_{\widehat{t_1}})$. Harish-Chandra has proven

\begin{theorem}\label{inv-orb-int}
Let $H_1,\dots,H_r$ be maximal a set of nonconjugate
$\theta$-stable Cartan subgroups. Let $H=H_j$ for some $j$. Then
for each $j$ there exists a continuous function $\Phi_{H|H_j}$ on
$H^{reg}\times \hat{H_j}$ such that for $h\in H^{reg}$ it holds
$$
'F_f^H(h)\= \sum_{j=1}^r
\int_{\hat{H_j}}\Phi_{H|H_j}(h,\widehat{h_j})\
\tr\pi_{\widehat{h_j}}(f)\ d\widehat{h_j}.
$$
Further $\Phi_{H|H_j}=0$ unless there is $g\in G$ such that
$gAg^{-1}\subset A_1$. Finally for $H_j=H$ the function can be
given explicitly as
\begin{eqnarray*}
\Phi_{H|H}(h,\hat{h}) &=& \rez{|W(G,H)|}\sum_{w\in
W(G,H)}\epsilon(w|T)\langle w\hat{h},h\rangle\\
 &=& \rez{|W(G,H)|} \ '\lap(h)\Theta_{\hat{h}}(h),
\end{eqnarray*}
where $\ '\lap = \lap_+ \lap_I$.
\end{theorem}

\prf \cite{HC-S}. \qed

\section{Regularized determinants and products}

Frequently we will use the notion of a
regularized determinant. For this let $A$ be a positive
operator, densely defined on some Hilbert space. We will
call the operator $A$ \emph{zeta-admissible}, if
\begin{itemize}
\item[-] there is a $p\geq 1$ such that $A^{-p}$ is of trace class and
\item[-] the trace of the heat operator $e^{-tA}$ admits an asymptotic
expansion
$$
\theta_A(t) \= {\rm tr}\ e^{-tA}\ \sim\ \sum_{k=0}^\infty c_{k}
t^{\alpha_k}
$$
as $t\rightarrow 0$ where $\alpha_k \in \R$, $\alpha_k \rightarrow
+\infty$.
\end{itemize}

The reader should notice that, from the fact that $A^{-p}$ is of
trace class it immediately follows that the Hilbert space $H$ is
separable and has a basis of $A$-eigenvectors, each eigenvalue
occurring with a finite multiplicity and the eigenvalues
accumulating at most at infinity.

Under these circumstances one defines the \emph{zeta
function} of $A$ as: $$ \zeta_A (s) \= {\rm tr}\ A^{-s}\ \
\ {\rm for}\ {\rm Re}\ s > p. $$ Then for ${\rm Re} s >>
0$ and $\lambda >> 0$ we consider the Mellin transform of
$\theta_A$: $$ M(z,\lambda )\= \int_0^\infty t^{z-1}
e^{-\lambda t} \theta_A(t)\ dt. $$ The asymptotic
expansion shows that for fixed $\lambda >0$ the function
$z\mapsto M(z,\lambda )$ with $M(z,\lambda ) = \Ga (z)
\zeta_{A+\la}(z)$ is holomorphic in $\C$ up to simple
poles at $-\alpha_k -n$, $k,n\geq 0$ of residue $$ {\rm
res}_{z=-\alpha_k-n}\ M(z,\lambda ) \= \la^nc_k (-1)^n/n!.
$$ (When two of those coincide the residues add up.) This
shows that $\zeta_{A+\la}(z)$ is regular at z=0. Extending
the case of finite dimension we define: $$ {\rm
det}(A+\lambda ) \= {\rm exp}(-\zeta_{A+\la}'(0)). $$
Using the asymptotic expansion one sees that the function
$\lambda \mapsto {\rm det}(A+\la )$ extends to an entire
function with zeroes given by the eigenvalues of -A, the
multiplicity of a zero at $\la_0$ being the multiplicity
of the eigenvalue $\la_0$.

We will extend the definition of the determinant ${\rm det}(A)$ to
the situation where not necessarily $A$ but  $A'=A\mid_{({\rm
ker}\ A)^\perp}$ is zeta admissible, then we will define $\zeta_A
= \zeta_{A'}$ and $\det'(A) = \det(A')$.

Evidently the determinant $\det(A+\la)$ only depends on the
sequence of eigenvalues $0<a_1\le a_2\le a_2\le\dots$ of $A$.
Further for any sequence of complex numbers there exists an
operator having the given eigenvalues. So we call the sequence
$(a_j)$ zeta admissible if the operator $A$ is and define the {\it
regularized product}\index{regularized product} as
$$
\widehat{\prod_j}(a_j+\la)\= \det(A+\la).
$$

\noindent \emph{Examples.} The theory of the heat equation shows
that an elliptic differential operator with positive definite
principal symbol is zeta admissible \cite{gilk}. Another class of
examples often occurs in number theory: Let $p$ be a nonconstant
polynomial over the reals with postive leading coefficient.
Assume that $a_n=p(n)$ for all but finitely many $n$, then the
operator defined by the sequence of eigenvalues $a_n$ is zeta
admissible.

\chapter{Euler-Poincar\'e functions}

\label{ep} In this chapter we generalize the construction
of pseudo-coefficients \cite{Lab} and Euler-Poincar\'e
functions to non-connected groups. Here $G$ will be a real
reductive group. It will be assumed that $G$ admits a
compact Cartan subgroup. It then follows that $G$ has
compact center.

\section{Existence}
Fix a maximal compact subgroup $K$ of $G$ and a Cartan $T$ of $G$
which lies inside  $K$. The group $G$ is called {\it orientation
preserving} \index{orientation preserving} if $G$ acts by
orientation preserving diffeomorphisms on the manifold $X=G/K$.
For example, the group $G=SL_2(\R)$ is orientation preserving but
the group $PGL_2(\R)$ is not. Recall the Cartan decomposition
$\g_0=\k_0\oplus\p_0$. Note that $G$ is orientation preserving if
and only if its maximal compact subgroup $K$ preserves
orientations on $\p_0$.

\begin{lemma} \label{orient}
The following holds:
\begin{itemize}
\item
Any connected group is orientation preserving.
\item
If $X$ carries the structure of a complex manifold which is left
stable by $G$, then $G$ is orientation preserving.
\end{itemize}
\end{lemma}

\prf The first is clear. The second follows from the fact that
biholomorphic maps are orientation preserving.
 \qed

Let $\t$\index{$\t$} be the complexified Lie algebra of the Cartan
subgroup $T$. We choose an ordering of the roots $\Phi(\g
,\t)$\index{$\Phi(\g ,\t)$} of the pair $(\g ,\t)$
\cite{wall-rg1}. This choice induces a decomposition $\p = \p_-
\oplus \p_+$,\index{$\p_\pm$} where $\p_\pm$ is the sum of the
positive/negative root spaces which lie in $\p$. As usual denote
by $\rho$\index{$\rho$} the half sum of the positive roots. The
chosen ordering induces an ordering of the {\it compact
roots}\index{compact root} $\Phi(\k ,\t)$ which form a subset of
the set of all roots $\Phi(\g,\t)$. Let $\rho_K$ denote the half
sum of the positive compact roots. Recall that a function $f$ on
$G$ is called {\it $K$-central} \index{$K$-central} if
$f(kxk^{-1})=f(x)$ for all $x\in G$, $k\in K$. For any
$K$-representation $(\rho,V)$ let $V^K$ denote the space of
$K$-fixed vectors, i.e., \index{$V^K$}
$$
V^K \=\{ v\in V | \rho(k)v=v\ \forall k\in K\}.
$$\index{$V^K$}
Let $(\tau,V_\tau)$ be a representation of $K$ on a finite
dimensional complex vector space $V_\tau$. Let $V_\tau^*$ be the
dual space then there is a representation $\breve{\tau}$ on
$V_{\breve{\tau}}:=V_\tau^*$ given by
$$
\breve{\tau}(k)\alpha(v)\= \alpha(\tau(k^{-1})v),
$$
for $k\in K$, $\alpha\in V_\tau^*$ and $v\in V_\tau$. This
representation is called the {\it contragredient} or {\it dual}
representation.\index{contragredient representation} \index{dual
representation} The restriction from $G$ to $K$ gives a ring
homomorphism of the representation rings:
$$
res_K^G:\Rep(G)\ra\Rep(K).
$$

\begin{theorem}\label{exist-ep}
(Euler-Poincar\'e functions) Let $(\tau ,V_\tau)$ a finite
dimensional representation of $K$. If $G$ is orientation
preserving or $\tau$ lies in the image of $res^G_K$, then there
is a compactly supported smooth $K$-central function $f_\tau$ on
$G$ such that for every irreducible unitary representation $(\pi
,V_\pi)$ of $G$ we have
$$
\tr\ \pi (f_\tau) \= \sum_{p=0}^{\dim (\p)} (-1)^p \dim (V_\pi
\otimes \wedge^p\p \otimes V_{\breve{\tau}})^K.
$$
We call $f_\tau$ an {\it Euler-Poincar\'e function} for $\tau$.

If, moreover, $K$ leaves invariant the decomposition
$\p=\p_+\oplus\p_-$ then there is a compactly supported smooth
$K$-central function $g_\tau$ on $G$ such that for every
irreducible unitary representation $(\pi ,V_\pi)$ we have
$$
\tr\ \pi (g_\tau) \= \sum_{p=0}^{\dim (\p_-)} (-1)^p \dim (V_\pi
\otimes \wedge^p\p_- \otimes V_{\breve{\tau}})^K.
$$
\end{theorem}

{\it Remark:} If the representation $\tau$ lies in the image of
$res_K^G$ or the group $G$ is connected then the theorem is well
known, \cite{CloDel}, \cite{Lab}. We will however have to apply
the theorem also in situations where neither condition is
satisfied.

\vspace{10pt}

\prf  We will concentrate on the case when $G$ is
orientation preserving, since the other case already is
dealt with in the literature. Suffice to say that it can be
treated similarly. Without loss of generality assume $\tau$
irreducible. Suppose given such a function $f$ which
satisfies the claims of the theorem except that it is not
necessarily $K$-central, then the function
$$
x\mapsto \int_K f(kxk^{-1})dk
$$
will satisfy all claims of the theorem. Thus one only needs to
construct a function having the claimed traces.

If $G$ is orientation preserving the adjoint action gives a
homomorphism $K\ra SO(\p)$. If this homomorphism happens to lift
to the double cover $Spin(\p)$ \cite{LawMich} we let $\tilde{G}=G$
and $\tilde{K}=K$. In the other case we apply the

\begin{lemma}\label{double-cover} If the homomorphism $K\ra
SO(\p)$ does not factor over the spin group $Spin(\p)$ then there is a double
covering $\tilde{G}\ra G$ such that with $\tilde{K}$ denoting the
inverse image of $K$ the induced homomorphism $\tilde{K}\ra
SO(\p)$ factors over $Spin(\p)\ra SO(\p)$. Moreover the kernel of
the map $\tilde{G}\ra G$ lies in the center of $\tilde{G}$
\end{lemma}

\prf At first $\tilde{K}$ is given by the pullback diagram:
$$
\begin{array}{ccc}
\tilde{K} & \ra & Spin(\p)\\
\downarrow & {} & \downarrow \alpha\\
K & \begin{array}{c}\Ad\\ \ra\\ {}\end{array} & SO(\p)
\end{array}
$$
that is, $\tilde{K}$ is given as the set of all $(k,g)\in K\times
Spin(\p)$ such that $\Ad(k)=\alpha(g)$. Then $\tilde{K}$ is a
double cover of $K$.

Next we use the fact that $K$ is a retract of $G$ to show that the
covering $\tilde{K}\ra K$ lifts to $G$. Explicitly let
$P=\exp(\p_0)$ then the map $K\times P\ra G, (k,p)\mapsto kp$ is a
diffeomorphism \cite{wall-rg1}. Let $g\mapsto
(\underline{k}(g),\underline{p}(g))$ be its inverse map. We let
$\tilde{G}\= \tilde{K}\times P$ then the covering $\tilde{K}\ra K$
defines a double covering $\beta :\tilde{G}\ra G$. We have to
install a group structure on $\tilde{G}$ which makes $\beta$ a
homomorphism and reduces to the known one on $\tilde{K}$. Now let
$k,k'\in K$ and $p,p'\in P$ then by
$$
k'p'kp\= k'k\ k^{-1}p'kp
$$
it follows that there are unique maps $a_K : P\times P\ra K$ and
$a_P :P\times P\ra P$ such that
\begin{eqnarray*}
\underline{k}(k'p'kp)&=& k'k a_K(k^{-1}p'k,p)\\
\underline{p}(k'p'kp)&=& a_P(k^{-1}p'k,p).
\end{eqnarray*}
Since $P$ is simply connected the map $a_K$ lifts to a map
$\tilde{a}_K: P\times P\ra \tilde{K}$. Since $P$ is connected
there is exactly one such lifting with $\tilde{a}_K(1,1)=1$. Now
the map
\begin{eqnarray*}
(\tilde{K}\times P)\times(\tilde{K}\times P)&\ra& \tilde{K}\times P\\
(k',p'),(k,p) &\mapsto&
(kk'\tilde{a}_K(k^{-1}p'k,p),a_P(k^{-1}p'k,p))
\end{eqnarray*}
defines a multiplication on $\tilde{G}=\tilde{K}\times P$ with
the desired properties.

Finally $\ker(\beta)$ will automatically be central because it is
a normal subgroup of order two.
\qed

Let $S$ be the spin representation of $Spin(\p)$ (see
\cite{LawMich}, p.36). It splits as a direct sum of two distinct
irreducible representations
$$
S\= S^+\oplus S^-.
$$
We will not go into the theory of the spin representation, we
only need the following properties:
\begin{itemize}
\item
The virtual representation
$$
(S^+- S^-)\otimes (S^+- S^-)
$$
is isomorphic to the adjoint representation on $\wedge^{even}\p
-\wedge^{odd}\p$ (see \cite{LawMich}, p. 36).
\item
If $K$ leaves invariant the spaces $\p_-$ and $\p_+$, as is the
case when $X$ carries a holomorphic structure fixed by $G$, then
there is a one dimensional representation $\epsilon$ of
$\tilde{K}$ such that
$$
(S^+-S^-)\otimes\epsilon \ \cong\
\wedge^{even}\p_--\wedge^{odd}\p_-.
$$
\end{itemize}
The proof of this latter property will be given in section
\ref{appA}.

\begin{theorem} \label{existh}
(Pseudo-coefficients) Assume that the group $G$ is
orientation preserving. Then for any finite dimensional
representation $(\tau ,V_\tau)$ of $\tilde{K}$ there is a
compactly supported smooth function $h_\tau$ on
$\tilde{G}$ such that for every irreducible unitary
representation $(\pi ,V_\pi)$ of $\tilde{G}$,
$$
\tr\ \pi (h_\tau) = \dim (V_\pi \otimes S^+ \otimes
V_{\breve{\tau}})^{\tilde{K}} - \dim (V_\pi \otimes S^- \otimes
V_{\breve{\tau}})^{\tilde{K}}.
$$
\end{theorem}

The functions given in this theorem are also known as
pseudo-coefficients \cite{Lab}. This result generalizes
the one in \cite{Lab} in several ways. First, the group
$G$ needn't be connected and secondly the representation
$\tau$ needn't be spinorial. The proof of this theorem
relies on the following lemma.

\begin{lemma}
Let $(\pi ,V_\pi)$ be an irreducible unitary representation of
$\tilde{G}$ and assume
$$
\dim (V_\pi \otimes S^+ \otimes V_{\breve{\tau}})^{\tilde{K}} -
\dim (V_\pi \otimes S^- \otimes V_{\breve{\tau}})^{\tilde{K}}
\neq 0,
$$
then the Casimir eigenvalue satisfies $\pi (C) = \breve{\tau}
(C_K) - B(\rho)+B(\rho_K)$.
\end{lemma}

\prf Let the $\tilde{K}$-invariant operator
 $$
 d_\pm : V_\pi\otimes S^\pm \ra V_\pi\otimes S^\mp
 $$
be defined by
 $$
 d_\pm : v\otimes s \mapsto \sum_i\pi(X_i)v\otimes c(X_i)s,
 $$
where $(X_i)$ is an orthonormal base of $\p$. The formula of
Parthasarathy, \cite{AtSch}, p. 55 now says
 $$
 d_-d_+\= d_+d_- \= \pi\otimes s^\pm(C_K) -\pi(C)\otimes 1
 -B(\rho)+B(\rho_K).
 $$
Our assumption leads to $ker(d_+ d_-) \cap \pi \otimes S(\tau)
\neq 0$, and therefore $0=\tau(C_K)-\pi(C) -B(\rho) +B(\rho_K)$.
\qed

For the proof of Theorem \ref{existh} let $(\tau,V_\tau)$ a finite
dimensional irreducible unitary representation of $\tilde{K}$ and
write $E_\tau$ for the $\tilde{G}$-homogeneous vector bundle over
${X}=\tilde{G}/\tilde{K}$ defined by $\tau$. The space of smooth
sections $\Ga^\infty (E_\tau)$ may be written as $\Ga^\infty
(E_\tau) = (C^\infty (\tilde{G}) \otimes V_\tau)^{\tilde{K}}$,
where $\tilde{K}$ acts on $C^\infty (\tilde{G})$ by right
translations. The Casimir operator $C$ of $\tilde{G}$ acts on
this space and defines a second order differential operator
$C_\tau$ on $E_\tau$. On the space of $L^2$-sections
$L^2(X,E_\tau)=(L^2(\tilde{G}) \otimes V_\tau)^{\tilde{K}}$ this
operator is formally selfadjoint with domain, say, the compactly
supported smooth functions and extends to a selfadjoint operator.
Consider a Schwartz function $f$ on $\R$ such that the Fourier
transform $\hat{f}$ has compact support. By general results on
hyperbolic equations (\cite{Tayl}, chap IV) it follows that the
smoothing operator $f(C_\tau)$, defined by the spectral theorem
has finite propagation speed. Since $f(C_\tau)$ is
$\tilde{G}$-equivariant it follows that $f(C_\tau)$ can be
represented as a convolution operator $\ph \mapsto \ph *
\breve{f}_\tau$, for some $\breve{f}_\tau \in
(C_c^\infty(\tilde{G})\otimes \End (V_\tau))^{\tilde{K}\times
\tilde{K}}$ and with $\tilde{f}_\tau$ denoting the pointwise
trace of $\breve{f}_\tau$ we have for $\pi \in \hat{\tilde{G}}$:
$\tr\ \pi(\tilde{f}_\tau) = f(\pi(C))\dim(V_\pi \otimes
V_\tau)^{\tilde{K}}$, where $\pi(C)$ denotes the Casimir
eigenvalue on $\pi$. This construction extends to virtual
representations by linearity.

Choose $f$ such that $f(\breve{\tau} (C_K) -B(\rho)+B(\rho_K))=1$.
Such an $f$ clearly exists. Let $\ga$ be the virtual
representation of $\tilde{K}$ on the space
$$
V_\ga = (S^+-S^-)\otimes V_\tau,
$$
then set $h_\tau =\tilde{f}_\ga$. Theorem \ref{existh} follows.
\qed

To get the first part of Theorem \ref{exist-ep} from Theorem
\ref{existh} one replaces $\tau$ in the proposition by the virtual
representation on $(S^+-S^-)\otimes V_\tau$. Since
$(S^+-S^-)\otimes (S^+-S^-)$ is as $\tilde{K}$ module isomorphic
to $\wedge^*\p$ we get the desired function, say $j$ on the group
$\tilde{G}$. Now if $\tilde{G}\ne G$ let $z$ be the nontrivial
element in the kernel of the isogeny $\tilde{G}\ra G$, then the
function
$$
f(x) \= \rez{2}(j(x)+j(zx))
$$
factors over $G$ and satisfies the claim.

To get the second part of the theorem one proceeds similarly
replacing $\tau$ by $\epsilon\otimes \tau$.
\qed

\section{Clifford algebras and Spin groups}\label{appA}
This section is solely given to provide a proof of the properties
of the spin representation used in the last section. We will
therefore not strive for the utmost generality but plainly state
things in the form needed. For more details the reader is
referred to \cite{LawMich}.

Let $V$ be a finite dimensional complex vector space and let
$q:V\ra\C$ be a non-degenerate quadratic form. We use the same
letter for the symmetric bilinear form:
$$
q(x,y) \= \rez{2}(q(x+y)-q(x)-q(y)).
$$
Let $SO(q)\subset GL(V)$ be the special orthogonal group of $q$.
The {\it Clifford algebra} \index{Clifford algebra} $Cl(q)$ will
be the quotient of the tensorial algebra
$$
TV \= \C\oplus V\oplus (V\otimes V)\oplus\dots
$$
by the two-sided ideal generated by all elements of the form
$v\otimes v+q(v)$, where $v\in V$.

This ideal is not homogeneous with respect to the natural
$\Z$-grading of $TV$, but it is homogeneous with respect to the
induced $\Z/2\Z$-grading given by the even and odd degrees. Hence
the latter is inherited by $Cl(q)$:
$$
Cl(q)\= Cl^0(q)\oplus Cl^1(q).
$$
For any $v\in V$ we have in $Cl(q)$ that $v^2=-q(v)$ and therefore
$v$ is invertible in $Cl(q)$ if $q(v)\ne 0$. Let $Cl(q)^\times$
be the group of invertible elements in $Cl(q)$. The algebra
$Cl(q)$ has the following universal property: For any linear map
$\ph :V\ra A$ to a $\C$-algebra $A$ such that $\ph(v)^2=-q(v)$
for all $v\in V$ there is a unique algebra homomorphism $Cl(v)\ra
A$ extending $\ph$.

Let $Pin(q)$ be the subgroup of the group $Cl(q)^\times$
generated by all elements $v$ of $V$ with $q(v)=\pm 1$. Let the
{\it complex spin group} \index{complex spin group} be defined by
$$
Spin(q) \= Pin(q)\cap Cl^0(q),
$$
i.e., the subgroup of $Pin(q)$ of those elements which are
representable by an even number of factors of the form $v$ or
$v^{-1}$ with $v\in V$. Then $Spin(q)$ acts on $V$ by $x.v =
xvx^{-1}$ and this gives a double covering: $Spin(q)\ra SO(q)$.

Assume the dimension of $V$ is even and let
$$
V \= V^+\oplus V^-
$$
be a {\it polarization}, \index{polarization} that is
$q(V^+)=q(V^-)=0$. Over $\C$ polarizations always exist for even
dimensional spaces. By the nondegeneracy of $q$ it follows that
to any $v\in V^+$ there is a unique $\hat{v}\in V^-$ such that
$q(v,\hat{v})=-1$. Further, let $V^{-,v}$ be the space of all
$w\in V^-$ such that $q(v,w)=0$, then
$$
V^- \= \C \hat{v}\oplus V^{-,v}.
$$
Let
$$
S\= \wedge^* V^- \= \C\oplus V^-\oplus \wedge^2
V^-\oplus\dots\oplus\wedge^{top}V^-,
$$
then we define an action of $Cl(q)$ on $S$ in the following way:
\begin{itemize}
\item
for $v\in V^-$ and $s\in S$ let
$$
v.s\= v\wedge s,
$$
\item
for $v\in V^+$ and $s\in \wedge^*V^{-,v}$ let
$$
v.s\= 0,
$$
\item
and for $v\in V^+$ and $s\in S$ of the form $s=\hat{v}\wedge s'$
with $s'\in \wedge^*V^{-,v}$ let
$$
v.s\= s'.
$$
\end{itemize}
By the universal property of $Cl(V)$ this extends to an action of
$Cl(q)$. The module $S$ is called the {\it spin
module}.\index{spin module} The induced action of $Spin(q)$
leaves invariant the subspaces
$$
S^+\=\wedge^{even}V^-,\ \ \ S^-\=\wedge^{odd}V^-,
$$
the representation of $Spin(q)$ on these spaces are called the
{\it half spin representations}.\index{half spin representations}
Let $SO(q)^+$ the subgroup of all elements in $SO(q)$ that leave
stable the decomposition $V=V^+\oplus V^-$. This is a connected
reductive group isomorphic to $GL(V^+)$, since, let $g\in
GL(V^+)$ and define $\hat{g}\in GL(V^-)$ to be the inverse of the
transpose of $g$ by the pairing induced by $q$ then the map
$Gl(v)\ra SO(q)^+$ given by $g\mapsto (g,\hat{g})$ is an
isomorphism. In other words, choosing a basis on $V^+$ and a the
dual basis on $V^-$ we get that $q$ is given in that basis by
$\matrix{0}{\1}{\1}{0}$. Then $SO(q)^+$ is the image of the
embedding
\begin{eqnarray*}
GL(V^-) &\hookrightarrow& SO(q)\\
A &\mapsto& \matrix{A}{0}{0}{^tA^{-1}}.
\end{eqnarray*}

Let $Spin(q)^+$ be the inverse image of $SO(q)^+$ in $Spin(q)$.
Then the covering $Spin(q)^+\ra SO(q)^+\cong GL(V^-)$ is the
``square root of the determinant'', i.e., it is isomorphic to the
covering $\tilde{GL}(V^-)\ra GL(V^-)$ given by the pullback
diagram of linear algebraic groups:
$$
\begin{array}{ccc}
\tilde{GL}(V^-)     & \ra   & GL(1)\\
\downarrow         & {} &\downarrow x\mapsto x^2\\
GL(V^-)         & \begin{array}{c}\det\\ \ra\\ {}\end{array}&
GL(1).
\end{array}
$$
As a set, $\tilde{GL}(V^-)$ is given  as the set of all pairs
$(g,z)\in GL(V^-)\times GL(1)$ such that $\det(g)=z^2$ and the
maps to $GL(V^-)$ and $GL(1)$ are the respective projections.

\begin{lemma}\label{epsilon}
There is a one dimensional representation $\epsilon$ of
$Spin(q)^+$ such that
$$
S^\pm\otimes\epsilon\ \cong\ \wedge^\pm V
$$
as $Spin(q)^+$-modules, where $\wedge^\pm$ means the even or odd
powers respectively.
\end{lemma}

\prf Since $Spin(q)^+$ is a connected reductive group over $\C$
we can apply highest weight theory. If the weights of the
representation of $Spin(q)^+$ on $V$ are given by
$\pm\mu_1,\dots,\pm\mu_m$, then the weights of the half spin
representations are given by
$$
\rez{2}(\pm\mu_1\pm\dots\pm\mu_m)
$$
with an even number of minus signs in the one and an odd number
in the other case. Let $\epsilon = \rez{2}(\mu_1+\dots +\mu_m)$
then $\epsilon$ is a weight for $Spin(q)^+$ and $2\epsilon$ is
the weight of, say, the one dimensional representation on
$\wedge^{top}V^+$. By Weyl's dimension formula this means that
$2\epsilon$ is invariant under the Weyl group and therefore
$\epsilon$ is. Again by Weyl's dimension formula it follows that
the representation with highest weight $\epsilon$ is one
dimensional. Now it follows that $S^+\otimes\epsilon$ has the
same weights as the representation on $\wedge^+ V$, hence must be
isomorphic to the latter. The case of the minus sign is analogous.
\qed

\section{Orbital integrals}
It now will be shown that $\tr\pi(f_\tau)$ vanishes for a
principal series representation $\pi$. To this end let $P=MAN$ be
a nontrivial parabolic subgroup with $A\subset \exp(\p_0)$. Let
$(\xi ,V_\xi)$ denote an irreducible unitary representation of $M$
and $e^\nu$ a quasicharacter of $A$. Let $\pi_{\xi ,\nu}:= {\rm
Ind}_P^G (\xi \otimes e^{\nu}\otimes 1)$.

\begin{lemma} \label{pivonggleichnull}
We have $\tr\pi_{\xi ,\nu}(f_\tau) =0$.
\end{lemma}

\prf By Frobenius reciprocity we have for any irreducible unitary
representation $\ga$ of $K$:
$$
\Hom_K(\ga ,\pi_{\xi ,\nu}|_K) \cong \Hom_{K_M}(\ga |_{K_M},\xi ),
$$
where $K_M := K\cap M$. This implies that $\tr\pi_{\xi
,\nu}(f_\tau)$ does not depend on $\nu$. On the other hand
$\tr\pi_{\xi ,\nu}(f_\tau)\ne 0$ for some $\nu$ would imply
$\pi_{\xi ,\nu}(C) =\breve{\tau} (C_K) -B(\rho)+B(\rho_K)$ which
only can hold for $\nu$ in a set of measure zero.
\qed

Recall that an element $g$ of $G$ is called {\it
elliptic}\index{elliptic} if it lies in a compact Cartan subgroup.
Since the following relies on results of Harish-Chandra which were
proven under the assumption that $G$ is of inner type, we will
from now on assume this.

\begin{theorem} \label{orbitalint}
Assume that $G$ is of inner type. Let $g$ be a semisimple element
of the group $G$. If $g$ is not elliptic, then the orbital
integral $\O_g(f_\tau)$ vanishes. If $g$ is elliptic we may assume
$g\in T$, where $T$ is a Cartan in $K$ and then we have
$$
\O_g(f_\tau) \= \tr\ \tau(g)\ c_g^{-1}|W(\t ,\g_g)| \prod_{\alpha
\in \Phi_g^+}(\rho_g ,\alpha),
$$
where $c_g$ is Harish-Chandra's constant, it does only depend on
the centralizer $G_g$ of g. Its value is given in \ref{notations}.

\end{theorem}

\prf The vanishing of $\O_g(f_\tau)$ for nonelliptic semisimple
$g$ is immediate by the lemma above and Theorem \ref{inv-orb-int}.
So consider $g\in T\cap G'$, where $G'$ denotes the set of regular
elements. Note that for regular $g$ the claim is $\CO_g(f_\tau)
=\tr \tau(g)$. Assume the claim proven for regular elements, then
the general result follows by standard considerations as in
\cite{HC-DS}, p.32 ff. where however different Haar-measure
normalizations are used that produce a factor $[G_g:G_g^0]$,
therefore these standard considerations are now explained. Fix
$g\in T$ not necessarily regular. Let $y\in T^0$ be such that $gy$
is regular. Then
\begin{eqnarray*}
\tr \tau (gy) &=& \int_{T\bs G} f_\tau(x^{-1}gyx) dx\\
    &=&  \int_{T^0\bs G} f_\tau(x^{-1}gyx) dx\\
    &=&  \int_{G_g\bs G}\int_{T^0\bs G_g} f_\tau(x^{-1}z^{-1}gyzx)\ dz\ dx\\
    &=&  \int_{G_g\bs G}\sum_{\eta :G_g /G_g^0}
\rez{[G_g:G_g^0]}\int_{T^0\bs G_g^0}
f_\tau(x^{-1}\eta^{-1}z^{-1}gyz\eta x)\ dz\ dx.
\end{eqnarray*}
The factor $\rez{[G_g:G_g^0]}$ comes in by the Haar-measure
normalizations. On $G_g^0$ consider the function
$$
h(y) = f(x^{-1}\eta^{-1} yg\eta x).
$$
Now apply Harish-Chandra's operator $\omega_{G_g}$ to $h$ then
for the connected group $G_g^0$ it holds
$$
h(1) = \lim_{y\ra 1} c_{G_g^0}\ \omega_{G_g^0}\ \CO_y^{G_g^0}(h).
$$
When $y$ tends to $1$ the $\eta$-conjugation drops out and the
claim follows.

So in order to prove the proposition one only has to consider the
regular orbital integrals. Next the proof will be reduced to the
case when the compact Cartan $T$ meets all connected components
of $G$. For this let $G^+=TG^0$ and assume the claim proven for
$G^+$. Let $x\in G$ then $xTx^{-1}$ again is a compact Cartan
subgroup. Since $G^0$ acts transitively on all compact Cartan
subalgebras it follows that $G^0$ acts transitively on the set of
all compact Cartan subgroups of $G$. It follows that there is a
$y\in G^0$ such that $xTx^{-1} =yTy^{-1}\subset TG^0 =G^+$, which
implies that $G^+$ is normal in $G$.

Let $\tau^+=\tau |_{G^+\cap K}$ and $f^+_{\tau^+}$ the
corresponding Euler-Poincar\'e function on $G^+$.

\begin{lemma}
$f^+_{\tau^+} = f_\tau |_{G^+}$
\end{lemma}

Since the Euler-Poincar\'e function is not uniquely determined the
claim reads that the right hand side is a EP-function for $G^+$.

\prf Let $\tau^+ =\tau |_{K^+}$, where $K^+ = TK^0 = K\cap G^+$.
Let $\ph^+\in(C_c^\infty(G^+)\otimes V_\tau)^{K^+}$, which may be
viewed as a function $\ph^+ : G^+ \ra V_\tau$ with $\ph^+(xk)
=\tau(k^{-1})\ph^+(x)$ for $x\in G^+$, $k\in K^+$. Extend $\ph^+$
to $\ph : G\ra V_\tau$ by $\ph(xk)=\tau(k^{-1})\ph^+(x)$ for
$x\in G^+$, $k\in K$. This defines an element of
$(C_c^\infty(G)\otimes V_\tau)^K$ with $\ph |_{G^+}=\ph^+$. Since
$C_\tau$ is a differential operator it follows $f(C_\tau)\ph
|_{G^+} =f(C_\tau)\ph^+$, so
$$
(\ph * \tilde{f}_\tau) |_{G^+} = \ph^+ * \tilde{f}_{\tau^+}.
$$
Considering the normalizations of Haar measures gives the lemma.
\qed

For $g\in T'$ we compute
\begin{eqnarray*}
\CO_g(f_\tau) &=& \int_{T\bs G} f_\tau (x^{-1} gx) dx\\
    &=& \sum_{y:G/G^+} \rez{[G:G^+]} \int_{T\bs G^+} f_\tau (x^{-1} y^{-1} g yx) dx,
\end{eqnarray*}
 where the factor $\rez{[G:G^+]}$ stems from normalization of
Haar measures and we have used the fact that $G^+$ is normal. The
latter equals
$$
\rez{[G:G^+]}\sum_{y:G/G^+}\CO^{G^+}_{y^{-1}xy}(f_\tau) =
\rez{[G:G^+]}\sum_{y:G/G^+}\CO^{G^+}_{y^{-1}xy}(f_{\tau^+}^+).
$$
Assuming the proposition proven for $G^+$, this is
$$
\rez{[G:G^+]}\sum_{y:G/G^+} \tr \tau(y^{-1}gy) =\tr\tau(g).
$$

From now on one thus may assume that the compact Cartan $T$ meets
all connected components of $G$. Let $(\pi,V_\pi)\in\hat{G}$.
Harish-Chandra has shown that for any $\ph\in C^\infty_c(G)$ the
operator $\pi(\ph)$ is of trace class and there is a locally
integrable conjugation invariant function $\Theta_\pi$ on $G$,
smooth on the regular set such that
$$
\tr \pi(\ph) = \int_G\ph(x)\Theta_\pi(x) dx.
$$
For any $\psi\in C^\infty(K)$ let
$\pi|_K(\psi)=\int_K\psi(k)\pi(k) dk$.

\begin{lemma} \label{character}
Assume $T$ meets all components of $G$. For any $\psi\in
C^\infty(K)$ the operator $\pi|_K(\psi)$ is of trace class and
for $\psi$ supported in the regular set $K' =K\cap G'$ we have
$$
\tr \pi|_K(\psi) = \int_K \psi(k) \Theta_\pi(k) dk.
$$
(For $G$ connected this assertion is in \cite{AtSch} p.16.)
\end{lemma}

\prf Let $V_\pi =\bigoplus_i V_\pi(i)$ be the decomposition of
$V_\pi$ into $K$-types. This is stable under $\pi|_K(\psi)$.
Harish-Chandra has proven $[\pi|_K :\tau]\le\dim \tau$ for any
$\tau\in\hat{K}$. Let $\psi =\sum_j\psi_j$ be the decomposition
of $\psi$ into $K$-bitypes. Since $\psi$ is smooth the sequence
$\norm{\psi_j}_1$ is rapidly decreasing for any enumeration of
the $K$-bitypes. Here $\norm{\psi}_1$ is the $L^1$-norm on $K$.
It follows that the sum $\sum_i\tr(\pi|_K(\psi)|V_\pi(i))$
converges absolutely, hence $\pi|_K(\psi)$ is of trace class.

Now let $S=\exp(\p_0)$ then $S$ is a smooth set of
representatives of $G/K$. Let $G.K=\cup_{g\in
G}gKg^{-1}=\cup_{s\in S}sKs^{-1}$, then, since $G$ has a compact
Cartan, the set $G.K$ has non-empty interior. Applying the Weyl
integration formula to $G$ and backwards to $K$ gives the
existence of a smooth measure $\mu$ on $S$ and a function $D$ with
$D(k)>0$ on the regular set such that
$$
\int_{G.K}\ph(x) dx = \int_S \int_K \ph(sks^{-1}) D(k) dk d\mu(s)
$$
for $\ph\in L^1(G.K)$. Now suppose $\ph\in C_c^\infty(G)$ with
support in the regular set. Then
\begin{eqnarray*}
\tr \pi(\ph) &=& \int_{G.K}\ph(x) \Theta_\pi(x) dx\\
    &=& \int_S \int_K \ph(sks^{-1}) D(k) \Theta_\pi(k)d\mu(s)\\
    &=& \int_K\int_S \ph^s(k) d\mu(s) D(k) \Theta_\pi(k) dk,
\end{eqnarray*}
where we have written $\ph^s(k)=\ph(sks^{-1})$. On the other hand
\begin{eqnarray*}
\tr\pi(\ph) &=& \tr \int_{G.K} \ph(x) \pi(x) dx\\
    &=& \tr \int_S\int_K \ph(sks^{-1})D(k) \pi(sks^{-1})dk d\mu(s)\\
    &=& \tr \int_S \pi(s) \pi|_K(\ph^s D)\pi(s)^{-1} d\mu(s)\\
    &=& \int_S \tr \pi |_K(\ph^s D) d\mu(s)\\
    &=& \tr\pi |_K\left( \int_S \ph^s d\mu(s) D\right).
\end{eqnarray*}
This implies the claim for all functions $\psi\in C_c^\infty(K)$
which are of the form
$$
\psi(k) = \int_S\ph(sks^{-1})d\mu(s) D(k)
$$
for some $\ph\in C_c^\infty(G)$ with support in  the regular set.
Consider the map
$$
\begin{array}{cccc}
F:& S\times K' & \ra & G.K'\\
{}& (s,k) & \mapsto & sks^{-1}
\end{array}
$$
Then the differential of $F$ is an isomorphism at any point and by
the inverse function theorem $F$ locally is a diffeomorphism. So
let $U\subset S$ and $W\subset K'$ be open sets such that
$F|_{U\times W}$ is a diffeomorphism. Then let $\alpha\in
C_c^\infty(U)$ and $\beta\in C_c^\infty(W)$, then define
$$
\phi(sks^{-1}) = \alpha(s)\beta(k)\ \ \ {\rm if}\ s\in U,\ k\in W
$$
and $\ph(g)=0$ if $g$ is not in $F(U\times W)$. We can choose the
function $\alpha$ such that $\int_S\alpha (s) d\mu(s) =1$. Then
$$
\int_S\ph(sks^{-1}) d\mu(s) D(k) = \beta(k) D(k).
$$
Since $\beta$ was arbitrary and $D(k) >0$ on $K'$ the lemma
follows.
\qed

Let $W$ denote the virtual $K$-representation on $\wedge^{even}\p
\otimes V_\tau - \wedge^{odd}\p\otimes V_\tau$ and write $\chi_W$
for its character.

\begin{lemma} \label{fin-lin-comb}
Assume $T$ meets all components of $G$, then for any
$\pi\in\hat{G}$ the function $\Theta_\pi\chi_W$ on $K'=K\cap G'$
equals a finite integer linear combination of $K$-characters.
\end{lemma}

\prf It suffices to show the assertion for $\tau =1$. Let $\ph$
be the homomorphism $K\ra O(\p)$ induced by the adjoint
representation, where the orthogonal group is formed with respect
to the Killing form. We claim that $\ph(K)\subset SO(\p)$, the
subgroup of elements of determinant one. Since we assume $K=K^0
T$ it suffices to show $\ph(T)\subset SO(\p)$. For this let $t\in
T$. Since $t$ centralizes $\t$ it fixes the  decomposition $\p
=\oplus_\alpha \p_\alpha$ into one dimensional root spaces. So
$t$ acts by a scalar, say $c$ on $\p_\alpha$ and by $d$ on
$\p_{-\alpha}$. There is $X\in\p_\alpha$ and $Y\in\p_{-\alpha}$
such that $B(X,Y)=1$. By the invariance of the Killing form $B$
we get
$$
1 = B(X,Y) = B(\Ad(t)X,\Ad(t)Y) = cdB(X,Y) =cd.
$$
So on each pair of root spaces $\Ad(t)$ has determinant one hence
also on $\p$.

Replacing $G$ by a double cover if necessary, which doesn't effect
the claim of the lemma, we may assume that $\ph$ lifts to the spin
group $\Spin(\p)$. Let $\p=\p^+\oplus\p^-$ be the decomposition
according to an ordering of $\phi(\t ,\g)$. This decomposition is
a polarization of the quadratic space $\p$ and hence the spin
group acts on $S^+=\wedge^{even}\p^+$ and $S^-=\wedge^{odd}\p^+$
in a way that the virtual module $(S^+-S^-)\otimes (S^+-S^-)$
becomes isomorphic to $W$. For $K$ connected the claim now follows
from \cite{AtSch} (4.5). An inspection shows however that the
proof of (4.5) in \cite{AtSch}, which is located in the appendix
(A.12), already applies when we only assume that the homomorphism
$\ph$ factors over the spin group. \qed

We continue the proof of the proposition. Let $\hat{T}$ denote the
set of all unitary characters of $T$. Any regular element
$\hat{t}\in \hat{T}$ gives rise to a discrete series
representation $(\omega ,V_\omega)$ of $G$. Let
$\Theta_{\hat{t}}=\Theta_\omega$ be its character which, due to
Harish-Chandra, is known to be a function on $G$. Harish-Chandra's
construction gives a bijection between the set of discrete series
representations of $G$ and the set of $W(G,T)=W(K,T)$-orbits of
regular characters of $T$.

Let $\Phi^+$ denote the set of positive roots of $(\g,\t)$ and let
$\Phi^+_c$ ,$\Phi^+_n$ denote the subsets of compact and
noncompact positive roots. For each root $\alpha$ let $t\mapsto
t^\alpha$ denote the corresponding character
 on $T$.
Define
\begin{eqnarray*}
'\Delta_c(t) &=& \prod_{\alpha\in\Phi^+_c}(1-t^{-\alpha})\\
'\Delta_n &=& \prod_{\alpha\in\Phi^+_n}(1-t^{-\alpha})
\end{eqnarray*}
and $'\Delta = '\Delta_c '\Delta_n$. If $\hat{t}\in \hat{T}$ is
singular, Harish-Chandra has also constructed an invariant
distribution $\Theta_{\hat{t}}$ which is a virtual character on
$G$. For $\hat{t}$ singular let $W(\hat{t})\subset W(\g,\t)$ be
the isotropy group. One has $\Theta_{\hat{t}} = \sum_{w\in
W(\hat{t})} \epsilon(w)\Theta'_{w,\hat{t}}$ with
$\Theta'_{w,\hat{t}}$ the character of an induced representation
acting on some Hilbert space $V_{w,\hat{t}}$ and $\epsilon(w)\in
\{\pm 1\}$. Let $\CE_2(G)$ denote the set of discrete series
representations of $G$ and $\CE_2^s(G)$ the set of $W(G,T)$-orbits
of singular characters.

By Theorem \ref{inv-orb-int} the proposition will follow from the

\begin{lemma}\label{trace-tau}
For $t\in T$ regular we have
$$
\tr \tau(t)\=
\rez{|W(G,T)|}\sum_{\hat{t}\in\hat{T}}\Theta_{\hat{t}}(f_\tau)
\Theta_{\hat{t}}(t).
$$
\end{lemma}

\prf Let $\ga$ denote the virtual $K$-representation on
$(\wedge^{even}\p -\wedge^{odd}\p)\otimes V_\tau$. Harish-Chandra
has shown (\cite{HC-S} Theorem 12) that for any
$\hat{t}\in\hat{T}$ there is an irreducible unitary representation
$\pi_{\hat{t}}^0$ such that $\Theta_{\hat{t}}$ coincides up to
sign with the character of $\pi_{\hat{t}}^0$ on the set of
elliptic elements of $G$ and $\pi_{\hat{t}}^0=\pi_{\hat{t}'}^0$ if
and only if there is a $w\in W(G,T)=W(K,T)$ such that
$\hat{t}'=w\hat{t}$.

Further (\cite{HC-S}, Theorem 14) Harish-Chandra has shown that
the family
$$
\left( \frac{'\lap(t)\Theta_{\hat{t}}(t)}
            {\sqrt{|W(G,T)|}}\right)_{\hat{t}\in \hat{T}/W(G,T)}
$$
forms an orthonormal basis of $L^2(T)$. Here we identify
$\hat{T}/W(G,T)$ to a set of representatives in $\hat{T}$ to make
$\Theta_{\hat{t}}$ well defined.

Consider the function $g(t)=\frac{\tr\ga(t) \
'\lap_c(t)}{\overline{\ '\lap_n(t)}} = \tr\tau(t) \ '\lap(t)$. Its
coefficients with respect to the above orthonormal basis are
\begin{eqnarray*}
\langle g, \frac{\ '\lap\Theta_{\hat{t}}}
            {\sqrt{|W(G,T)|}}\rangle &=&
            \rez{\sqrt{|W(G,T)|}}\int_T\tr\ga(t) |
            '\lap_c(t)|^2\overline{\Theta_{\hat{t}}(t)} dt\\
 &=& \sqrt{|W(G,T)|}\int_K \tr\ga(k)\overline{\Theta_{\hat{t}}(k)}
 dk
\end{eqnarray*}
where we have used the Weyl integration formula for the group $K$
and the fact that $W(G,T)=W(K,T)$. Next by Lemma
\ref{fin-lin-comb} this equals
$$
\sqrt{|W(G,T)|}
\dim((\wedge^{even}\p-\wedge^{odd}\p)\otimes\breve{\tau}\otimes\pi_{\hat{t}}^0)^K
\= \sqrt{|W(G,T)|}\Theta_{\hat{t}}(f_\tau).
$$
Hence
\begin{eqnarray*}
g(t) &=& \tr\tau(t) \ '\lap(t)\\
    &=& \sum_{\hat{t}\in \hat{T}/W(G,T)}\Theta_{\hat{t}}(f_\tau)
    \ '\lap(t) \Theta_{\hat{t}}(t)\\
    &=& \rez{|W(G,T)|}\sum_{\hat{t}\in \hat{T}}\Theta_{\hat{t}}(f_\tau)
    \ '\lap(t) \Theta_{\hat{t}}(t).
\end{eqnarray*}
The lemma and the proposition are proven. \qed

\begin{corollary}
If $\tilde{g}\in\tilde{G}$ is semisimple and not elliptic then
$\CO_{\tilde{g}}(g_\tau)=0$. If $\tilde{g}$ is elliptic regular
then
$$
\CO_{\tilde{g}}(g_\tau)\=\frac{\tr\tau(\tilde{g})} {\tr(\tilde{g}
| S^+-S^-)}.
$$
\end{corollary}

\prf Same as for the last proposition with $g_\tau$ replacing
$f_\tau$.
\qed

\begin{proposition}
Assume that $\tau$ extends to a representation of the group $G$
on the same space. For the function $f_{\tau}$ we have for any
$\pi \in \hat{G}$:
$$
\tr \ \pi(f_{{\tau}}) \= \sum_{p=0}^{\dim \ \g /\k}(-1)^p \dim \
{\rm Ext}_{(\g ,K)}^p (V_{{\tau}} ,V_\pi),
$$
i.e., $f_{{\tau}}$ gives the Euler-Poincar\'e numbers of the $(\g
,K)$-modules $(V_{{\tau}} ,V_\pi)$, this justifies the name
Euler-Poincar\'e function.
\end{proposition}

\prf By definition it is clear that
$$
\tr \ \pi (f_{{\tau}}) \= \sum_{p=0}^{\dim \ \p} (-1)^p \dim \
H^p(\g ,K,V_{\breve{\tau}} \otimes V_\pi).
$$
The claim now follows from \cite{BorWall}, p. 52. \qed

\chapter{The Selberg trace formula}
In this chapter we will fix the basic notation and set up
the trace formula. For compactly supported functions this
formula is easily deduced. In the sequel we however will
need it for functions with noncompact support and
therefore will have to show more general versions of the
trace formula.

Let $G$ denote a real reductive group.

\section{The trace formula}
Let $\Ga\subset G$\index{$\Ga$} be a discrete subgroup such that
the quotient manifold $\Ga\bs G$ is compact. We say that $\Ga$ is
{\it cocompact}\index{cocompact} in $G$. Examples are given by
nonisotropic arithmetic groups \cite{marg}. Since $G$ is
unimodular the Haar measure on $G$ induces a $G$-invariant measure
on $\Ga\bs G$, so we can form the Hilbert space $L^2(\Ga\bs G)$
\index{$L^2(\Ga\bs G)$} of square integrable measurable functions
on $\Ga\bs G$ modulo null functions. More generally, let
$(\omega,V_\omega)$\index{$(\omega,V_\omega)$} be a finite
dimensional unitary representation of $\Ga$ and let $L^2(\Ga\bs
G,\omega)$\index{$L^2(\Ga\bs G,\omega)$} be the Hilbert space of
all measurable functions $f : G\ra V_\omega$ such that $f(\ga
x)=\omega(\ga)f(x)$ for all $\ga\in\Ga$ and all $x\in G$ and such
that
$$
\int_{\Ga\bs G} \norm{f(x)}^2 dx\ <\ \infty
$$
modulo null functions. The scalar product of $f,g\in L^2(\Ga\bs
G,\omega)$ is
$$
(f,g) \= \int_{\Ga\bs G} \langle f(x),g(x)\rangle dx,
$$
where $\langle .,.\rangle$ is the scalar product on $V_\omega$.
Let $C^\infty(\Ga\bs G,\omega)$\index{$C^\infty(\Ga\bs G,\omega)$}
be the subspace consisting of smooth functions.

The group $G$ acts unitarily on $L^2(\Ga\bs G,\omega)$ by
$$
R(y)\ph(x)\= \ph(xy)
$$
for $x,y\in G$ and $\ph\in L^2(\Ga\bs G,\omega)$. Let $\pi$ be
any unitary representation of $G$, any $f\in L^1(G)$ defines a
bounded operator
$$
\pi(f) \= \int_G f(x) \pi(x) dx
$$
on the space of $\pi$. \index{$\pi(f)$} We apply this to the case
$\pi=R$. Let $C_c^\infty(G)$ be the space of all smooth functions
of compact support on $G$. Let $f\in C_c^\infty(G)$, $\ph\in
L^2(\Ga\bs G,\omega)$. Fix a fundamental domain $\CF\subset G$
for the $\Ga$-action on $G$ and compute formally at first:
\begin{eqnarray*}
R(f)\ph(x) &=& \int_G f(y) \ph(xy) dy\\
    &=& \int_G f(x^{-1}y) \ph(y) dy\\
    &=& \sum_{\ga\in\Ga} \int_{\ga\CF}f(x^{-1}y)\ph(y)dy\\
    &=& \sum_{\ga\in\Ga}\int_\CF f(x^{-1}\ga y) \omega(\ga)\ph(y) dy\\
    &=& \int_{\Ga\bs G} k_f(x,y)\ph(y) dy,
\end{eqnarray*}
where $k_f(x,y) =\sum_{\ga\in\Ga}f(x^{-1}\ga y)\omega(\ga)$.
Since $f$ has compact support the latter sum is locally finite
and therefore defines a smooth Schwartz kernel on the compact
manifold $\Ga\bs G$. This implies that the operator $R(f)$ is a
smoothing operator and hence of trace class.

Since the convolution algebra $C_c^\infty$ contains an approximate
identity we can infer that the unitary $G$-representation $R$ on
$L^2(\Ga\bs G,\omega)$ decomposes into a direct sum of
irreducibles, i.e.,
$$
L^2(\Ga\bs G,\omega) \ \cong\ \bigoplus_{\pi\in\hat{G}} N_{\Ga
,\omega}(\pi)\pi
$$
with finite multiplicities $N_{\Ga ,\omega}(\pi)\in\N_0$.
\index{$N_{\Ga ,\omega}(\pi)$} Moreover, the trace of $R(f)$ is
given by the integral over the diagonal, so
$$
\tr R(f) \= \int_{\Ga\bs G} \tr k_f(x,x) dx,
$$
where the trace on the right hand side is the trace in
$\End(V_\omega)$. We plug in the sum for $k_f$ and rearrange that
sum in that we first sum over all conjugacy classes in the group
$\Ga$. We write $\Ga_\ga$ and $G_\ga$ for the centralizers of
$\ga\in\Ga$ in $\Ga$ and in $G$ resp.
\begin{eqnarray*}
\tr R(f) &=& \sum_{\ga\in\Ga} \tr\omega(\ga) \int_\CF f(x^{-1}\ga x) dx\\
    &=& \sum_{[\ga]}\sum_{\sigma\in\Ga /\Ga_\ga} \tr\omega(\ga)
        \int_\CF f((\sigma x)^{-1}\ga \sigma x) dx\\
    &=& \sum_{[\ga]}\tr\omega(\ga)
        \int_{\Ga_\ga\bs G} f(x^{-1}\ga x) dx\\
    &=& \sum_{[\ga]}\tr\omega(\ga)\vol(\Ga_\ga\bs G_\ga) \CO_\ga(f),
\end{eqnarray*}
where $\CO_\ga(f)$ is the orbital integral. We have proved the following
proposition.

\begin{proposition}
For $f\in C_c^\infty(G)$ we have by absolute convergence of all
sums and integrals:
$$
\sum_{\pi\in\hat{G}} N_{\Ga,\omega}(\pi)\tr\pi(f) \= \sum_{[\ga]}
\tr\omega(\ga) \vol(\Ga_\ga\bs G)\CO_\ga(f).
$$
\end{proposition}
\qed

All this is classical and may be found at various places. For our
applications we will need to extend the range of functions $f$ to
be put into the trace formula. For this sake we prove:

\begin{proposition}\label{convergence-of-tf}
Suppose $f\in C^{2N}(G)$ with $N>\frac{\dim G}{2}$ and suppose
that $Df\in L^1(G)$ for any $D\in U(\g)$ with $\deg(D)\le 2N$.
 Then the trace formula is valid for
$f$.
\end{proposition}

\prf We first prove this under the additional assumption that either $f\ge
0$ or that the geometric side of
the trace formula converges for $\omega =1$ and with $f$ replaced
by the absolute value $|f|$.

 We consider $U(\g)$ as the algebra of all left invariant
differential operators on $G$. Choose a left invariant Riemannian
metric on $G$ and let $\lap$ denote the corresponding Laplace
operator. Then $\lap\in U(\g)$ and thus it makes sense to write
$R(\lap)$, which is an elliptic differential operator of order
$2$ on the compact manifold $\Ga\bs G$, essentially selfadjoint
and non-negative. The theory of pseudodifferential operators
implies that $R(\lap +1)^{-N}$ has a $C^1$-kernel and thus is of
trace class. Let $g=(\lap +1)^Nf$ then $g\in L^1(G)$, so $R(g)$
is defined and gives a continuous linear operator on the Hilbert
space $L^2(\Ga \bs G,\ph)$. We infer that $R(f)=R(\lap
+1)^{-N}R(g)$ is of trace class.

Let $\chi :[0,\infty [\ra [0,1]$ be a monotonic $C^{2N}$-function
with compact support, $\chi\equiv 1$ on $[0,1]$ and
$|\chi^{(k)}(t)|\le 1$ for $k=1,\dots ,2N$. Let
$h_n(x):=\chi(\frac{dist(x,e)}{n})$ for $x\in G$, $n\in \N$. Then
$|Dh_n(x)|\le \frac{C_D}{n}$ for any $D\in\g U(\g)$. Let
$f_n=h_nf$ then $f_n\ra f$locally uniformly.

{\it Claim.} For the $L^1$-norm on $G$ we have
$$
\parallel Df_n -Df \parallel_1 \ra 0
$$
as $n\ra \infty$ for any $D\in U(\g)$ of degree $\le 2N$.

Proof of the claim: By the Poincar\'e-Birkhoff-Witt theorem
$D(f_n)=D(h_nf)$ is a sum of expressions of the type
$D_1(h_n)D_2(f)$ and $D_1$ can be chosen to be the identity
operator or in $\g U(\g)$. The first case gives the summand
$h_nD(f)$ and it is clear that $\norm{Df-h_nDf}_1$ tends to zero
as $n$ tends to infinity. For the rest assume $D_1\in \g U(\g)$.
Then $|D_1(h_n)(x)|\le \frac{c}{n}$ hence
$\norm{D_1(h_n)D_2(f)}_1$ tends to zero because $D_2(f)$ is in
$L^1(G)$. The claim follows.

To prove the lemma we estimate the operator norm as
$$
\norm{R((\lap +1)^Nf_n)-R((\lap +1)^Nf)}\ \le\ \norm{(\lap
+1)^Nf_n-(\lap +1)^Nf}_1
$$
the latter tends to zero according to the claim proven above.
Denoting the trace norm by $\norm{.}_{tr}$ we infer
$$
\norm{R(f_n)-R(f)}_{tr}
$$ $$
= \norm{R(\lap +1)^{-N}(R((\lap +1)^Nf_n)-R((\lap +1)^Nf))}_{tr}
$$ $$
\le \norm{R(\lap +1)^{-N}}_{tr}\norm{R((\lap +1)^Nf_n)-R((\lap
+1)^Nf)}
$$
which tends to zero. Therefore $\tr R(f_n)$ tends to $\tr R(f)$
as $n\ra \infty$. It follows
\begin{eqnarray*}
\sum_{\pi\in\hat{G}} N_{\Ga ,\omega}(\pi)\tr\pi (f) &=& \tr R(f)\\
    &=& \lim_n \tr R(f_n)\\
    &=& \lim_n \sum_{[\ga]} \tr\omega(\ga)\vol(\Ga_\ga\bs G_\ga) \CO_\ga(f_n)\\
    &=& \sum_{[\ga]} \tr\omega(\ga)\vol(\Ga_\ga\bs G_\ga) \CO_\ga(f)
\end{eqnarray*}
which follows by dominated convergence if $f$ satisfies the second
condition in the proposition. If $f$ satisfies the first
condition we first conclude the statement for $\omega=1$ by
monotonic convergence, since then $f_n$ converges to $f$
monotonically from below. Then also the convergence of the
geometric side of the trace formula for $\omega=1$ follows and
the general statement is inferred by dominated convergence again.

This shows the proposition in the case that the additional assumption holds.
Now let $f$ be arbitrary. First there is $\tilde f \in \CC^{2 d}(G)$ such
that $\tilde f\ge |f|$. To construct such a function one proceeds as
follows. First choose a smooth function $b:\R\ra [0,\infty)$ with
$b(x)=|x|$ for $|x|\ge 1$ and $b(x)\ge|x|$ for every $x$. Next let
$b(\eps,x)=\eps b(x/\eps)$ for $\eps>0$. Choose a function $h\in \CC^{2
d}(G)$ with $0<h(x)\le 1$ and let $\tilde f(x)=b(h(x),f(x))$. Then $\tilde
f\in\CC^{2 d}(G)$ and $\tilde f\ge |f|$. By (a) the trace formula is valid
for $\tilde f$, hence the geometric side converges for $\tilde f$, so it
converges for $|f|$, so, by (b) the trace formula is valid for
$f$.
\qed

\chapter{The Lefschetz formula}

In this chapter $G$ will be a connected semisimple Lie
group with finite center.

\section{Setting up the formula}\label{4.1}
Let $K\subset G$ be a maximal compact subgroup with Cartan
involution $\theta$. Let $X=G/K$ denote the symmetric space
attached to $G$. Let $H\subset G$ be a Cartan subgroup. Modulo
conjugation we may assume that $H$ is stable under $\theta$. Then
$H=AB$, where $A$ is a connected split torus and $B$ is a subgroup
of $K$. The double use of the letter $B$ for the group and a
bilinear form on the Lie algebra will not cause any confusion. Fix
a parabolic subgroup $P$ of $G$ with split component $A$. Then
$P=MAN$ where $N$ is the unipotent radical of $P$ and $M$ is
reductive with compact center and finite component group. The
choice of the parabolic $P$ amounts to the same as a choice of a
set of positive roots $\Phi^+=\Phi^+ (\g ,\a)$ in the root system
$\Phi(\g,\a)$ such that for the Lie algebra $\n=Lie_\C(N)$ we have
$\n=\bigoplus_{\alpha\in\Phi^+}\g_\alpha$. Let
$\bar{\n}=\bigoplus_{\alpha\in\Phi^+}\g_{-\alpha}$, $\bar\n_0=\bar\n\cap\g_0$ and
$\bar{N}=\exp(\bar{\n}_0)$ then $\bar{P}=MA\bar{N}$ is the parabolic
{\it opposite} to $P$.\index{opposite parabolic} The root space
decomposition then writes as
$\g=\a\oplus\m\oplus\n\oplus\bar{\n}$. Let $\rho_P$ be the half of
the sum of the positive roots, each weighted with its
multiplicity, i.e., for $a\in A$ we have $a^{2\rho_P}
=\det(a|\n)$. Let $A^-\subset A$\index{$A^-$} denote the negative
Weyl chamber corresponding to that ordering, i.e., $A^-$ consists
of  all $a\in A$ which act contractingly on the Lie algebra
${\n}$. Further let $\overline{A^-}$ be the closure of $A^-$ in
$G$, this is a manifold with boundary. Let $K_M$ be a maximal
compact subgroup of $M$. We may suppose that $K_M=M\cap K$ and
that $K_M$ contains $B$. Fix an irreducible unitary
representation $(\tau ,V_\tau)$ of $K_M$.

Let $L$ be a Lie group. An element $x$ of
$L$ is called \emph{neat}\index{neat} if  for every
 finite dimensional representation $\eta$ of $L$ the linear map $\eta(x)$ 
 has no nontrivial root of unity as an eigenvalue. 
A subset $A$ of $L$ is called neat if each of its members are. Every neat subgroup is torsion free.
Every arithmetic group has a subgroup of finite index which is neat
\cite{Bor}.

\begin{lemma}
Let $x\in L$ be semisimple and neat. Let $L_x$ denote the centralizer of $x$
in
$L$. Then for each $k\in\N$ the connected components of $L_x$ and $L_{x^k}$ coincide.

\end{lemma}

\prf
It suffices to show that the Lie algebras coincide. The Lie algebra of $L_x$ is just
the fixed space of $\Ad(x)$ in Lie$(L)$. Since $\Ad(x)$ is semisimple and does not have a
root of unity for an eigenvalue this fixed space coincides with the fixed space of
$\Ad(x^k)$. The claim follows.
\qed

Let $\Ga\subset G$ denote a cocompact discrete subgroup which is
 neat. Since $\Ga$ is torsion free it acts fixed point free
on the contractible space $X$ and hence $\Ga$ is the fundamental
group of the Riemannian manifold
$$
X_{\Ga} = \Ga \bs X = \Ga \bs G/K
$$\index{$X_\Ga$}
 it follows that we have a
canonical bijection of the homotopy classes of loops:
$$
[S^1 : X_{\Ga} ] \ra \Ga / {\rm conjugacy}.
$$
For a given class $[\ga]$ let $X_\ga$\index{$X_\ga$} denote the
union of all closed geodesics in the corresponding class in $[S^1
: X_\Ga ]$. Then $X_\ga$ is a smooth submanifold of $X_{\Ga_H}$
\cite{DKV}, indeed, it follows that
$$
X_\ga\ \cong\ \Ga_\ga \bs G_\ga/K_\ga,
$$
where $G_\ga$ and $\Ga_\ga$ are the centralizers of $\ga$ in $G$
and $\Ga$ and $K_\ga$ is a maximal compact subgroup of $G_\ga$.
Further all closed geodesics in the class $[\ga]$ have the same
length $l_\ga$.\index{$l_\ga$}

\begin{lemma}\label{neatXgamma}
For $\ga\in\Ga$ and $n\in\N$ we have $X_{\ga^n}=X_\ga$.
\end{lemma}

\prf By the last lemma we have for the connected components, $G_\ga^0=G_{\ga^n}^0$.
By definition one has that $X_{\ga^n}$ is a subset of $X_\ga$. Since both are connected
submanifolds of $X_\Ga$ they are equal if their dimensions are the same.
We have
\begin{eqnarray*}
\dim\, X_\ga &=& \dim\, G_\ga/K_\ga\\
&=& \dim\, G_\ga^0/K_\ga^0\\
&=& \dim\, G_{\ga^n}^0/K_{\ga^n}^0\\
&=& \dim\, G_{\ga^n}/K_{\ga^n}\\
&=&\dim\, X_{\ga^n}. 
\end{eqnarray*}\qed

Let $\CE_P(\Ga)$ denote the set of all conjugacy classes $[\ga]$
in $\Ga$ such that $\ga$ is in $G$ conjugate to an element $a_\ga
b_\ga$ of $A^- B$.

Take a class $[\ga]$ in $\CE_P(\Ga)$. Then there is a conjugate
$H_\ga =A_\ga B_\ga$ of $H$ that contains $\ga$. Then the
centralizer $\Ga_{\ga}$ projects to a lattice $\Ga_{A,\ga}$in the
split part $A_\ga$. Let $\la_\ga$ be the covolume of this lattice.

Let $(\omega,V_\omega)$ be a finite dimensional unitary
representation of $\Ga$.

For any module $V$ of the Lie algebra ${\n}$ let $H^q({\n},V)$,
$q=0,\dots ,\dim\n$ denote the {\it Lie algebra cohomology}
\index{Lie algebra cohomology} \cite{BorWall}. If $\pi\in\hat{G}$
then $H^q({\n},\pi_K)$ is an admissible $(\a\oplus\m,K_M)$-module
of finite length \cite{HeSch}.

\begin{theorem}\label{lefschetz}  (Lefschetz formula)
Let $\ph$ be smooth and compactly supported on ${A^-}$. Assume
that $M$ is orientation preserving or that $\tau$ lies in the
image of the restriction map $\res_{K_M}^M$. Then the expression
$$
\sum_{\pi \in \hat{G}} N_{\Ga,\omega} (\pi) \sum_{p,q}(-1)^{p+q}
 \int_{A^-} \ph (a)
\tr(a|(H^q(\n,\pi_K)\otimes\wedge^p\p_M \otimes
V_{\breve{\tau}})^{K_M})da,
$$
henceforth referred to as the global side, equals
$$
(-1)^{\dim N}\sum_{[\ga]\in {\cal E}_P(\Ga)}\la_\ga \
\chi(A_\ga\bs X_\ga)\ \tr\omega(\ga)\ \frac{\ph(a_\ga)\tr
\tau(b_\ga)}{\det(1- a_\ga b_\ga |{\n})},
$$
called the local side, where $\chi(A_\ga\bs X_\ga)$ denotes the
Euler characteristic of the manifold $A_\ga\bs X_\ga$.
\end{theorem}

The proof will be given in section \ref{prf-lef}

We will give a reformulation of the theorem and for this
we need the following notation. Let $V$ be a complex
vector space on which $A$ acts linearly. For each
$\la\in\a^*$ let $V^\la$ denote the generalized
$\la$-eigenspace in $V$, i.e., $V^\la$ consists of all
$v\in V$ such that there is $n\in \N$ with
$$
(a-a^{\la})^n v\= 0
$$
for every $a\in A$. The theorem above is equivalent to the
following

\begin{theorem}\label{4.2}
Assume that $M$ is orientation preserving or that $\tau$
lies in the image of the restriction map $\res_{K_M}^M$.
Then we have the following identity of distributions on
$A^-$.
$$
\sum_{\pi\in\hat G} N_\Ga(\pi)\sum_{\la\in\a^*}m_\la(\pi)\, (\cdot)^\la\=
\sum_{[\ga]\in\CE_P(\Ga)} c_\ga \delta_{a_\ga}.
$$
Here $m_\la(\pi)$ equals
$$
\sum_{p,q}(-1)^{p+q+\dim N}
\dim\left(H^q(\n,\pi_K)^\la\otimes\wedge^p\p_M \otimes
V_{\breve{\tau}}\right)^{K_M}.
$$
The sum indeed is finite for each $\la\in\a^*$. Further,
for $[\ga]\in {\cal E}_P(\Ga)$ we set
$$
c_\ga\=\la_\ga \ \chi(A_\ga\bs X_\ga)\ \tr\omega(\ga)\
\frac{\tr \tau(b_\ga)}{\det(1- a_\ga b_\ga |{\n})}.
$$
\end{theorem}

\begin{corollary} \label{ruelle-lef}
Let $\ph$ be smooth and compactly supported on ${A^-}$. Then
$$
\sum_{\pi \in \hat{G}} N_{\Ga,\omega} (\pi)
\sum_{p,q,r}(-1)^{p+q+r+\dim N}\times\hspace{150pt}
$$ $$
\hspace{50pt} \int_{A^-} \ph (a)
\tr(a|(H^q(\n,\pi_K)\otimes\wedge^p\p_M \otimes\wedge^r\n\otimes
V_{\breve{\tau}})^{K_M})da,
$$
equals
$$
\sum_{[\ga]\in {\cal E}_P(\Ga)}\la_\ga \ \chi(A_\ga\bs X_\ga)\
\tr\omega(\ga)\ \ph(a_\ga)\tr \tau(b_\ga).
$$
\end{corollary}

\prf Let $\wedge^r\n=\bigoplus_{j\in I_r} V_j$ be the
decomposition of the adjoint action of the group $M$ into
irreducible representations. On $V_j$ the torus $A$ acts
by a character $\la_j$. We apply the theorem to $\tau$
replaced by $\tau\oplus V_j|_{K_M}$ and $\ph(a)$ replaced
by $\ph(a)\la_j(a)$. We sum these over $j$ and take the
alternating sum with respect to $r$. On the local side we
apply the identity
$$
\sum_{r=0}^{\dim N}(-1)^r\tr (a_\ga b_\ga |\wedge^r\n) =
\det(1-a_\ga b_\ga|\n)
$$
to get exactly the local side of the corollary. On the global side
we need to recall that the $K_M$-module $\wedge^r\n$ is self dual.
The corollary follows. \qed

\section{Euler characteristics}
Let $L$ be a real reductive group and suppose there is a finite
subgroup $E$ of the center of $L$ and a reductive and Zariski-connected linear group
 $\CL$ over $\R$ such that $L/E$ is
isomorphic to a subgroup of $\CL(\R)$ of finite index. Note that
these conditions are satisfied whenever $L$ is a Levi component of
a connected semisimple group $G$ with finite center. Let $K_L$ be
a maximal compact subgroup of $L$ and let $\Ga$ be a cocompact
discrete subgroup of $L$. Fix a nondegenerate invariant bilinear
form on the Lie algebra $\l_0$ of $L$ such that $B$ is negative
definite on the Lie algebra of $K_L$ and positive definite on its
orthocomplement. Let $\theta$ be the Cartan involution fixing
$K_L$ pointwise then the form $-B(X,\theta(Y))$ is positive
definite and thus defines a left invariant metric on $L$. For any
closed subgroup $Q$ we get a left invariant metric on $Q$. The
volume element of that metric gives a Haar measure, called the
{\it standard volume} \index{standard volume} with respect to $B$
on $Q$.

Let now $H$ be a $\theta$-stable Cartan subgroup, then $H=AB$,
where $A$ is a connected split torus and $B\subset K_L$ is a
Cartan of $K_L$. On the space $L/H$ we have a pseudo-Riemannian
structure given by the form $B$. The Gauss-Bonnet construction
(\cite{dieu} sect. 24 or see below) generalizes to
Pseudo-Riemannian structures to give an Euler-Poincar\'e measure
$\eta$ on $L/H$. Define a (signed) Haar-measure by
$$
\mu_{EP} \= \eta\otimes ({\rm normalized\ Haar\ measure\ on}\ H).
$$
Let $W=W(L,H)$ denote the Weyl group and let $W_\C=W(L_\C,H_\C)$
be the Weyl group of the complexifications. Let the {\it generic
Euler characteristic}\index{generic Euler characteristic} be
defined by
$$
\chi_{gen}(\Ga \bs L/K_L)\= \frac{\mu_{EP}(\Ga\bs L)}{|W|}.
$$
Write $X_\Ga = \Ga\bs L/K$.

\begin{lemma}
If $H$ is compact and $\Ga$ is torsion-free, then the generic
Euler characteristic equals the ordinary Euler characteristic,
i.e., $\chi_{gen}(X_\Ga)=\chi(X_\Ga)$.
\end{lemma}

\prf In this case we have $H=B$. The Gauss-Bonnet Theorem tells us
$$
\eta(\Ga\bs L/B)\= \chi(\Ga\bs L/B).
$$
Now $\Ga\bs L/B\ra \Ga\bs L/K_L$ is a fiber bundle with fiber
$K_L/B$, therefore
 we get
$$
\chi(\Ga\bs L/B)\= \chi(\Ga\bs L/K_L)\chi(K_L/B).
$$
Finally the Hopf-Samelson formula says $\chi(K_L/B)=|W|$. \qed

For the next proposition assume that $A$ is central in $L$, $L=AL_1$ where $L_1$ has compact
center and that $\Ga$ is neat. Let $C$ denote the center of $L$, then
$A\subset C$ and $C=AB_C$, where $B_C=B\cap C$. Let $L'$ be the
derived group of $L$ and let $\Ga'=L'\cap \Ga C$ and
$\Ga_C=\Ga\cap C$ then by Lemma 3.3 in \cite{Wolf} we infer that
$\Ga_C$ is a torsion free cocompact subgroup of $C$ and $\Ga'$ is
a cocompact discrete subgroup of $L'$. Note that $\Ga$ being
 neat and $\CL$ Zariski connected implies that $\Ga'$ is
torsion-free modulo the center of $L'$. Let $\Ga_A=A\cap \Ga_C
B_C$ the projection of $\Ga_C$ to $A$. Then $\Ga_A$ is discrete
and cocompact in $A$.

\begin{proposition}
Assume $\Ga$ is  neat and $A$ is central in $L$ then
$A/\Ga_A$ acts freely on $X_\Ga$ and $\chi_{gen}(X_\Ga)=\chi(A\bs
X_\Ga)\vol(A/\Ga_A)$.
\end{proposition}

\prf The group $A_\Ga=A/\Ga_A$ acts on $\Ga \bs L/B$ by
multiplication from the right. We claim that this action is free,
i.e., that it defines a fiber bundle
$$
A_\Ga\ra\Ga\bs L/B\ra \Ga\bs L/H.
$$
To see this let $\Ga xaB=\Ga xB$ for some $a\in A$ and $x\in L$,
then $a=x^{-1}\ga xb$ for some $\ga\in\Ga$ and $b\in B$. Writing
$\ga$ as $\ga' a_\ga b_\ga$ with $\ga'\in\Ga'$ and $a_\ga\in A$
and $b_\ga\in B$ we conclude that $\ga'$ must be elliptic, hence
in the center of $L'$, which is a subset of $B$. From $a=a_\ga
b_\ga b$ it follows that $a=a_\ga$, hence the claim.

In the same way we see that we get a fiber bundle
$$
A_\Ga\ra\Ga\bs L/K_L\ra A\Ga\bs L/K_L.
$$
We now apply the Gauss-Bonnet theorem to conclude
\begin{eqnarray*}
\chi_{gen}(X_\Ga) &=& \vol(A/\Ga_A)\frac{\eta(\Ga\bs L/H)}{|W|}\\
    &=& \vol(A/\Ga_A)\frac{\chi(\Ga\bs L/H)}{\chi(K_L/B)}\\
    &=& \vol(A/\Ga_A)\frac{\chi(A\Ga\bs L/B)}{\chi(K_L/B)}\\
    &=& \vol(A/\Ga_A)\chi(A\Ga\bs L/K_L).
\end{eqnarray*}
The proposition is proven.
\qed

We will compute $\chi_{gen}(X_\Ga)$ in terms of root systems. Let
$\Phi$ denote the root system of $(\l,\h)$, where $\l$ and $\h$
are the complexified Lie algebras of $L$ and $H$. Let $\Phi_n$ be
the set of noncompact imaginary roots and choose a set $\Phi^+$ of
positive roots such that for $\alpha\in\Phi^+$ nonimaginary we
have that $\alpha^c\in\Phi^+$. Let $\nu=\dim L/K_L-\rank L/K_L$,
let $\rho$ denote the half of the sum of all positive roots. For
any compact subgroup $U$ of $L$ let $v(U)$ denote the standard
volume.

\begin{theorem} \label{3.3}
The generic Euler number satisfies
\begin{eqnarray*}
\chi_{gen}(X_\Ga)&=& \frac{(-1)^{|\Phi_n^+|}|W_\C|
                \prod_{\alpha\in\Phi^+}(\rho,\alpha)v(K_L)}
            {(2\pi)^{|\Phi^+|}|W|2^{\nu/2}v(B)}
\vol(\Ga\bs L)\\
    &=& c_L^{-1}|W_\C|\prod_{\alpha\in\Phi^+}(\rho,\alpha)\vol(\Ga\bs L),
\end{eqnarray*}
where $c_L$ is Harish-Chandra's constant, i.e.,
$$
c_L\= (-1)^{|\Phi_n^+|}
(2\pi)^{|\Phi^+|}2^{\nu/2}\frac{v(B)}{v(K_L)}|W|.
$$
So, especially in the case when $\Ga$ is neat and $A$ is
central we get
$$
\chi(A\bs X_\Ga)\= \frac{|W_\C|\prod_{\alpha\in\Phi^+}
(\rho,\alpha)}{c_L\vol(A/\Ga_A)}\vol(\Ga\bs L).
$$
\end{theorem}

\prf On the manifold $L/H$ the form $B$ gives the structure of a
pseudo-Riemannian manifold. Let $P$ denote the corresponding
$SO_{p,q}$ fiber bundle, where $(p,q)$ is the signature of $B$ on
$\l_0/\h_0$. Let $\ph :H\ra SO_{p,q}$ denote the homomorphism
induced by the adjoint representation. We the have $P=L\times_\ph
SO_{p,q}$. A connection on $P$ is given by the $L$-invariant
connection $1$-form
$$
\omega \left( A+\sum_\alpha c_\alpha X_\alpha \right) \= A,\ \ \
\ A\in so(p,q),
$$
where we have used $B_eP\cong
so(p,q)\oplus\left((\oplus_\alpha\g_\alpha)\cap\g_0\right)$. By an
inspection in local charts one finds the following formula for the
$L$-invariant $2$-form $d\omega$:
$$
d\omega(e)\left( A+\sum_\alpha c_\alpha X_\alpha,A'+ \sum_\alpha
c_\alpha' X_\alpha\right) = [A,A'] -\sum_\alpha c_\alpha
c_\alpha' \ph_*H_\alpha.
$$
Let $\Omega=d\omega$ and let $Pf$ be the Pfaffian as in
\cite{dieu} 24.46.10. Let $J$ be the diaginal matrix
$\diag(1,\dots,1,-1,\dots,-1)$ with $p$ ones and $q$ minus
ones and let $P(X)=Pf(JX)$. Write this as
$$
P(X)\= (2\pi)^{-m}\sum_{h,k}\chi_{h,k}X_{h_1,k_1}\dots
X_{h_m,k_m},
$$
where $m-\rez{2}(p+q)\in\N$. Let
$$
F(\Omega)\= (2\pi)^{-m} \sum_{h,k}
\chi_{h,k}\Omega_{h_1,k_1}\wedge\dots\wedge\Omega_{h_m,k_m}.
$$
There is a unique $G$-invariant form $F_B(\Omega)$ on $L/H$ which
at the origin is the pullback of $F(\Omega)$ with respect to a
section $s: L/H\ra P$. Writing $\Omega =\Omega_1 -\Omega_2$ with
$$
\Omega_1\left( A+\sum_\alpha c_\alpha X_\alpha\ ,\ A'+ \sum_\alpha
c_\alpha' X_\alpha\right) = [A,A']
$$
we get
$$
F_B(\Omega)\= F_B(-\Omega_2)\= (-1)^m F_B(\Omega_2).
$$
We call this form $\eta$. On the space $\C X_\alpha +\C
X_{-\alpha}$ we have that
$$
\ph_*Y=-\alpha(Y)i\matrix{{}}{1}{1}{{}}
$$
with respect to the
basis $(X_\alpha,X_{-\alpha})$. Let $\alpha_1,\dots,\alpha_m$ be
an enumeration of $\Phi^+$, then
$$
\eta\= (2\pi)^{-m} \sum_{\sigma\in Per(\Phi)}\sum_{h,k}\chi_{h,k}
(\ph_*H_{\sigma(\alpha_1)})_{h_1,k_1}\dots
(\ph_*H_{\sigma(\alpha_m)})_{h_m,k_m} \tilde{\omega}_\sigma,
$$
where
$$
\tilde{\omega}_\sigma\= dX_{\sigma(\alpha_1)}\wedge
dX_{-\sigma(\alpha_1)} \wedge\dots\wedge
dX_{\sigma(\alpha_m)}\wedge dX_{-\sigma(\alpha_m)}.
$$
We end up with
\begin{eqnarray*}
\eta &=& (2\pi)^m(-1)^m \sum_{{\alpha\in\Phi^+},\
{\sigma\in\Per(\Phi^+)}}
(\alpha,\sigma(\alpha))\tilde{\omega}_\sigma\\
    &=& (2\pi)^m(-1)^m |W_\C|\prod_{\alpha\in\Phi^+} (\rho,\alpha)
\tilde{\omega}_\sigma.
\end{eqnarray*}
Comparing $\eta$ with the standard measure given by the form $B$
gives the claim.
\qed

\section{Proof of the Lefschetz formula}\label{prf-lef}
The notations are as in section \ref{4.1}. Let $G$ act on itself
by conjugation, write $g.x = gxg^{-1}$, write $G.x$ for the orbit,
so $G.x = \{ gxg^{-1} | g\in G \}$ as well as $G.S = \{ gsg^{-1} |
s\in S , g\in G \}$ for any subset $S$ of $G$. We are going to
consider functions that are supported in the set $G.(MA^-)$. By
Theorem \ref{exist-ep} there exists an Euler-Poincar\'e function
$f_\tau^M\in C_c^\infty(M)$ to the representation $\tau$.

For $g\in G$ and $V$ any complex vector space on which $g$ acts
linearly let $E(g|V)\subset \R^*_+$ be defined by
$$
E(g|V) \ :=\ \{ |\mu | : \mu {\rm \ is\ an\ eigenvalue\ of}\ g
{\rm \ on\ } V\}.
$$
Let $\la_{min}(g|V):= \min(E(g|V))$ and  $\la_{max}(g):=
\max(E(g|V))$ the minimum and maximum.

For $am\in AM$ define
$$
\la(am) := \frac{\la_{min}(a|\bar{\n})}{\la_{max}(m|\g)^2}.
$$
Note that $\la_{max}(m|\g)$ is always $\ge 1$ and that
$\la_{max}(m|\g)\la_{min}(m|\g)=1$. We will consider the set
$$
(AM)^\sim \ :=\ \{ am\in AM | \la(am)>1 \}.
$$
Let $M_{ell}$ denote the set of elliptic elements in $M$.

\begin{lemma} \label{MA}
The set $(AM)^\sim$ has the following properties:

\begin{itemize}
\item[1.]
$A^-M_{ell}\subset (AM)^{\sim}$
\item[2.]
$am\in (AM)^\sim \Rightarrow a\in A^-$
\item[3.]
$am, a'm' \in (AM)^\sim, g\in G\ {\rm with}\ a'm'=gamg^{-1}
\Rightarrow a=a', g\in AM$.
\end{itemize}
\end{lemma}

\prf The first two are immediate. For the third let $am, a'm' \in
(AM)^\sim$ and $g\in G$ with $a'm'=gamg^{-1}$. Observe that by the
definition of $(AM)^\sim$ we have
\begin{eqnarray*}
\la_{min}(am | \bar{\n}) &\ge& \la_{min}(a| \bar{\n})\la_{min}(m |
\g)\\
 &>& \la_{max}(m|\g)^2\la_{min}(m|\g)\\
 &=& \la_{max}(m|\g)\\
 &\ge & \la_{max}(m | \a +\m +\n)\\ &\ge&\la_{max}(am | \a +\m +\n)
\end{eqnarray*}
that is, any eigenvalue of $am$ on $\bar{\n}$ is strictly bigger
than any eigenvalue on $\a +\m +\n$. Since $\g = \a +\m +\n
+\bar{\n}$ and the same holds for $a'm'$, which has the same
eigenvalues as $am$, we infer that $\Ad(g)\bar{\n} =\bar{\n}$. So
$g$ lies in the normalizer of $\bar{\n}$, which is
$\bar{P}=MA\bar{N} =\bar{N}AM$. Now suppose $g=nm_1a_1$ and
$\hat{m} =m_1mm_1^{-1}$ then
$$
gamg^{-1} \= na\hat{m}n^{-1} \= a\hat{m}\
(a\hat{m})^{-1}n(a\hat{m})\ n^{-1}.
$$
Since this lies in $AM$ we have $(a\hat{m})^{-1}n(a\hat{m}) =n$
which since $am\in (AM)^\sim$ implies $n=1$. The lemma is proven.
\qed

Fix a smooth function $\eta$ on $N$ which has compact support, is
positive, invariant under $K_M$ and satisfies $\int_N\eta(n) dn
=1$. Extend the function $\ph$ from $A^-$ to a conjugation
invariant smooth function $\tilde{\ph}$ on $AM$ such that
$\tilde{\ph}(am)=\ph(a)$ whenever $m$ is elliptic and such that
there is a compact subset $C\subset A^-$ such that $\tilde{\ph}$
is supported in $CM\cap(AM)^\sim$. This is clearly possible since
the support of $\ph$ is compact. It follows that the function
$$
am\ \mapsto\ f_\tau^M(m)\frac{\tilde{\ph}(am)}{\det(1-am|\n)}
$$
is smooth and compactly supported on $AM$. Given these data let
$f = f_{\eta ,\tau ,\ph} : H\ra \C$ be defined by
$$
f (kn ma (kn)^{-1}) := \eta (n) f_\tau^M(m)
\frac{\tilde{\ph}(am)}{\det(1-(ma)|\n)},
$$
for $k\in K, n\in N, m\in M, a\in\overline{A^-}$. Further $f(x)=0$
if $x$ is not in $G.(AM)^\sim$.

\begin{lemma} \label{welldef}
The function $f$ is well defined.
\end{lemma}

\prf By the decomposition $G=KP=KNMA$ every element $x\in
G.(AM)^\sim$ can be written in the form $kn ma (kn)^{-1}$. Now
suppose two such representations coincide, that is
$$
kn ma (kn)^{-1}\= k'n' m'a' (k'n')^{-1}
$$
then by Lemma \ref{MA} we get $(n')^{-1} (k')^{-1}kn\in MA$, or
$(k')^{-1}k\in n'MAn^{-1}\subset MAN$, hence $(k')^{-1}k\in K\cap
MAN=K\cap M=K_M$. Write $(k')^{-1}k=k_M$ and $n''=k_Mnk_M^{-1}$,
then it follows
$$
n'' k_Mmk_M^{-1} a (n'')^{-1}\= n'm'a'(n')^{-1}.
$$
Again by Lemma \ref{MA} we conclude $(n')^{-1}n''\in MA$, hence
$n'=n''$ and so
$$
k_Mmk_M^{-1} a \= m'a',
$$
which implies the well-definedness of $f$.
\qed

We will plug $f$ into the trace formula. For the geometric side
let $\ga \in \Ga$. We have to calculate the orbital integral:
$$
\CO_\ga (f) = \int_{G_\ga \bs G} f(x^{-1}\ga x) dx.
$$
by the definition of $f$ it follows that $\CO_\ga(f)=0$ if
$\ga\notin G.(AM)^\sim$. It remains to compute $\CO_{am}(f)$ for
$am\in(AM)^\sim$. Again by the definition of $f$ it follows
$$
\CO_{am}(f) \=
\CO_m^M(f_\tau^M)\frac{\tilde{\ph}(am)}{\det(1-ma|\n)},
$$
where $\CO_m^M$ denotes the orbital integral in the group $M$.

Since only elliptic elements have nonvanishing orbital integrals
at $f_\tau^M$ it follows that only those conjugacy classes $[\ga]$
contribute for which  $\ga$ is in $G$ conjugate to $a_\ga b_\ga\in
A^-B$. Recall that Theorem \ref{3.3} says
$$
\vol(\Ga_\ga\bs G_\ga)\= \chi(A_\ga\bs X_\ga) \la_\ga
 \frac{c_{G_\ga}}
     {|W_{\ga,\C}|\prod_{\alpha\in\Phi_\ga^+}(\rho_\ga,\alpha)}.
$$
By Theorem \ref{orbitalint} we on the other hand get
$$
\CO_\ga(f) \=
\frac{|W_{\ga,\C}|\prod_{\alpha\in\Phi_\ga^+}(\rho_\ga,\alpha)}
     {c_{G_\ga}}
     \tr\tau(b_\ga)\frac{\ph(a_\ga)}
                        {\det(1-a_\ga t_\ga|\n)},
$$
so that
$$
\vol(\Ga_\ga\bs G_\ga)\, \CO_\ga(f)\= \chi(A_\ga\bs X_\ga) \la_\ga
\tr\tau(b_\ga)\frac{\ph(a_\ga)}
                        {\det(1-a_\ga b_\ga|\n)}.
$$
It follows that the geometric side of the trace formula coincides
with the geometric side of the Lefschetz formula.

Now for the spectral side let $\pi\in\hat{G}$. We want to compute
$\tr\pi(f)$. Let $\Theta_\pi^G$ be the locally integrable
conjugation invariant function  on $G$ such that
$$
\tr\pi(f)\= \int_G f(x) \Theta_\pi^G(x) dx.
$$
To evaluate $\tr \pi(f)$ we will employ the Hecht-Schmid
character formula \cite{HeSch}. For this let
$$
(AM)^- = {\rm interior\ in\ } MA {\rm \ of\ the\ set}
$$ $$
\left\{ g\in MA | \det (1-ga | \n) \geq 0\ {\rm for\ all\ } \a\in
A^- \right\}.
$$
The main result of \cite{HeSch} is that for $ma\in (AM)^- \cap
G^{reg}$, the regular set, we have
$$
\Theta_\pi^G(am) = \frac{\sum_{p=0}^{\dim \n}(-1)^p \Theta_{H_p(\n
,\pi_K)}^{MA}(am)}{\det (1-am |\n)}.
$$
Let $f$ be supported on $G.(MA^-)$, then the Weyl integration
formula implies that
$$
\int_G f(x) dx = \int_{G/MA} \int_{MA^-} f(gmag^{-1})
|det(1-ma|\n\oplus\bar{\n})\,dadmdg.
$$
So that for $\pi\in\hat{G}$:
\begin{eqnarray*}
\tr \pi(f) &=& \int_G \Theta_\pi^G (x) f (x) dx\\
    &=& \int_{MA^-} \Theta_\pi^G (ma) f_\tau^M(m)\tilde{\ph}(am)
|\det(1-ma|\bar{\n})|\,dadm\\
    &=&  \int_{MA^-}f_\tau^M(m)\sum_{p=0}^{\dim
    N}(-1)^p\Theta_{H_p(\n,\pi_K)}^{AM}(am)\\
    &{}&\hspace{20pt}\times
    \frac{|\det(1-am|\bar{\n})|}
         {\det(1-am|\n)}
    \tilde{\ph}(am)\ dadm.
\end{eqnarray*}
Now we find that
\begin{eqnarray*}
|\det(1-am|\bar{\n})| &=& (-1)^{\dim N} \det(1-am|\bar{\n})\\
    &=& (-1)^{\dim N} a^{-2\rho_P}\det(a^{-1}-m|\bar{\n})\\
    &=& (-1)^{\dim N} a^{-2\rho_P}\det((am)^{-1}-1|\bar{\n})\\
    &=& a^{-2\rho_P}\det(1-(am)^{-1}|\bar{\n})\\
    &=& a^{-2\rho_P}\det(1-am|{\n})
\end{eqnarray*}
so that
$$
\tr\pi(f)\= \int_{MA^-}f_\tau^M(m)\sum_{p=0}^{\dim
N}(-1)^p\Theta_{H_p(\n,\pi_K)}^{AM}(am)a^{-2\rho_P}\tilde{\ph}(am)dadm.
$$
We have an isomorphism of $(\a\oplus\m,K_M)$-modules \cite{HeSch}
$$
H_p(\n,\pi_K)\cong H^{\dim N-p}(\n,\pi_K)\otimes \wedge^{top}\n.
$$
This implies
$$
\sum_{p=0}^{\dim N}(-1)^p\Theta_{H_p(\n,\pi_K)}^{AM}(am)
a^{-2\rho_P}\= (-1)^{\dim N} \sum_{p=0}^{\dim
N}(-1)^p\Theta_{H^p(\n,\pi_K)}^{AM}(am).
$$
And so
$$
\tr\pi(f)\= \int_{MA^-}f_\tau^M(m) \sum_{p=0}^{\dim N}(-1)^{p+\dim
N} \Theta_{H^p(\n,\pi_K)}^{AM}(am)\tilde{\ph}(am)\ dadm.
$$
Let $B=H_1,\dots,H_n$ be the conjugacy classes of Cartan subgroups
in $M$. By the Weyl integration formula the integral over $M$ is a
sum of expressions of the form
$$
\int_{H_j}\int_{M/H_j}
f_\tau^M(x^{-1}hx)\Theta_{H^*(\n,\pi_K)}^{MA}(x^{-1}hax)
\tilde{\ph}(x^{-1}hax)dx dh
$$ $$
=\int_{H_j}\int_{M/H_j}
f_\tau^M(x^{-1}hx)\Theta_{H^*(\n,\pi_K)}^{MA}(ha)
\tilde{\ph}(ha)dx dh
$$ $$
=\int_{H_j}\CO_h^M(f_\tau^M)\Theta_{H^*(\n,\pi_K)}^{MA}(ha)
\tilde{\ph}(ha) dh,
$$
where we have used the conjugacy invariance of
$\Theta_{H^*(\n,\pi_K)}^{MA}$ and $\tilde{\ph}$. The orbital
integrals $\CO_h^M(f_\tau^M)$ are nonvanishing only for $h$
elliptic, so only the summand with $H_j=H_1=B$ survives. In this
term we may replace $\tilde{\ph}(ha)$ by $\ph(a)$ so that we get
$$
\tr\pi(f)\= \int_{MA^-}f_\tau^M(m)\sum_{p=0}^{\dim N}(-1)^{p+\dim
N} \Theta_{H^p(\n,\pi_K)}^{MA}(am) \ph(a)da dm
$$ $$
= \sum_{p,q\ge 0} (-1)^{p+q+\dim N}\int_{A^-} \ph(a)
\tr(a|(H^q(\n,\pi_K)\otimes\wedge^p\p_M\otimes
V_{\breve{\tau}})^{K_M}) da.
$$
The theorem follows. \qed

\section{Geometric interpretation}

A {\it foliation} $F$ on a $n$-dimensional connected smooth
manifold $M$ is a decomposition of $M$ into submanifolds, called
leaves: $M=\bigcup_\alpha L_\alpha$ such that the following holds:
The leaves $L_\alpha$ all have the same dimension $k$ and to each
$m\in M$ there is a smooth local chart $(x,y):
U\ra\R^k\times\R^{n-k}$ such that the connected components of
$U\cap L_\alpha$ for each $\alpha$ are all of the form $\{ y=
const\}$.

Then the tangent bundles of the $L_\alpha$ define a subbundle
$T_F$ of the tangent bundle $TM$ which is integrable, i.e., for
any two smooth sections $X,Y$ of $T_F\subset TM$ their Lie product
$[X,Y]$ again is a section of $T_F$. The Lie product is defined
since we can view $X$ and $Y$ as vector fields on $M$. Every
integrable subbundle of $TM$ comes from a unique foliation in this
way.

Let $\CF$ be a foliation and let $\CS_\CF$ be the sheaf of germs
of smooth functions on $M$ which are locally constant on each leaf
of $\CF$. The {\it tangential cohomology}\index{tangential
cohomology} $H^.(\CF)$ of $\CF$ is by definition the sheaf
cohomology of the sheaf $\CS_\CF$. The tangential cohomology can
be computed by a de Rham complex
$$
0\ra C^\infty(M)\begin{array}{c}{d_0}\\ {\ra}\\ {}\end{array}
\Ga^\infty(T_\CF^*)\begin{array}{c}{d_1}\\ {\ra}\\ {}\end{array}
\Ga^\infty(\wedge^2T_\CF^*)\ra\dots
$$
The spaces $\Ga^\infty(\wedge^qT_\CF^*)$ are Fr\'echet spaces but
the cohomology spaces $H^q(\CF) = \ker(d^q)/\im(d^{q-1})$ are not
in general Hausdorff, since the image of $d^{q-1}$ need not be
closed. So let
$$
\bar{H}^q(\CF)\= \ker d^q/{\rm closure\ of\ (im\ }\ d^{q-1})
$$
be the {\it reduced tangential cohomology}\index{reduced
tangential cohomology}.

Let $V$ be a locally convex complex vector space, and let
$(V_i)_{i\in I}$ be a family of closed subspaces. We say that $V$
is the {\it topological direct sum}\index{topological direct sum}
of the $V_i$ if every $v\in V$ can be written in a unique way as
a convergent sum
$$
v\=\sum_{i\in I} x_i,
$$
with $x_i\in V_i$.

Now let $T:V\ra V$ be a continuous linear operator. For
$\la\in\C$ let
$$
{\rm Eig}^\infty_T(\la)\= \bigcup_{n=1}^\infty
\ker\left((T-\la{\rm Id})^n\right)
$$
be the generalized eigenspace. The operator $T$ is called {\it
traceable}\index{traceable operator} if $V$ is the topological
direct sum of the generalized eigenspaces of $V$ and
$$
\sum_{\la\ne 0} \dim({\rm Eig}^\infty_T(\la)) |\la |\ <\ \infty.
$$
In this case we define the {\it trace}\index{trace of a traceable
operator} of $T$ by
$$
\tr T\= \sum_{\la\ne 0}\dim\left({\rm Eig}^\infty_T(\la)\right)\, \la.
$$
Note that every linear operator on a finite dimensional space is
traceable and the trace in that case coincides with the usual
trace.

Let $E$ be a smooth vector bundle over $M$ with a flat leafwise
connection, i.e., a differential operator of order one
$$
\nabla : \Ga^\infty(E)\ra\Ga^\infty(T_\CF^*\otimes E)
$$
such that $\nabla(fs)=d_0(f)\otimes s+f\nabla(s)$ then $\nabla$
extends in the usual manner to
$$
\nabla : \Ga^\infty(\wedge^pT_\CF^*\otimes
E)\ra\Ga^\infty(\wedge^{p+1}T_\CF^*\otimes E)
$$
and flatness means that $\nabla^2=0$. This then defines a complex
whose reduced cohomology we will denote by $\bar{H}^q(\CF\otimes
E)$.

Back to our previous setting the representation $\tau$ defines a
homogeneous vector bundle $E_\tau$ over $G/K_M$ and by homogeneity
this pushes down to a locally homogeneous bundle over $\Ga \bs
G/K_M = {}_MX_\Ga$. The tangent bundle $T(_MX_\Ga)$ can be
described in this way as stemming from the representation of $K_M$
on
$$
\g /\k_M \cong \a \oplus \p_M \oplus \n \oplus\bar{\n}.
$$
We get a splitting into subbundles
$$
T(_MX_\Ga) = T_c \oplus T_n \oplus T_u \oplus T_s.
$$
These bundles can be characterized by dynamical properties: the
action of $A\cong\R^r$, which we will denote by $\phi$, is
furnished with a positive time direction given by the positive
Weyl chamber $A^+=\{ a\in A| a^{-1}\in A^-\}$. Then $T_s$, the
\emph{stable part} is characterized by the fact that $A^+$ acts
contractingly on $T_s$. On the \emph{unstable part} $T_u$ the
opposite chamber $A^-$ acts contractingly. $T_c$, the \emph{central
part} is spanned by the flow $A$ itself and $T_n$ is an additive
\emph{neutral part}. Note that $T_n$ vanishes if we choose $H$ to
be the maximal split torus. The bundle $T_n\oplus T_u$ is
integrable, so it defines a foliation $\CF$. To this foliation we
have the reduced tangential cohomology $\bar{H}^*(\CF)$. If
$\tau$ is a finite dimensional representation of $M$ it defines a
flat homogeneous bundle over $MN/K_M$, so a bundle with a flat
leafwise connection over $G/K_M$. We can then define the twisted
tangential cohomology $\bar H^q(\CF\otimes\tau)$ as the
cohomology of the de Rham complex for the foliation with
coefficients in the tangentially flat bundle given by $\tau$. The
flow $A$ acts on the tangential cohomology whose alternating sum
we will consider as a virtual $A$-module. For any $\ph\in
C_c^\infty(\overline{A^-})$ we define $L_\ph^q =
\int_{A^-}\ph(a)(a|\bar H^q(\CF\otimes\tau)) da$ as an operator
on $\bar H^q(\CF\otimes\tau)$.  Then we have

\begin{proposition}
Let $\tau$ be a finite dimensional representation of $M$. Under
the assumptions of Theorem \ref{lefschetz} the operator
$L_\ph^q=(L_\ph |\bar{H}^q(\CF\otimes \tau))$  is traceable  for
each $q$. The spectral side of Theorem \ref{lefschetz}
 can be written as
$$
\sum_q(-1)^q \tr (L_\ph |\bar{H}^q(\CF\otimes \tau)).
$$
\end{proposition}

\prf Follows from the Lefschetz formula. \qed

According to a general philosophy of C. Deninger \cite{den} the
local terms of a Lefschetz formula should be formed by the same
principal as the global term. We will show this on our present
case.

The torus $A$ acts on $\Ga\bs G/K_M$ by $a.(\Ga xK_M)=\Ga xaK_M$.
Let $E$ be the set of all compact $A$-orbits in $\Ga\bs G/K_M$
modulo homotopy. Note that in the rank one case $\Ga\bs G/K_M$
just gives the sphere bundle and $E$ is the set of all closed
orbits. In general, $\Ga\bs G/K_M$ is a piece of the sphere bundle
associated with $A$. For each $e\in E$ let $X_e$ be the union of
all $A$-orbits in the class $e$. It is known that $X_e$ is a
smooth submanifold of $\Ga\bs G/K_M$.

\begin{theorem}
In the above notation we have
$$
\sum_{q,r} (-1)^{q+r+\rank\CF}\
 \tr(L_\ph |\bar{H}^q(\CF\otimes\wedge^rT_s))\hspace{70pt}
$$ $$
\hspace{70pt}= \sum_{e\in E} \sum_{q,r}(-1)^{q+r+\rank\CF|_{X_e}}\
 \tr(L_\ph |\bar{H}^q(\CF|_{X_e}\otimes\wedge^r(TX_e)_s)).
$$
\end{theorem}
Note that the stable bundle of $X_e$ is the zero bundle. Therefore the sum over $r$ on the right hand side
has only one summand, namely for $r=0$.

\prf The global side of this theorem coincides with the global
side of Corollary \ref{ruelle-lef} times $(-1)^{\dim
N+\dim\p_M}=(-1)^{\dim N}$.

For the local side consider the following equivalence relation on
the set of $\Ga$-conjugacy classes $[\ga]$. We say that two
classes $[\ga]$ and $[\eta]$ are equivalent if they have
representatives $\ga$ and $\eta$ satisfying $G_\ga=G_\eta$. Let
$E_P(\Ga)$ be the set of equivalence classes of $[\ga]$ which lie
in $\CE_P(\Ga)$. Then $E_P(\Ga)$ coincides with the set $E$ above.
It emerges that for $[\ga]\in\CE_P(\Ga)$ the quantities $\la_\ga$
and $\chi(A_\ga\bs X_\ga)$ only depend on the class $e=e_\ga$ in
$E_P(\Ga)$. Therefore we will write these as $\la_e$ and
$\chi(A_e\bs X_e)$. Let furthermore $\Ga_e$ denote the lattice
$\Ga_{A_\ga}$ of $A_\ga$, where we have chosen a fixed
representative $\ga$. The local side of the formula of Corollary
\ref{ruelle-lef} can be written as
$$
\sum_{e\in E_P(\Ga)} \la_e\chi(A_e\bs X_e) \sum_{[\ga]\in
e}\ph(a_\ga)
$$
Now the form $B$ defines a bilinear pairing $(.,.):A_\ga\times
A_\ga\ra \R^\times_+$. Let $[\ga]\in e$ then $e$ consists of all
classes $[\sigma]$ with $\sigma\in\Ga_\ga\cap A^-M_{ell}$. Since
$\Ga$ is neat it follows that the projection
$\Ga_\ga\ra\Ga_{A_\ga}$ is an isomorphism. We will identify these
two groups in the sequel. Since $\ph$ vanishes on non-regular
elements we may replace the sum over $e$ by a sum over $\Ga_\ga$
or $\Ga_{A_\ga}$. Let $\Ga_\ga^\perp$ be the lattice in $A_\ga$
which is dual to $\Ga_\ga$. The Poisson summation formula tells us
that
\begin{eqnarray*}
\la_e\sum_{[\ga]\in e}\ph(a_\ga) &=& \la_\ga
\sum_{\sigma\in\Ga_\ga}\ph(a_\sigma)\\
 &=& \sum_{\eta\in\Ga_\ga^\perp} \hat{\ph}(a_\eta),
 \end{eqnarray*}
 where
 $$
 \hat{\ph}(a)\= \int_A\ph(y)e^{-2\pi i(a,y)}dy
 $$
 is the Fourier transform. We get that the local side equals
 $$
 \sum_e \chi(A_e\bs X_e) \sum_{\eta\in\Ga_e^\perp} \int_A\ph(y)e^{-2\pi
 i(\eta,y)}dy.
 $$
\begin{eqnarray*}
&=& \sum_e\chi(A_e\bs X_e) \sum_{\eta\in\Ga_e^\perp}\langle
L_\ph\epsilon_\eta,\epsilon_\eta\rangle\\
 &=& \sum_e\chi(A_e\bs X_e) \tr(L_\ph|C^\infty (\Ga_e\bs A_e)),
 \end{eqnarray*}
where we have written $\epsilon_\eta(a)=e^{2\pi i(\eta,a)}$. Note
that the $\epsilon_\eta$ form an orthonormal basis of
$L^2(\Ga_e\bs A_e)$. The trace is to be interpreted as the trace
of the nuclear operator $L_\ph$.

Now note that the foliation $\CF$ restricted to
$$
X_\ga =\Ga_\ga\bs G_\ga/K_\ga = \Ga_\ga A_\ga M_\ga /K_{M_\ga}
$$
has leaves $\Ga_\ga aM_\ga$ parametrized by $ \Ga_\ga\bs A_\ga$.
The space of leaves is Hausdorff and coincides with $ \Ga_\ga\bs
A_\ga$. The de Rham complex computing the tangential cohomology
coincides with the smooth family of full de Rham complexes of the
leaves. The cohomology of each leaf is the cohomology of
$M_\ga/K_{M_\ga}$. Therefore
$$
\chi(A_\ga\bs X_\ga)\tr(L_\ph|C^\infty(\Ga_\ga\bs A_\ga))
$$
equals
$$
\sum_q(-1)^q\tr(L_\ph|\bar{H}^q(\CF|_{X_\ga}))
$$
and the theorem follows. \qed

\chapter{Meromorphic continuation}
In this chapter $G$ will be a connected semisimple Lie
group with finite center. Fix a maximal compact subgroup
$K$ with Cartan involution $\theta$.

\section{The zeta function} \label{constzeta}

Fix a $\theta$-stable Cartan subgroup $H$ of split rank 1. Note
that such a $H$ doesn't always exist. It exists only if the
absolute ranks of $G$ and $K$ satisfy the relation:
$$
{\rm rank}\ G -{\rm rank}\ K \leq 1.
$$
This certainly holds if $G$ has a compact Cartan subgroup or if
the real rank of $G$ is one. In the case ${\rm rank}\ G ={\rm
rank}\ K$, i.e., if $G$ has a compact Cartan there will in general
be several $G$-conjugacy classes of split rank one Cartan
subgroups. In the case ${\rm rank}\ G -{\rm rank}\ K = 1$,
however, there will only be one. The number $FR(G):={\rm rank}\ G
-{\rm rank}\ K $ is called the {\it fundamental
rank}\index{fundamental rank} of $G$.

Write $H=AB$ where $A$ is the connected split component and
$B\subset K$ is compact. Choose a parabolic subgroup $P$
with Langlands decomposition $P=MAN$. Then $K_M =K\cap M$
is a maximal compact subgroup of $M$. Let
$\bar{P}=MA\bar{N}$ be the opposite parabolic. As in the
last section let $A^-=\exp(\a_0^-)\subset A$ be the
negative Weyl chamber of all $a\in A$ which act
contractingly on $\n$. Fix a finite dimensional
representation $(\tau ,V_\tau)$ of $K_M$.

Let $H_1\in \a_0^-$ be the unique element with $B(H_1)=1$.

Let $\Ga\subset G$ be a cocompact discrete subgroup which is
 neat. An element $\ga \neq 1$ will be called {\it
primitive}\index{primitive element} if $\tau \in \Ga$ and $\tau^n
=\ga$ with $n\in \N$ implies $n=1$. Every $\ga \neq 1$ is a power
of a unique primitive element. Obviously primitivity is a property
of conjugacy classes. Let $\CE_P^p(\Ga)$ denote the subset of
$\CE_P(\Ga)$ consisting of all primitive classes. Recall the
length $l_\ga$ of any geodesic in the class $[\ga]$. If $\ga$ is
conjugate to $am\in A^-M_{ell}$ then $l_\ga=l_a$, where $l_a=|\log
a|=\sqrt{B(\log a)}$. Let $A_\ga$ be the connected split component
of $G_\ga$, then $A_\ga$ is conjugate to $A$.

In this section we are going to prove the following

\begin{theorem}\label{genSelberg} Suppose that $\tau$ is in
the image of the restriction map $res_{K_M}^M$ or that $M$
is orientation preserving. Let $\Ga$ be  neat and
$(\omega ,V_\omega)$ a finite dimensional unitary
representation of $\Ga$. For $\Re(s)>>0$ define the {\it
generalized Selberg zeta function}:\index{Selberg zeta
function}
$$
Z_{P,\tau,\omega}(s) \= \prod_{[\ga]\in {\cal E}_P^p(\Ga)}
\prod_{N\geq 0} \det\left(1-e^{-sl_\ga}\ga \left|
\begin{array}{c}V_\omega \otimes V_\tau \otimes
S^N(\n)\end{array}\right.\right)^{\chi(A_\ga\bs X_\ga)},
$$
where $S^N(\n)$ denotes the $N$-th symmetric power of the space
$\n$ and $\ga$ acts on $V_\omega \otimes V_\tau \otimes S^N({\n})$
via $\omega(\ga) \otimes \tau(b_\ga) \otimes Ad^N(a_\ga b_\ga)$,
here $\ga \in \Ga$ is conjugate to $a_\ga b_\ga \in A^-B$.

Then $Z_{P,\tau,\omega}$ has a meromorphic continuation to the
entire plane. The vanishing order of $Z_{P,\tau ,\omega}(s)$ at a
point $s=\la (H_1)$, $\la \in \a^*$, is
$$
(-1)^{\dim N}\sum_{\pi \in \hat{G}}N_{\Ga ,\omega}(\pi)
\sum_{p,q}(-1)^{p+q} \dim
\left(\begin{array}{c}H^q(\n,\pi_K)^\la\otimes \wedge^p\p_M \otimes
V_{\breve{\tau}}\end{array}\right)^{K_M}.
$$
Further, all poles and zeroes of the function $Z_{P,\tau
,\omega}(s+|\rho_0|)$ lie in $\R \cup i\R$.
\end{theorem}

The proof will rely on the following

\begin{proposition} \label{existf}
For $j\in\N$ and $s\in\C$ with $j,\Re(s)>>0$ there is a function
$f$ on $G$ which is $j-(2\dim N +\dim K)$-times continuously
differentiable, satisfies the conditions of Proposition
\ref{convergence-of-tf} ( so the trace formula is valid for $f$)
and for $g\in G$ semisimple we have that $\CO_g(f)=0$ whenever $g$
is not conjugate to an element in $A^-B$ and for $ab\in A^-B$ it
holds
$$
\CO_{ab}(f)\= \O_{b}^M(f_\tau^M)\ \frac{l_a^{j+1}e^{-sl_a}}{\det
(1- a b|\n)}.
$$
Especially it follows that
$$
\sum_{\pi \in \hat{G}} N_{\Ga,\omega} (\pi)\ \tr\pi(f)
$$
equals
$$
 \sum_{[\ga]} \vol (\Ga_\ga \bs G_\ga )\ \tr\omega(\ga)\
\O_{b_\ga}^M(f_\tau^M)\ \frac{l_\ga^{j+1}e^{-sl_\ga}}{\det (1-
a_\ga b_\ga|\n)},
$$
where the sum runs over all classes $[\ga]$ such that $\ga$ is
conjugate to an element $a_\ga b_\ga$ of $A^- B$ and the
convergence is locally uniform in $s$.
\end{proposition}

\prf Let $C=\supp(f_\tau^M)\cup B$, then $C$ is a compact subset
of $G$. Thus there exists $a_0\in A^-$ such that $a_0
\overline{A^-}C$ is disjoint to $S=\partial(AM)^\sim$.

\begin{lemma}\label{exist-Q}
There exists a $(j-\dim N)$-times continuously differentiable,
conjugation invariant function $Q$ on $AM$ such that
\begin{itemize}
\item
$\supp\ Q\subset $ closure of $(AM)^\sim$ and
\item
for $am\in A^-M_{ell}\cup a_0\overline{A^-}C$ we have
$$
Q(am)\= \frac{l_a^{j+1}}{\det(1-am|\n)}
$$
\end{itemize}
\end{lemma}

\prf Let $Z$ be the center of $G$ and let $G'$ be the quotient
group $G/Z$. For any subset $H$ let $H'$ be its image in $G'$.
Likewise, for any element $x\in G$ let $x'$ be its image in $G'$.
Suppose we have a function $Q'$ on $(AM)'=A'M'$ satisfying the
above conditions with $(AM)^\sim$ and $C_1$ replaced with their
images in $G'$. Then the function
$$
Q(am)\= Q'(a'm')
$$
satisfies the claim. We are thus reduced to the case that $G$ has
trivial center. We will assume that for the length of this proof.
Then $G$ coincides with the image under the adjoint representation
and so there is a linear algebraic group $\CG$ over $\R$ such that
$G$ is isomorphic to the connected component of $\CG(\R)$. For the
parabolic $P=MAN$ it follows that there is a Zariski connected
linear algebraic subgroup $\CM$ of $\CG$ such that $M=\CM(\R)\cap
G$ (see \cite{Borel-lingroups} Cor. 11.12).

Let $\CM(\C)$ be the complex points of $\CM$. Then $\CM(\C)$ is
connected in the Hausdorff topology and the exponential map $\exp
: \m\ra \CM(\C)$ is surjective as follows from the Jordan
decomposition \cite{Borel-lingroups} I.4, since this implies that
any element $m$ of $\CM(\C)$ can be written as a product
$m=m_sm_u=m_um_s$ where $m_s$ is semisimple and $m_u$ is
unipotent. Semisimple and unipotent elements are in the image of
the exponential map and, since they commute, also their product.
This implies the surjectivity of the exponential map.

The reductive algebraic group $\CM$ acts on itself by conjugation
and the quotient $\CM/conj$ is a smooth affine variety over $\C$.
The set of complex points $\CM/conj(\C)$ therefore is a smooth
manifold. We will show

\begin{lemma} \label{exist-psi}
There is a smooth function $\psi$ on $A^-\times \CM/conj(\C)$ or a
conjugation invariant smooth function $\psi$ on $A^-\times
\CM(\C)$ such that
\begin{itemize}
\item
$\psi \equiv 1$ on $A^- M_{ell}\cup a_0\overline{A^-}C$ and
\item
$\supp \psi\cap A^-M\subset (AM)^\sim$ and
\item
the function $am\mapsto l_a^{j+1}\psi(am)$ extends to a $j$-times
continuously differentiable function on $AM$ which is zero on
$A^+M$.
\end{itemize}
\end{lemma}

\prf First note that the image of $M_{ell}$ in $\CM/conj(\C)$
coincides with the image on $B$ and therefore is compact. The set
$(AM)^\sim$ satisfies $\exp(X)\in (AM)^\sim \Rightarrow
\exp(tX)\in(AM)^\sim$ for any $t>0$, therefore $(AM)^\sim$ is an
``open cone''. In the following picture the $x$-axis represents
$M$ and the $y$-axis is $A$:

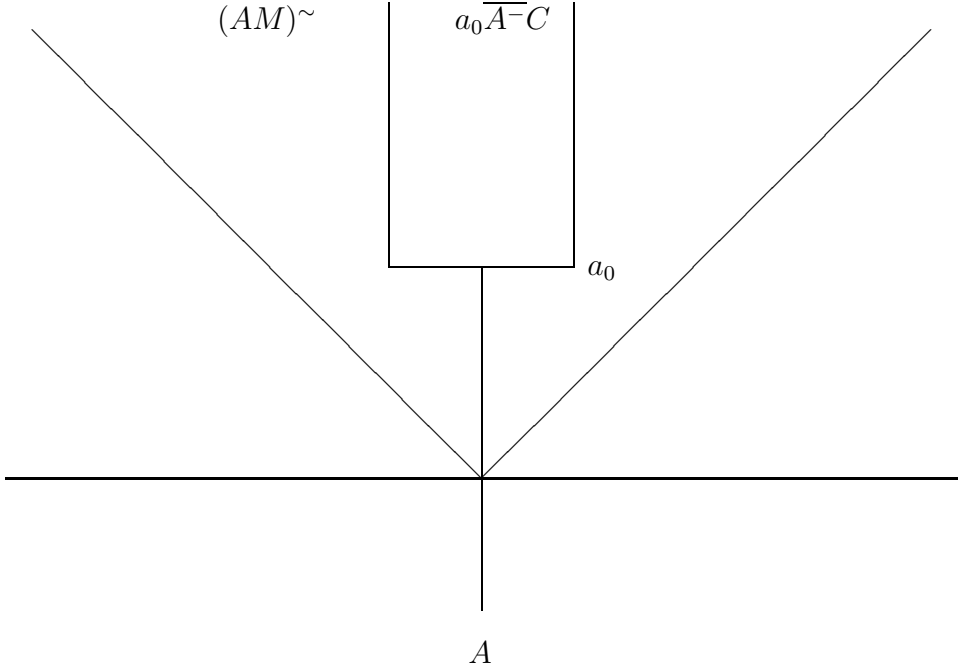
\begin{picture}(300,300)(-150,-120)
\put(-180,0){\line(1,0){360}} \put(0,-50){\line(0,1){130}}
\put(-35,80){\line(1,0){70}}  \put(0,0){\line(1,1){170}}
\put(0,0){\line(-1,1){170}} \put(-35,80){\line(0,1){100}}
\put(35,80){\line(0,1){100}} \put(-10,170){$a_0\overline{A^-}C$}
\put(-100,170){$(AM)^\sim$} \put(40,77){$a_0$} \put(-5,-70){$A$}
\end{picture}

For any $a_1\in A^-$ it is easy to give a function as in Lemma
\ref{exist-psi} on $a_1A^-$ since there the last property of
$\psi$, as insisted in the lemma, is immediate. A partition of
unity argument then shows that we only have to worry about small
values of $a$, and for those the set $a_0\overline{A^-}C$ doesn't
play a role.

Choose a Riemannian metric on $\CM/conj(\C)$ then, by the
compactness of $M_{ell}/conj$ there is a $c>0$ such that
$$
{\rm dist}(aM_{ell},\partial(AM)^\sim)\ \ge \ c\ l_a,
$$
from which the Lemma \ref{exist-psi} is deduced by means of the
following

\begin{lemma}\label{exist-psi2}
Let $Y$ be a smooth Riemannian manifold and let $T\subset Y$ be a
compact subset. Let $S\subset \R^+\times Y$ be closed and assume
that there is $c>0$ such that for any $t>0$
$$
\dist(T,S_t)\ \ge\ ct,
$$
where $S_t=S\cap(\{ t\} \times Y)$.

Then there is a smooth function $\psi : \R^+\times Y\ra [0,1]$
such that $\psi \equiv 1$ on $\R^+\times T$ and $\psi \equiv 0$ on
$S$ and the function $(t,y)\mapsto t^{j+1}\psi(t,y)$ extends to a
$j$-times continuously differentiable function on $\R\times Y$
which is zero on $\R^-\times Y$.
\end{lemma}

\prf At first we show that is suffices to consider the case when
$T$ lies in a chart neighborhood of $Y$. For this let
$C=C_1\cup\dots\cup C_N$ be a decomposition into compact subsets
where each of the $C_j$ lies in a chart neighborhood of $Y$.
Assume the lemma proven for each $C_j$, so for each $j$ there
exists a function $\psi_j$ satisfying the claim of the lemma for
$C_j$. Let $\chi :\R\ra [0,1]$ be smooth, $\chi \equiv 0$ on
$\R^-$ and $\chi \equiv 1$ on $[1,\infty)$. Let
$$
\psi\= \chi(\psi_1+\dots +\psi_N),
$$
then $\psi$ satisfies the claims of the lemma.

Now we are reduced to the case when $T$ lies in  a chart
neighborhood of $Y$. We can replace $Y$ by that neighborhood and
assume $Y=\R^n$. By compactness of $T$ we max also assume that the
metric of $Y$ is the euclidean one. Let
$$
 \tilde\psi_t(y)\= \left\{ \begin{array}{cl} 1 & {\rm
 if\ dist}(y,T)\le \frac{c}{4},\\ 0& {\rm otherwise.}
 \end{array}\right.
$$
Let $h:\R^n\ra [0,\infty)$ be smooth, with support in the ball
with radius one around $0$ and such that $\int_{\R^n} h(y) dy =1$.
For $t>0$ the function $h_t(y) = h(4y/ct)$ has support in the ball
with radius $ct/4$ around zero. Let
$$
\psi(t,y)\= \tilde{\psi}_t*h_t(y),
$$
where $*$ means convolution in $\R^n$. Then $\psi$ satisfies the
claims of  the lemma. Lemma \ref{exist-psi2} and thus Lemma
\ref{exist-psi} are proved. \qed

\pagebreak Finally, the function $Q$ of Lemma \ref{exist-Q} then
is gotten by setting
$$
Q(am)\= \frac{l_a^{j+1}\psi(am)}{\det(1-am|\n)}
$$
\qed

Now the group $G$ may have nontrivial center again. In order to
continue the proof of Proposition \ref{existf} choose any smooth
$\eta : N \ra \R$ which has compact support, is non-negative,
invariant under $K\cap M$ and such that $\int_N \eta(n) dn =1$.

Given these data we define $f =f_{\eta, \tau ,j,s} : G\ra \C$  by
$$
f(kn ma (kn)^{-1}) \= \eta(n) f_\tau^M(m) Q(am)e^{-sl_a},
$$
for $k\in K$, $n\in N$, $am\in (AM)^{\sim}$. Further $f(x)=0$ if
$x$ is not in $G.(MA)^{\sim}$.

The well-definedness of $f$ is proved in the same way as in Lemma
\ref{welldef}. We have to show that the function $f$ for
$\Re(s)>>0$ and $j>>0$ satisfies the conditions of Proposition
\ref{convergence-of-tf}. For this consider the map
$$
\begin{array}{cccc}
F\ :& K\times N\times M\times A &\ra& G\\
{}& (k,n,m,a)&\mapsto&knam(kn)^{-1}.
\end{array}
$$
Now $f$ is a $j-\dim(\n)$-times continuously differentiable
function on $K\times N\times M\times A$ which factors over $F$.
To compute the order of differentiability as a function on $G$ we
have to take into count the zeroes of the differential of $F$. So
we compute the differential of $F$. Let at first $X\in \k$, then
$$
DF(X) f(knam(kn)^{-1}) = \frac{d}{dt}|_{t=0} f(k
\exp(tx)namn^{-1}\exp(-tX)k^{-1}),
$$
which implies the equality
$$
DF(X)_x = (\Ad(k(\Ad(n(am)^{-1}n^{-1})-1)X)_x,
$$
when $x$
equals $knam(kn)^{-1}$. Similarly for $X\in\n$ we get that
$$
DF(X)_x = (\Ad(kn)(\Ad((am)^{-1})-1)X)_x
$$
and for $X\in
\a\oplus \m$ we finally have $DF(X)_x = (\Ad(kn)X)_x$.
From this it becomes clear that $F$, regular on $K\times
N\times (MA)^\sim$, may on the boundary have vanishing
differential of order $\dim(\n)+\dim(k)$. Together we get
the claimed degree of differentiability. So we assume
$j\ge 2\dim(\n)+\dim(\k)$ from now on. In order to show
that $f$ goes into the trace formula for $s$ and $j$ large
we fix $N>\frac{\dim G}{2}$. We have to show

\hspace{10pt}
$Df\in L^1(G)$ for any $D\in U(\g)$ of degree $\le 2N$.

For this we recall the map $F$ and our
computation of its differential. Let $\q\subset \k$ be a
complementary space to $\k_M$. On the regular set $DF$ is
surjective and it becomes bijective on
$\q\oplus\n\oplus\a\oplus\m$. Fix $x=knam(kn)^{-1}$ in the
regular set and let $DF_x^{-1}$ denote the inverse map of $DF$
which maps to $\q\oplus\n\oplus\a\oplus\m$. Introducing norms on
the Lie algebras we get an operator norm for $DF_x^{-1}$ and the
above calculations show that $\norm{DF_x}\le P(am)$, where $P$ is
a class function on $M$, which, restricted to any Cartan $H=AB$
of $AM$ is a linear combination of quasi-characters. Supposing
$j$ and $\Re(s)$ large enough we get for $D\in U(\g)$ with
$\deg(D)\le 2N$:
$$
|Df(knam(kn)^{-1})| \le \sum_{D_1} P_{D_1}(am) |D_1f(k,n,a,m)|,
$$
where the sum runs over a finite set of $D_1\in
U(\k\oplus\n\oplus\a\oplus\m)$ of degree $\le 2N$ and $P_{D_1}$
is a function of the type of $P$. On the right hand side we have
considered $f$ as a function on $K\times N\times A\times M$.

We want to show  that if $j,\Re(s)>>0$ then $(D_1f) P P_{D_1}$
lies in $L^1(K\times N\times A\times M)$ for any $D_1$. The only
difficulty in this is the derivatives of the function $Q$ which we
know only little about. But fortunately we have arranged things so
that we do not have to worry about those. At first it is clear
that we only have to care about what happens outside a compact
set, so we only have to consider the points $x=knam(kn)^{-1}$ with
$a\in a_0\overline{A^-}$ and $m\in\supp(f_\tau^M)$. In this range
we have $Q(am)=\frac{l_a^{j+1}}{\det(1-am|\n)}$ and the
derivatives are rational functions of quasicharacters on $A$ and
will be in $L^1(G)$ if only $\Re(s)$ is large enough. This implies
a). We have proven that the trace formula is valid for $f$.

We have to compute the orbital integrals to see that they give
what is claimed in the proposition. Since $f$ vanishes on all
$g\in G$ which are not conjugate to some $am\in A^-M$ it remains
to compute the orbital integrals $\CO_{am}(f)$. The group $\Ga$,
being cocompact, only contains semisimple elements \cite{marg}
therefore it suffices to compute $\CO_{am}(f)$ for semisimple $m$.
We have the decomposition $G=KNAM$ and, in this order, the
Haar-measure of $G$ is just the product of the Haar measures of
the smaller subgroups. Since $f(am)=0$ unless $a$  is regular we
restrict to $am \in A^-M$, $m$ semisimple. At first we consider
the case, when $m$ is not elliptic and we want to show that the
orbital integral $\CO_{am}(f)$ vanishes in this case. We show more
sharply that the integral
$$
\int_{G/AM_m}f(xamx^{-1})\ dx
$$
vanishes whenever $f$ is nontrivial on the orbit of $am$. Since
$AM_m$ is a subgroup of $G_{am}$ this implies the assertion. The
latter integral equals
$$
\int_{KN} \int_{M/M_m} f(kn a m'm(m')^{-1}(kn)^{-1})\ dm' dk dn,
$$
which by definition of $f$ equals
$$
\CO_m^M(f_\tau^M) Q(am)e^{-sl_a}.
$$
Thus we see that $\CO_{am}(f)=0$ unless $m$ is elliptic since this
holds for the $M$-orbital integrals over $f_\tau^M$. Now we
consider $am\in A^-M$ with $m$ elliptic. In this case we have the
equality
$$
G_{am}\= AM_m.
$$
So that as above we get
$$
\CO_{am}(f) \= \CO_m^M(f_\tau^M) \frac{l_a^{j+1}e^{-sl_a}}
                       {\det(1-am|\n)}.
$$
The proposition follows.
\qed

Now to the proof of Theorem \ref{genSelberg}. Let
$\ga\in\CE_P(\Ga)$ be conjugate to $a_\ga b_\ga\in A^-B$. By
Theorem \ref{3.3} we have
$$
\chi(A_\ga\bs X_\ga)\=
\frac{|W(\m,\t)|\prod_{\alpha\in\Phi^+(\m,\t)}(\rho_M,\alpha)}
     {c_Ml_\ga}
\vol(\Ga_\ga\bs G_\ga).
$$
By Theorem \ref{orbitalint} we have
$$
\CO_{b_\ga}(f_\tau^M) \= \tr\tau(b_\ga)
\frac{|W(\m,\t)|\prod_{\alpha\in\Phi^+(\m,\t)}(\rho_M,\alpha)}
     {c_M}l_a.
$$
Therefore
$$
\CO_{b_\ga}(f_\tau^M)\vol(\Ga_\ga\bs G_\ga)\= \chi(A_\ga\bs
X_\ga)\tr\tau(b_\ga).
$$
Since $\Ga$ is  neat it follows  $X_{\ga^n} \cong X_{\ga}$ by Lemma \ref{neatXgamma}
for any $n$. This implies that  the geometric side of the trace
formula:
$$
\sum_{[\ga]} \frac{l_{\ga_0} \chi_{_1}(X_\ga) {\tr\ \tau(b_\ga)}
        l_\ga^{j+1} e^{-sl_\ga}\tr(\omega(\ga))}
     {\det(1-a_\ga b_\ga | \n)}
$$
equals
$$
(-1)^{j+1} (\der{s})^{j+2} \log Z_{P,\tau,\omega} (s).
$$

We now consider the spectral side.  For $\pi \in \hat{G}$ we get
$$
\tr\ \pi(f) \= \int_G \Theta_\pi^G(x) f(x)\ dx
$$
and this equals
$$
\int_{MA^-} \Theta_\pi^G(ma)\ f_\tau^M(m)\
 \left| {\det(1-ma|{\n}\oplus \bar{\n})}\right| \ Q(am)e^{-sl_a}\ da dm.
$$
Repeating the computations of the last section we get
$$
\tr\ \pi(f) \=  \sum_{p=0}^{\dim N}(-1)^{p+\dim N}\int_{MA^-}
f_{\tau}^M(m)\ \Theta_{H^p(\n, \pi_K)}^{MA}(ma)\ dm
l_a^{j+1}e^{-sl_a}\ da.
$$
Using the properties of $f_\tau^M$ we get that $ \tr\ \pi (f)$
equals
$$
\sum_{p,q\ge 0}(-1)^{p+q+\dim N}\int_{A^-} \tr (a|(H^p(\n,\pi_K)
\otimes \wedge^q \p_M \otimes V_{\tau})^{K_M})\ l_a^{j+1}
e^{-sl_a} \ da
$$
Now let the $A$-decomposition of the virtual $A$-module
$(H^*(\n,\pi_K) \otimes \wedge^* \p_M \otimes V_{\tau})^{K_M}$ be
$\bigoplus_{\la\in\La_\pi} E_\la,$ where $E_\la$ is the biggest
virtual subspace on which $a-e^{\la(\log a)}$ acts nilpotently.

Let $m_\la\in\Z$ denote the virtual dimension of $E_\la$.
Integrating over $\a^-$ instead over $A^-$ we get
\begin{eqnarray*}
\tr\ \pi (\Phi) &\=& (-1)^{\dim N}\int_0^\infty \sum_\la m_\la
e^{(\la(H_1)-s)t} t^{j+1}\ dt
\\
        &\=& (-1)^{j+1+\dim N} (\der{s})^{j+1} \sum_\la m_\la \rez{s-\la(H_1)}
\end{eqnarray*}

From this it follows that $Z_{P,\tau ,\omega}(s)$ extends to a
meromorphic function on the plane and that the order of $Z_{P,\tau
,\omega}(s)$ at a point $s=\la (H_1)$, $\la \in \a^*$, is
$$
(-1)^{\dim N} \sum_{\pi \in \hat{G}}N_{\Ga ,\omega}(\pi)
\sum_{p,q}(-1)^{p+q} \dim \left(\begin{array}{c}H^q(\n,\pi_K)^\la
\otimes \bigwedge^p\p_M \otimes
V_{\tau}\end{array}\right)^{K_M}.
$$

To prove the last assertion on the zeroes and poles recall that in
the tensor product $H^q(\n,\pi^0)\otimes \wedge^p\p_M \otimes
V_{\tau}$ the group $A$ only acts on the first tensor factor. Let
$\la\in\a^*$ be such that
$$
\sum_p (-1)^p\dim \left( H^q(\n ,\pi_K)^\la\otimes \wedge^p\p_M
\otimes V_\tau\right)^{K_M}\ne 0.
$$
Then there is a $M$-irreducible subquotient $\eta$ of $H^q(\n
,\pi_K)$ satisfying the same assertion, especially we get
$\tr\eta(f_\tau^M)\ne 0$. A {\it standard
representation}\index{standard representation} $\xi$ of $M$ is a
representation induced from a tempered one. Any character
$\Theta_\xi$ for $\xi$ irreducible admissible is an integer linear
combination of standard representations. On the other hand by
Lemma \ref{pivonggleichnull} we have that $\tr\xi(f_\tau^M)=0$ for
$\xi$ induced, so that there is a discrete series or limit of
discrete series representation $\xi$ having the same infinitesimal
character $\la_\eta$ as $\eta$. Especially it follows that
$\la_\eta$ is real.

By Corollary 3.32 of \cite{HeSch} it follows that $\la +\la_\eta
=w\la_\pi +\rho-\rho_M$, where $\la_\pi$ is the infinitesimal
character of $\pi$ and $w$ is an element of the Weyl group $W(\g
,\h)$. Let $C$ be the Casimir operator of $G$, then $\pi(C)$ is a
real number. We compute $\pi(C) = B(\la_\pi) -B(\rho) =
B(w\la_\pi) -B(\rho) = B(\la +\la_\eta -\rho +\rho_M) -B(\rho) =
B(\la -\rho_0) + B(\la_\eta -\rho |_M+\rho_M)-B(\rho)$. Since
$\la_\eta$ is real it follows that $B(\la_\eta -\rho |_M+\rho_M)$
is real and so then is $B(\la -\rho_0)$. Since $\a$ is one
dimensional, $\la-\rho_0$ must be real or purely imaginary. \qed

\section{The functional equation}

Assume that the Weyl group $W(G,A)$ is nontrivial, then is has
order two. To give the reader a feeling of this condition consider
the case $G=SL_3(\R)$. In that case the Weyl group $W(G,A)$ is
trivial. On the other hand, consider the case when the fundamental
rank of $G$ is $0$; this is the most interesting case to us since
only then we have several conjugacy classes of splitrank-one
Cartan subgroups. In that case it follows that the dimension of
all irreducible factors of the symmetric space $X=G/K$ is even,
hence the point-reflection at the point $eK$ is in the connected
component of the group of isometries of $X$. This reflection can
be thought of as an element of $K$ which induces a nontrivial
element of the Weyl group $W(G,A)$. So we see that in the most
important case we have $|W(G,A)|=2$.

Let $w$ be the nontrivial element. It has a representative in $K$
which we also denote by $w$. Then $wK_Mw^{-1}=K_M$ and we let
$\tau^w$ be the representation given by
$\tau^w(k)=\tau(wkw^{-1})$. It is clear from the definitions that
$$
Z_{P,\tau,\omega}\= Z_{\bar{P},\tau^w,\omega}.
$$
We will show a functional equation for $Z_{P,\tau,\omega}$ this
needs some preparation. Assume $G$ admits a compact Cartan
$T\subset K$, then a representation $\pi\in\hat{G}$ is called
\emph{elliptic} if $\Theta_\pi$ is nonzero on the compact Cartan. Let
$\hat{G}_{ell}$ be the set of elliptic elements in $\hat{G}$ and
denote by $\hat{G}_{ds}$ the subset of discrete series
representations. Further let $\hat{G}_{lds}$ denote the set of all
discrete series and all limits of discrete series representations.
In Theorem \ref{genSelberg} we have shown that the vanishing order
of $Z_{P,\tau,\omega}(s)$ at the point $s=\mu (H_1)$, $\mu\in\a^*$
is
$$
(-1)^{\dim N}\sum_{\pi\in\hat{G}} N_{\Ga ,\omega}(\pi) m(\pi ,\tau
,\mu),
$$
where
$$
m(\pi ,\tau ,\mu) \= \sum_{p,q} (-1)^{p+q} \dim
\left(\begin{array}{c}H^q(\n ,\pi_K)(\mu)\otimes \wedge^p\p_M\otimes
V_{\breve{\tau}}\end{array}\right)^{K_M}.
$$

Any character $\Theta_\pi$ for $\pi\in\hat{G}$ is an integer
linear combination of characters of standard representations. From
this it follows that for $\pi\in\hat{G}$ the character restricted
to the compact Cartan $T$ is
$$
\Theta_\pi \mid_T \= \sum_{\pi'\in\hat{G}_{lds}} k_{\pi ,\pi'}
\Theta_{\pi'} \mid_T,
$$
with integer coefficients $k_{\pi ,\pi'}$.

\begin{lemma}
There is a $C>0$ such that for $\Re (\mu(H_1))<-C$ the order of
$Z_{P,\tau,\omega}(s)$ at $s=\mu(H_1)$ is
$$
(-1)^{\dim N}\sum_{\pi\in\hat{G}_{ell}} N_{\Ga ,\omega}(\pi)
\sum_{\pi'\in\hat{G}_{lds}} k_{\pi ,\pi'} m(\pi' ,\tau ,\mu).
$$
\end{lemma}

\prf For any $\pi\in\hat{G}$ we know that if $\Theta_\pi
|_{AB_M}\neq 0$ then in the representation of $\Theta_\pi$ as
linear combination of standard characters there must occur
lds-characters and characters of representations $\pi_{\xi ,\nu}$
induced from $P$. Since $\Theta_{\pi_{\xi
,\nu}}=\Theta_{\pi_{^w\xi ,-\nu}}$ any contribution of $\pi_{\xi
,\nu}$ for $\Re (s) <<0$ would also give a pole or zero of
$Z_{P,\tau^w,\omega}$ for $\Re(s)>>0$. In the latter region we do
have an Euler product, hence there are no poles or zeroes. \qed

Now consider the case $FR(G)=1$, so there is no compact Cartan,
hence no discrete series. By a finite genus argument it follows
that

\begin{theorem} Assume that the fundamental rank of $G$ is $1$,
then there is a polynomial $P$ of degree $\leq \dim G+\dim N$ such
that
$$
Z_{P,\tau,\omega}(s) \= e^{P(s)} Z_{P,\tau^w,\omega}(2| \rho_0|
-s).
$$
\end{theorem}
\qed

Now assume $FR(G)=0$ so there is a compact Cartan subgroup $T$. As
Haar measure on $G$ we take the Euler-Poincar\'e measure. The sum
in the lemma can be rearranged to
$$
(-1)^{\dim N} \sum_{\pi'\in\hat{G}_{lds}}  m(\pi' ,\tau
,\mu)\sum_{\pi\in\hat{G}_{ell}} N_{\Ga ,\omega}(\pi) k_{\pi
,\pi'}.
$$
We want to show that the summands with $\pi'$ in the limit of the
discrete series add up to zero. For this suppose $\pi'$ and
$\pi''$ are distinct and belong to the limit of the discrete
series. Assume further that their Harish-Chandra parameters agree.
By the Paley-Wiener theorem \cite{CloDel} there is a smooth
compactly supported function $f_{\pi' ,\pi''}$ such that for any
tempered $\pi\in\hat{G}$:
$$
\tr \pi(f_{\pi' ,\pi''}) \= \left\{ \begin{array}{cl} 1& {\rm if}\
\pi =\pi'\\
                        -1  &{\rm if}\ \pi =\pi''\\
                        0   & {\rm else.}
                    \end{array}\right.
$$
Plugging $f_{\pi' ,\pi''}$ into the trace formula it follows
$$
\sum_{\pi\in\hat{G}_{ell}} N_{\Ga ,\omega}(\pi)k_{\pi ,\pi'}
\=\sum_{\pi\in\hat{G}_{ell}} N_{\Ga ,\omega}(\pi)k_{\pi ,\pi''},
$$
so that in the above sum the summands to $\pi'$ and $\pi''$ occur
with the same coefficient. Let $\pi_0$ be the induced
representation whose character is the sum of the characters of the
$\pi''$, where $\pi''$ varies over all lds-representations with
the same Harish-Chandra parameter as $\pi'$. Then for $\Re
(\mu(H_1))<-C$ we have $m(\pi_0,\tau,\mu)=0$. Thus it follows that
the contribution of the limit series vanishes.

Plugging the pseudo-coefficients \cite{Lab} of the discrete series
representations into the trace formula gives for
$\pi\in\hat{G}_{ds}$:
$$
\sum_{\pi'\in\hat{G}_{ell}} k_{\pi' ,\pi} N_{\Ga ,\omega}(\pi') \=
\dim \omega (-1)^{\frac{\dim X}{2}}\chi(X_\Ga) d_{\pi},
$$
where $d_{\pi}$ is the formal degree of $\pi$.

The infinitesimal character $\la$ of $\pi$ can be viewed as an
element of the coset space $(\t^*)^{reg}/W_K$. So let $J$ denote
the finite set of connected components of $(\t^*)^{reg}/W_K$, then
we get a decomposition $\hat{G}_{ds} = \coprod_{j\in J}
\hat{G}_{ds,j}$.

In the proof of the theorem we used the Hecht-Schmid character
formula to deal with the global characters. On the other hand it
is known that global characters are given on the regular set by
sums of toric characters over the Weyl denominator. So on $H=AB$
the character $\Theta_\pi$ for $\pi\in\hat{G}$ is of the form $\CN
/D$, where $D$ is the Weyl denominator and the numerator $\CN$ is
of the form
$$
\CN (h) \= \sum_{w\in W(\t,\g)} c_w h^{w\la},
$$
where $\la\in\h^*$ is the infinitesimal character of $\pi$.
Accordingly, the expression $m(\pi ,\tau ,\mu)$ expands as a sum
$$
m(\pi ,\tau ,\mu) \= \sum_{w\in W(\t ,\g)} m_w(\pi ,\tau ,\mu).
$$

\begin{lemma}
Let $\pi ,\pi'\in \hat{G}_{ds,j}$ with infinitesimal characters
$\la ,\la'$ which we now also view as elements of $(\h^*)^+$, then
$$
m_w(\pi ,\tau ,\mu) \= m_w(\pi',\tau ,\mu +w(\la'-\la)|_\a).
$$
\end{lemma}

\prf In light of the preceding it suffices to show the following:
Let $\tau_\la, \tau_{\la'}$ denote the numerators of the global
characters of $\pi$ and $\pi'$ on $\h^+$. Write
$$
\tau_\la (h) \= \sum_{w\in W(\t ,\g)}c_w h^{w\la}
$$
for some constants $c_w$. Then we have
$$
\tau_{\la'} (h) \= \sum_{w\in W(\t ,\g)}c_w h^{w\la'}.
$$
To see this, choose a $\la''$ dominating both $\la$ and $\la'$,
then apply the Zuckerman functors $\ph_{\la''}^\la$ and
$\ph_{\la''}^{\la'}$. Proposition 10.44 of \cite{Knapp} gives the
claim. \qed

Write $\pi_\la$ for the discrete series representation with
infinitesimal character $\la$. Let $d(\la):= d_{\pi_\la}$ be the
formal degree then $d(\la)$ is a polynomial in $\la$, more
precisely from \cite{AtSch} we take
$$
d(\la) \= \prod_{\alpha\in\Phi^+(\t ,\g)} \frac{(\alpha ,\la
+\rho)}{(\alpha ,\rho)},
$$
where the ordering $\Phi^+$ is chosen to make $\la$ positive.

Putting things together we see that for $\Re(s)$ small enough the
order of $Z_{P,\tau,\omega}(s)$ at $s=\mu(H_1)$ is
\begin{eqnarray*}
\CO (\mu) &\=& \dim \omega (-1)^{\frac{\dim X}{2}} \chi (X_\Ga)\\
    & {}&\sum_{j\in J} \sum_{w\in W(\t ,\g)} \sum_{\pi \in \hat{G}_{ds,j}}
    d(\la_\pi) m_w(\pi_j,\tau ,\mu +w(\la_j-\la_\pi)|_\a),
\end{eqnarray*}
where $\pi_j\in\hat{G}_{ds,j}$ is a fixed element. The function
$\mu \mapsto m_w(\pi_j ,\tau ,\mu)$ takes nonzero values only for
finitely many $\mu$. Since further $\la\mapsto d(\la)$ is a
polynomial it follows that the regularized product
$$
D_{P,\tau ,\omega}(s) := \widehat{\prod_{\mu ,\CO(\mu)\neq 0}}
(s-\mu(H_1))^{\CO(\mu)}
$$
exists. We now have proven the following theorem.

\begin{theorem}
With
$$
\hat{Z}_{P,\tau ,\omega}(s) := Z_{P,\tau,\omega}(s) D_{H,\tau
,\omega}(s)^{-1}
$$
we have
$$
\hat{Z}_{P,\tau ,\omega}(2| \rho_0| -s) \= e^{Q(s)}
\hat{Z}_{P,\tau ,\omega}(s),
$$ where $Q$ is a polynomial.
\qed
\end{theorem}

\begin{proposition} \label{extformel}
Let $\tau$ be a finite dimensional representation of $M$ then the
order of $Z_{P,\tau , \omega}(s)$ at $s=\la(H_1)$ is
$$
(-1)^{\dim(N)} \sum_{\pi \in \hat{G}}N_{\Ga,\omega}(^\theta\pi)
\sum_{q=0}^{\dim(\m\oplus{\n} /\k_M)}(-1)^q \dim(H^q(\m\oplus\n,
K_M ,\pi_K \otimes V_{\breve{\tau}})^\la).
$$
This can also be expressed as
$$
(-1)^{\dim(N)} \sum_{\pi \in \hat{G}}N_{\Ga,\omega}(^\theta\pi)
\sum_{q=0}^{\dim(\m\oplus\n /\k_M)}(-1)^q \dim({\rm Ext}_{(\m
\oplus {\n} ,K_M)}^q (V_{{\tau}} ,V_\pi)^\la).
$$
\end{proposition}

\prf Extend $V_{{\tau}}$ to a $\m \oplus \n$-module by letting
$\n$ act trivially. We then get
$$
H^p(\n,\pi_K) \otimes V_{{\tau}} \cong H^p(\n,\pi_K \otimes
V_{\breve{\tau}}).
$$

The $(\m ,K_M)$-cohomology of the module $H^p(\n,\pi_K \otimes
V_{\tau})$ is the cohomology of the complex $(C^*)$ with
\begin{eqnarray*}
C^q &\=& {\rm Hom}_{K_M}(\wedge^q\p_M ,H^p(\n,\pi_K)\otimes
V_{\breve{\tau}})
\\
        &\=& (\wedge^q\p_M \otimes H^p(\n,\pi_K)\otimes V_{\breve{\tau}})^{K_M},
\end{eqnarray*}
since $\wedge^p\p_M$ is a self-dual $K_M$-module. Therefore we
have an isomorphism of virtual $A$-modules:
$$
\sum_q (-1)^q (H^p(\n,\pi_K)\otimes \wedge^q\p_M \otimes
V_{\breve{\tau}} )^{K_M} \cong \sum_q (-1)^q H^q (\m
,K_M,H^p(\n,\pi_K \otimes V_{\breve{\tau}})).
$$

Now one considers the Hochschild-Serre spectral sequence in the
relative case for the exact sequence of Lie algebras
$$
0 \ra \n\ra \m \oplus \n\ra \m \ra 0
$$
and the $(\m\oplus \n,K_M)$-module $\pi \otimes V_{\breve{\tau}}$.
We have
$$
E_2^{p,q} \= H^q(\m ,K_M ,H^p(\n,\pi_K\otimes V_{\breve{\tau}}))
$$
and
$$
E_\infty^{p,q} \= {\rm Gr}^q(H^{p+q}(\m \oplus \n,K_M ,\pi_K
\otimes V_{\breve{\tau}})).
$$
Now the module in question is just
$$
\chi(E_2) \= \sum_{p,q} (-1)^{p+q} E_2^{p,q}.
$$
Since the differentials in the spectral sequence are
$A$-homomorphisms this equals $\chi(E_\infty)$. So we get an
$A$-module isomorphism of virtual $A$-modules
$$
\sum_{p,q} (-1)^{p+q} (H^p(\n,\pi_K)\otimes \wedge^q\p_M \otimes
V_{\breve{\tau}})^{K_M} \cong \sum_j (-1)^j H^j (\m\oplus \n,K_M,
\pi_K \otimes V_{\breve{\tau}}).
$$
The second statement is clear by \cite{BorWall} p.16. \qed

\section{The Ruelle zeta function}
The generalized Ruelle zeta function can be described in terms of
the Selberg zeta function as follows.

\begin{theorem}
Let $\Ga$ be  neat and choose a parabolic $P$ of splitrank
one. For $\Re(s)>>0$ define the zeta function
$$
Z_{P,\omega}^R(s) \= \prod_{[\ga]\in {\cal E}_H^p(\Ga)}
\det\left(\begin{array}{c}1-e^{-sl_\ga}\omega(\ga)\end{array}\right)^{\chi_{_1}(X_\ga)},
$$
then $Z_{P,\omega}^R(s)$ extends to a meromorphic function on
$\C$. More precisely, let $\n =\n_\alpha \oplus \n_{2\alpha}$ be
the root space decomposition of $\n$ with respect to the roots of
$(\a ,\g)$ then
$$
Z_{P,\omega}^R(s) = \prod_{q=0}^{\dim \n_\alpha} \prod_{p=0}^{\dim
\n_{2\alpha}}
Z_{P,(\wedge^q\n_\alpha)\otimes(\wedge^p\n_{2\alpha}),\omega}(s+(q+2p)|\alpha|)^{(-1)^{p+q}}.
$$
\end{theorem}

In the case when ${\rm rank}_\R G =1$ this zeta function coincides
with the Ruelle zeta function of the geodesic flow of $X_\Ga$.

\prf For any finite dimensional virtual representation $\tau$ of
$M$ we compute
$$
\log Z_{P,\tau,\omega}(s) = \sum_{[\ga]\in\CE_P^p(\Ga)}
\chi_{_1}(X_\ga)
        \sum_{N\ge 0} \tr(\log(1-e^{-sl_\ga}\ga)|\omega\otimes\tau\otimes S^N(\n))
$$
\begin{eqnarray*}
&=& \sum_{[\ga]\in\CE_P^p(\Ga)} \chi_{_1}(X_\ga)
        \sum_{N\ge 0} \sum_{n\ge 1} \rez{n} e^{-sl_\ga n}\tr(\ga |\omega\otimes\tau\otimes S^N(\n))\\
&=& \sum_{[\ga]\in\CE_P(\Ga)} \chi_{_1}(X_\ga)
\frac{e^{-sl_\ga}}{\mu(\ga)}
    \tr\omega(\ga)
    \frac{\tr\tau(b_\ga)}{\det(1-(a_\ga b_\ga)^{-1}|\n)}\\
&=& \sum_{[\ga]\in\CE_P(\Ga)} \chi_{_1}(X_\ga)
\frac{e^{-sl_\ga}}{\mu(\ga)}
    \tr\omega(\ga)
    \frac{\tr\tau(b_\ga)}{\tr((a_\ga b_\ga)^{-1}|\wedge^*\n)}.
\end{eqnarray*}
Since $\n$ is an $M$-module defined over the reals we conclude
that the trace $ \tr((a_\ga m_\ga)^{-1}|\wedge^*\n) $ is a real
number. Therefore it equals its complex conjugate which is
$\tr(a_\ga^{-1} m_\ga|\wedge^*\n)$. Now split into the
contributions from $\n_{\alpha}$ and $\n_{2\alpha}$. The claim now
becomes clear. \qed

\chapter{The Patterson conjecture}
In this chapter $G$ continues to be a connected semisimple
Lie group with finite center.

Fix a parabolic $P=MAN$ of splitrank one, ie, $\dim A=1$. Let $\nu \in \a^*$ and let
$({{\sigma}},V_\sigma)$ be a finite dimensional complex
representation of $M$. Consider the principal series
representation $\pi_{\sigma ,\nu}$ on the Hilbert space
$\pi_{\sigma ,\nu}$ of all functions $f$ from $G$ to
$V_\sigma$ such that $f(xman)=a^{-(\nu +\rho)}
\sigma(m)^{-1} f(x)$ and such that the restriction of $f$
to $K$ is an $L^2$-function. Let $\pi_{\sigma,\nu}^\infty$ denote the Fr\'echet space of smooth vectors
and
$\pi_{\sigma,\nu}^{-\infty}$ for its continuous dual. For
$\nu
\in
\a^*$ let
$\bar{\nu}$ denote its complex conjugate with respect to
the real form $\a_0^*$.

We will now formulate the Patterson conjecture \cite{BuOl}.

\begin{theorem} Suppose that $|W(G,A)|=2$. 
Then the cohomology  groups $H^p(\Ga ,\pi_{\sigma ,\nu}^{-\infty} \otimes
V_\omega)$ are finite dimensional for all $p\geq 0$ and the vanishing order of the Selberg zeta function
$Z_{P,\tau,\omega}(s+|\rho_0|)$ at $s=\nu(H_1)$ equals
$$
\ord_{s=\nu(H_1)}Z_{P,\tau,\omega}(s+|\rho_0|) \= \chi_{_1}(\Ga ,\pi_{\breve\sigma ,{-\nu}}^{-\infty}
\otimes V_{\omega})
$$
if $\nu \neq 0$ and
$$
\ord_{s=0}Z_{P,\tau,\omega}(s+|\rho_0|)\= \chi_{_1}(\Ga ,\hat{H}_{\breve\sigma,0}^{-\infty}\otimes
V_{\omega}),
$$
where $\hat{H}_{\sigma ,0}^{-\infty}$ is a certain nontrivial
extension of $\pi_{\sigma ,0}^{-\infty}$ with itself.

Further $\chi (\Ga ,\pi_{\sigma ,\nu}^{-\infty}\otimes
V_{\breve{\omega}})$ vanishes.
\end{theorem}

\prf Let $\a_0$ be the real
Lie algebra of $A$ and $\a_0^-$ the negative Weyl chamber with
$\exp(\a_0^-)=A^-$. Let $H_1\in \a_0^-$ be the unique element of
norm $1$. Recall that the vanishing order equals
$$
\sum_{\pi\in\hat{G}} N_{\Ga ,\omega}(\pi) \sum_{q=0}^{\dim
N}\sum_{p=0}^{\dim\p_M} (-1)^{p+q}
\dim(H_q({\n},\pi_K)^{\nu-\rho_0}\otimes \wedge^p\p_M\otimes
V_{\breve{\sigma}})^{K_M}.
$$
Set $\la=\nu-\rho_0$.

\begin{lemma}
If $\la\ne -\rho_0$ then $\a$ acts semisimply on
$H_q({\n},\pi_K)^\la$. The space $H_\rho({\n},\pi_K)^{-\rho_0}$
is annihilated by $(H+\rho_0(H))^2$ for any $H\in\a$.
\end{lemma}

Since $W(G,A)$ is nontrivial, the proof in \cite{BuOl} Prop. 4.1. carries over to our situation.
\qed

To prove the theorem, assume first $\la\ne -\rho_0$. Let
$\pi\in\hat{G}$, then
$$
\sum_{q=0}^{\dim N}\sum_{p=0}^{\dim\p_M} (-1)^{p+q} \dim(
H_q({\n},\pi_K)^\la \otimes \wedge^p\p_M\otimes
V_{\breve{\sigma}} )^{K_M}
$$ $$
=\ \sum_{q=0}^{\dim N} \sum_{p=0}^{\dim\p_M} (-1)^{p+q} \dim(
H^0(\a ,H^p(\m,K_M,H_q({\n},\pi_K) \otimes V_{{\breve{\sigma}}
,-\la}))),
$$
where $V_{{\breve{\sigma}} ,-\la}$ is the representation
space of the representation ${\breve{\sigma}}\otimes
(-\la)$. We want to show that this equals
$$
\sum_{q=0}^\infty \sum_{p=0}^\infty (-1)^{p+q+r}p
\dim H^p(\a\oplus\m,K_M, H_q({\n},\pi_K)\otimes
V_{{\breve{\sigma}} ,-\la}).
$$ 

Since the Hochschild-Serre spectral sequence degenerates for a one dimensional Lie
algebra we get that
$$
\dim H^p(\a\oplus\m,K_M,V)
$$
equals
$$
\dim H^0(\a,H^{p-1}(\m,K_M,V))+\dim H^0(\a,H^{p}(\m,K_M,V)).
$$
This implies
$$
\sum_{p,q\ge 0} (-1)^{p+q+r}p \sum_{b=p-1}^p
\dim H^0(\a,H^b(\m,K_M, H_q({\n},\pi_K)\otimes
V_{{\breve{\sigma}} ,-\la}))
$$ $$
\= \sum_{b,q\ge 0}(-1)^q \sum_{p=b}^{b+1} p\, \dim H^0(\a,H^b(\m,K_M,
H_q({\n},\pi_K)\otimes V_{{\breve{\sigma}} ,-\la})).
$$
It follows that the vanishing order equals
$$
\sum_{q=0}^\infty \sum_{p=0}^\infty (-1)^{p+q+1} p
\dim H^p(\a\oplus\m,K_M, H_q({\n},\pi_K)\otimes
V_{{\breve{\sigma}} ,-\la}).
$$ $$
=\ \sum_{q=0}^\infty \sum_{p=0}^\infty (-1)^{p+q+1} \dim Ext_{\a\oplus\m,K_M}^p(
H_q({\n},\breve{\pi_K}), V_{{\breve{\sigma}} ,-\la}).
$$

For a $(\g,K)$-module $V$ and a $(\a\oplus\m,M)$-module $U$ we
have \cite{HeSch}:
$$
Hom_{\g,K}(V,Ind_P^G(U)) \ \cong\ Hom_{\a\oplus\m,K_M}(H_0(\n
,V),U\otimes\C_{\rho_0}),
$$
where $\C_{\rho_0}$ is the one dimensional $A$-module
given by $\rho_0$. 
Thus
$$
\sum_{\pi\in\hat{G}} N_{\Ga,\omega}(\pi) \sum_{p=0}^\infty
\sum_{q=0}^\infty (-1)^{p+q+1} \dim
Ext_{\a\oplus\m,K_M}^p
(H_q(\n,\breve{\pi_K}),\pi_{{\breve{\sigma}},-\la +\rho_P,K}^\infty)
$$
\begin{eqnarray*}
    &=& \sum_{p=0}^\infty (-1)^{p+1} \dim Ext_{\g,K}^p (C^\infty(\Ga\bs G,\breve{\omega})_K,
    \pi_{{\breve{\sigma}}, -\la,K}^\infty)
\end{eqnarray*}
Diualizing shows that this equals
$$
\sum_{p\ge0} (-1)^{p+1}\dim Ext_{\g,K}^p(\pi_{\sigma,\la,K}^\infty,C^\infty(\Ga\bs
G,\omega)_K).
$$
Next by \cite{BorWall}, Chap. I,
$$
Ext_{\g,K}^p(\pi_{\sigma,\la,K}^\infty,C^\infty(\Ga\bs
G,\omega)_K)\ \cong\ H^p(\g,K,\Hom_\C(\pi_{\sigma,\la,K}^\infty,C^\infty(\Ga\bs G,\omega)_K).
$$
For any two smooth $G$-representations $V,W$ the restriction map gives an isomorphism
$$
\Hom_{ct}(V,W)_K\ \cong\ \Hom_\C(V_K,W_K)_K,
$$
where $\Hom_{ct}$ means continuous homomorphisms. Therefore, using the classical identification of
$(\g,K)$ with differentiable and continuous cohomology as in
\cite{BorWall} we get
\begin{eqnarray*}
Ext_{\g,K}^q(\pi_{\sigma,\la,K}^\infty,C^\infty(\Ga\bs G,\omega)_K) &\cong &
H^p(\g,K,\Hom_{ct}(\pi_{\sigma,\la}^\infty,C^\infty(\Ga\bs G,\omega))_K)\\
&\cong &
H^p(\g,K,\Hom_{ct}(\pi_{\sigma,\la}^\infty,C^\infty(\Ga\bs G,\omega)))\\
&\cong &
H_d^p(G,\Hom_{ct}(\pi_{\sigma,\la}^\infty,C^\infty(\Ga\bs G,\omega)))\\
&\cong &
H_{ct}^p(G,\Hom_{ct}(\pi_{\sigma,\la}^\infty,C^\infty(\Ga\bs G,\omega)))\\
&\cong &
Ext_G^p(\pi_{\sigma,\la}^\infty, C^\infty(\Ga\bs G,\omega))\\
&\cong &
Ext_\Ga^p(\pi_{\sigma,\la}^\infty,\omega)\\
&\cong&
H^p(\Ga,\pi_{\breve\sigma,-\la}^{-\infty}\otimes\omega).
\end{eqnarray*}
This gives the claim for $s\ne 0$. In the case $s=0$ we have to replace $H^0(\a,\cdot)$ by
$H^0(\a^2,\cdot)$, where $\a^2$ means the subalgebra of $U(\a)$ generated by $H^2$, $H\in\a$. Then
the induced representation $Ind_P^G(U)$ is replaced by a suitable self-extension.
\qed

\chapter{Holomorphic torsion}
\section{Holomorphic torsion}

Let $E=(E_0 \rightarrow E_1 \rightarrow \dots \rightarrow E_n)$ be
an elliptic complex over a compact smooth manifold $M$. That is,
each $E_j$ is a vector bundle over $M$ and there is differential
operators $d_j:\Ga^\infty(E_j)\ra\Ga^\infty(E_{j+1})$ of order one
such that the sequence of principal symbols
$$
\sigma_E(x,\xi)\ :\ 0\ra E_{0,x}\ra E_{1,x}\ra\dots\ra E_{n,x}\ra
0
$$
is exact whenever $\xi\in T_xM^*$ is nonzero.

Assume each $E_k$ is equipped with a Hermitian metric. Then we can
form the Laplace operators $\lap_j= d_jd_j^*+d_j^*d_j$ as second
order differential operators. When considered as unbounded
operators on the spaces $L^2(M,E_k)$ these are known to be zeta
admissible. Now define the \emph{torsion} of $E$ as $$ \tau_1 (E)
\= \prod_{k=0}^n \det'(\lap_k)^{k(-1)^{k+1}}.
$$ Note that this definition differs by an exponent 2 from the
original one \cite{RS-RT}.

The corresponding notion in the combinatorial case, i.e. for a
finite CW-complex and the combinatorial Laplacians, was introduced
by Reidemeister in the 1930's who used it to distinguish homotopy
equivalent spaces which are not homeomorphic. The torsion as
defined above, also called {\it analytic torsion}\index{analytic
torsion}, was defined by Ray and Singer in \cite{RS-RT} where they
conjectured that the combinatorial and the analytical torsion
should coincide. This conjecture was later proven independently by
J. Cheeger and W. M\"uller.

For a compact smooth Riemannian manifold $M$ and $E\rightarrow M$
a flat Hermitian vector bundle the complex of $E$-valued forms on
$M$ satisfies the conditions above so that we can define the
torsion $\tau (E)$ via the de Rham complex. Now assume further,
$M$ is K\"{a}hlerian and $E$ holomorphic then we may also consider
the torsion $T(E)$ of the Dolbeault complex $\bar{\partial}\ :\
\Omega^{0,.}(M,E) \rightarrow \Omega^{0,.+1}(M,E)$. This then is
called the {\it holomorphic torsion}\index{holomorphic torsion}.
The holomorphic torsion is of significance in Arakelov theory
\cite{soule} where it is used as normalization factor for families
of Hermitian metrics.

We now define \emph{$L^2$-torsion}. For the following see also
\cite{L}. Let $M$ denote a compact oriented smooth manifold, $\Ga$
its fundamental group and $\tilde{M}$ its universal covering. Let
$E=E_0\rightarrow \dots \rightarrow E_n$ denote an elliptic
complex over $M$ and $\tilde{E}=\tilde{E}_0\rightarrow \dots
\rightarrow \tilde{E}_n$ its pullback to $\tilde{M}$. Assume all
$E_k$ are equipped with Hermitian metrics.

Let $\tilde{\triangle}_p$ and $\triangle_p$ denote the
corresponding Laplacians. The ordinary torsion was defined via the
trace of the complex powers $\triangle_p^s$. The $L^2$-torsion
will instead be defined by considering the complex powers of
$\tilde{\triangle}_p$ and applying a different trace functional.
Write ${\cal F}$ for a fundamental domain of the $\Ga$-action on
$\tilde{M}$ then as a $\Ga$-module we have $$ L^2(\tilde{E}_p)\
\cong\ l^2(\Ga) \otimes L^2(\tilde{E}_p\mid_{\cal F})\ \cong\
l^2(\Ga) \otimes L^2(E_p). $$ The von Neumann algebra $VN(\Ga)$
generated by the right action of $\Ga$ on $l^2(\Ga)$has a
canonical trace making it a type ${\rm II}_1$ von Neumann algebra
if $\Ga$ is infinite \cite{GHJ}. This trace and the canonical
trace on the space $B(L^2(E))$ of bounded linear operators on
$L^2(E)$ define a trace ${\rm tr}_\Ga$ on $VN(\Ga) \otimes
B(L^2(E))$ which makes it a type ${\rm II}_\infty$ von Neumann
algebra. The corresponding dimension function is denoted
$\dim_\Ga$.  Assume for example, a $\Ga$-invariant operator $T$ on
$L^2(E)$ is given as an integral operator with a smooth kernel
$k_T$, then a computation shows
$$
\tr_\Ga(T) \= \int_{\cal F} \tr(k_T(x,x))\ dx.
$$
It follows for the heat operator $e^{-t\tilde{\lap}_p}$ that $$
\tr_\Ga e^{-t\tilde{\lap}_p} \= \int_{\cal F}\tr <x\mid
e^{-t\tilde{\lap}_p}\mid x> dx. $$ From this we read off that
$\tr_\Ga e^{-t\tilde{\lap}_p}$ satisfies the same small time
asymptotic as $\tr e^{-t\lap_p}$.

Let $\tilde{\lap}_p' =
\tilde{\lap}_p|_{\ker(\tilde{\lap}_p)^\perp}$. Unfortunately very
little is known about large time asymptotic of
$\tr_\Ga(e^{-t\tilde{\lap}_p'})$ (see \cite{LL}). Let $$
NS(\lap_p) \= \sup \{ \alpha \in \R \mid \tr_\Ga e^{-t
\tilde{\lap}_p'} \= O(t^{-\alpha/2})\ {\rm as}\ t\rightarrow
\infty \} $$ denote the \emph{Novikov-Shubin invariant} of $\lap_p$
(\cite{GrSh}, \cite{LL}).

Then $NS(\lap_p)$ is always $\geq 0$; in this section we will \emph{assume} that the Novikov-Shubin invariant of $\lap_p$ is positive.
This is in general an unproven conjecture. In the cases of our
concern in later sections, however, the operators in question are
homogeneous and it can be proven then that their Novikov-Shubin
invariants are in fact positive. We will consider the integral $$
\zeta_{\lap_p}^1(s) \= \rez{\Ga (s)} \int_0^1 t^{s-1} \tr_\Ga
e^{-t\tilde{\lap}_p'}\ dt, $$ which converges for $\Re (s) >>0$
and extends to a meromorphic function on the entire plane which is
holomorphic at $s=0$, as is easily shown by using the small time
asymptotic (\cite{BGV},Thm 2.30).

Further the integral
$$
\zeta_{\lap_p(s)}^2(s) \= \rez{\Ga(s)}\int_1^\infty
t^{s-1}\tr_\Ga e^{-t\tilde{\lap}_p'}\ dt
$$
converges for $\Re (s)<\rez{2}NS(\lap_p)$, so in this region we
define the \emph{$L^2$-zeta function} of $\lap_p$ as
$$
\zeta_{\lap_p}^{(2)} (s) \= \zeta_{\lap_p}^1(s) +
\zeta_{\lap_p}^2(s).
$$
Assuming the Novikov-Shubin invariant of $\lap_p$ to be positive
we define the \emph{$L^2$-determinant} of $\lap_p$ as
$$
{\det}^{(2)}(\lap_p) \= \exp \left( \left. -\frac{d}{d{s}}
\right|_{s=0} \zeta_{\lap_p}^{(2)} (s)\right) .
$$

Now let the $L^2$-torsion be defined by $$ T^{(2)}(E) \=
\prod_{p=0}^n {\det}^{(2)}(\lap_p)^{p(-1)^{p+1}}. $$ Again let $M$
be a K\"ahler manifold and $E\ra M$ a flat Hermitian holomorphic
vector bundle then we will write $T^{(2)}_{hol}(E)$ for the
$L^2$-torsion of the Dolbeault complex $\Omega^{0,*}(M,e)$.

\section{Computation of Casimir eigenvalues}
In this section we give formulas for Casimir eigenvalues of
irreducible representations of reductive groups which hold under
conditions on the $K$-types.

Let $G$ be a connected semisimple Lie group with finite center.
Fix a maximal compact subgroup $K$ with Cartan involution
$\theta$. Let $P=MAN$ be the Langlands decomposition of a
parabolic subgroup $P$ of $G$. Modulo conjugation we can assume
that $AM$ is stable under $\theta$ and then $K_M=K\cap M$ is a
maximal compact subgroup of $M$. Let $\m =\k_M\oplus\p_M$ be the
corresponding Cartan decomposition of the complexified Lie algebra
of $M$. Let $C_M$ denote the Casimir operator of $M$ induced by
the Killing form on $G$.

\begin{lemma} \label{Caseig1}
Let $(\sigma ,V_\sigma)$ be an irreducible finite dimensional
representation of $M$. Let $(\xi ,V_\xi)$ be an irreducible
unitary representation of $M$ and assume
$$
\sum_{p=0}^{\dim (\p)} (-1)^p \dim (V_\xi \otimes \wedge^p\p_M
\otimes V_{\sigma})^{K_M} \neq 0,
$$
then the Casimir eigenvalues satisfy
 $$
 \xi(C_M) = \sigma(C_M).
 $$
\end{lemma}

\prf Recall that the Killing form of $G$ defines a
$K_M$-isomorphism between $\p_M$ and its dual $\p_M^*$, hence in
the assumption of the lemma we may replace $\p_M$ by $\p_M^*$. Let
$\xi_K$ denote the $(\m ,K_M)$-module of $K_M$-finite vectors in
$V_\xi$ and let $C^q(\xi_K\otimes V_\sigma) =
\Hom_{\k_M}(\wedge^q\p_M ,\xi_K\otimes V_\sigma) =
(\wedge^q\p_M^*\otimes\xi_K\otimes V_\sigma)^{\k_M}$ the standard
complex for the relative Lie algebra cohomology $H^q(\m ,\k_M
,\xi_K\otimes V_\sigma)$. Further
$(\wedge^q\p_M^*\otimes\xi_K\otimes V_\sigma)^{K_M}$ forms the
standard complex for the relative $(\m ,K_M)$-cohomology $H^q(\m
,K_M ,\xi_K\otimes V_\sigma)$. In \cite{BorWall}, p.28 it is shown
that
$$
H^q(\m ,K_M ,\xi_K\otimes V_\sigma)=H^q(\m ,\k_M
,\xi_K\otimes V_\sigma)^{K_M/K_M^0}.
$$
Our assumption
implies $\sum_q(-1)^q \dim H^q(\m ,K_M ,\xi_K\otimes
V_\sigma) \ne 0$, therefore there is a $q$ with $0\ne
H^q(\m ,K_M ,\xi_K\otimes V_\sigma)=H^q(\m ,\k_M
,\xi_K\otimes V_\sigma)^{K_M/K_M^0}$, hence $H^q(\m ,\k_M
,\xi_K\otimes V_\sigma)\ne 0$. Now Proposition 3.1 on page
52 of \cite{BorWall} says that $\pi(C_M)\ne\sigma(C_M)$
implies that $H^q(\m,\k_M,\xi_K\otimes V_\sigma)=0$ for
all $q$. The claim follows. \qed

Let $X=G/K$ be the symmetric space to $G$ and assume that $X$ is
Hermitian\index{Hermitian symmetric space}, i.e., $X$ has a
complex structure which is stable under $G$. Let $\theta$ denote
the Cartan involution fixing $K$ pointwise. Since $X$ is
hermitian it follows that $G$ admits a compact Cartan subgroup
$T\subset K$. We denote the real Lie algebras of $G$, $K$ and $T$
by $\g_0, \k_0$ and $\t_0$ and their complexifications by $\g
,\k$ and $\t$. We will fix a scalar multiple $B$ of the Killing
form as explained in the chapter on notation. As well, we will
write $B$ for its diagonal, so $B(X) = B(X,X)$. Denote by $\p_0$
the orthocomplement of $\k_0$ in $\g_0$ with respect to $B$; then
via the differential of {\rm exp} the space $\p_0$ is isomorphic
to the real tangent space of $X=G/K$ at the point $eK$. Let $\Phi
(\t ,\g)$ denote the system of roots of $(\t ,\g)$, let $\Phi_c
(\t ,\g) = \Phi (\t ,\k)$ denote the subset of compact roots and
$\Phi_{nc} = \Phi -\Phi_c$ the set of noncompact roots. To any
root $\alpha$ let $\g_\alpha$ denote the corresponding root
space. Fix an ordering $\Phi^+$ on $ \Phi = \Phi (\t ,\g)$ and
let $\p_\pm = \bigoplus_{\alpha \in \Phi_{nc}^+}\g_{\pm \alpha}$.
Then the complexification $\p$ of $\p_0$ splits as $\p = \p_+
\oplus \p_-$ and the ordering can be chosen such that this
decomposition corresponds via {\rm exp} to the decomposition of
the complexified tangent space of $X$ into holomorphic and
antiholomorphic part. By Lemma \ref{double-cover} we can,
replacing $G$ by a double cover if necessary, assume that the
adjoint homomorphism $K\ra SO(\p)$ factors over $Spin(\p)$.

\begin{lemma} \label{Caseig2}
Let $(\tau ,V_\tau)$ denote an irreducible representation of $K$.
Assume $X$ Hermitian and let $(\pi ,W_\pi)$ be an irreducible
unitary representation of $G$ and assume that
$$
\sum_{p=0}^{\dim \p_-} (-1)^p \dim(W_\pi \otimes
\wedge^p\p_-\otimes V_{{\tau}})^K \neq 0
$$
then we have
$$
\pi(C) \= {\tau}\otimes\epsilon(C_K) -B(\rho)+B(\rho_K),
$$
where $\epsilon$ is the one dimensional representation of $K$
satisfying $\epsilon\otimes\epsilon\cong\wedge^{top}\p_+$.
\end{lemma}

\prf This follows from Lemma \ref{epsilon}. \qed

\section{The local trace of the heat kernel}
On a Hermitian globally symmetric space the heat operator is given
by convolution with a function on the group of isometries. In this
section we determine the trace of this function on any irreducible
unitary representation.

Let $H$ denote a $\theta$-stable Cartan subgroup of $G$ so $H=A B$
where $A$ is the connected split component and $B$ compact. The
dimension of $A$ is called the split rank of H. Let $\a$ denote
the complex Lie algebra of $A$. Then $\a$ is an abelian subspace
of $\p=\p_+ \oplus \p_-$. Let $X\mapsto X^c$ denote the complex
conjugation on $\g$ according to the real form $\g_0$. The next
lemma shows that $\a$ lies skew to the decomposition $\p =\p_+
\oplus \p_-$.

\begin{lemma} Let $Pr_\pm$ denote the projections from $\p$ to $\p_\pm$.
Then we have $\dim Pr_+(\a)$ = {\rm dim} $Pr_-(\a)$ = {\rm dim}
$\a$, or, what amounts to the same: $\a \cap \p_\pm = 0$.
\end{lemma}

Proof: $\a$ consists of semisimple elements whereas $\p_\pm$
consists of nilpotent elements only.
 \qed

Now let $\a'$ denote the orthocomplement of $\a$ in $Pr_+(\a)
\oplus Pr_-(\a)$. For later use we write ${\cal V} = \a \oplus \a'
= {\cal V}_+ \oplus {\cal V}_-$, where ${\cal V}_\pm = {\cal
V}\cap \p_\pm = Pr_\pm (\a)$.

Let $(\tau,V_\tau)$ be a finite dimensional representation of the
compact group $K$ and let $E_\tau = G\times_K V_\tau = (G\times
V)/K$ be the corresponding $G$-homogeneous vector bundle over $X$.
The sections of the bundle $E_\tau$ are as a $G$-module given by
the space of $K$-invariants
 $$
 \Ga^\infty(E_\tau) \= \left( C^\infty(G)\otimes V_\tau\right)^K,
 $$
where $K$ acts on $C^\infty(G)\otimes V_\tau$ by $k.(f(x)\otimes
v) = f(xk^{-1})\otimes \tau(k)v$. and $G$ acts on $\left(
C^\infty(G)\otimes V_\tau\right)^K$ by left translations on the
first factor.

Consider the convolution algebra $C_c^\infty(G)$ of compactly
supported smooth functions on $G$. Let $K\times K$ act on it by
right and left translations and on $\End_\C(V_\tau)$ by
$(k_1,k_2)T=\tau(k_1)T\tau(k_2^{-1})$. Then the space of
invariants
 $$
 \left(C_c^\infty(G)\otimes \End_\C(V_\tau)\right)^{K\times K}
 $$
is seen to be an algebra again and it acts on
$\Ga^\infty(E_\tau)=\left( C^\infty(G)\otimes V_\tau\right)^K$ by
convolution on the first factor and in the obvious way on the
second. This describes all $G$-invariant smoothing operators on
$E_\tau$ which extend to all sections. $G$-invariant smoothing
operators which not extend to all sections will be given by
Schwartz kernels in
 $$
 \left(C^\infty(G)\otimes \End_\C(V_\tau)\right)^{K\times K}.
 $$

Let $n=2m$ denote the real dimension of $X$ and for $0\leq p,q
\leq m$ let $\Omega^{p,q}(X)$ denote the space of smooth
$(p,q)$-forms on $X$. The above calculus holds for the space of
sections $\Omega^{p,q}(X)$. D.Barbasch and H. Moscovici have shown
in \cite{BM} that the heat operator $e^{-t\triangle_{p,q}}$ has a
smooth kernel $h_t^{p,q}$ of rapid decay in
 $$
 ( C^\infty (G)\otimes {\rm End}
 (\wedge^p\p_+ \otimes \wedge^q\p_-))^{K\times K}.
 $$
Now fix p and set for $t>0$ $$ f_t^p \= \sum_{q=0}^m q(-1)^{q+1}\
{\rm tr}\ h_t^{p,q}, $$ where tr means the trace in ${\rm End}
(\wedge^p\p_+ \otimes \wedge^q\p_-)$.

We want to compute the trace of $f_t^p$ on the principal series
representations. To this end let $H=A B$ as above and let $P$
denote a parabolic subgroup of $G$ with Langlands decomposition
$P=MA N$. Let $(\xi ,W_\xi)$ denote an irreducible unitary
representation of $M$, $e^\nu$ a quasicharacter of $A$ and set
$\pi_{\xi ,\nu} = Ind_P^G (\xi \otimes e^{\nu + \rho_{_P}}\otimes
1)$, where $\rho_{_P}$ is the half of the sum of the $P$-positive
roots.

Let $C$ denote the Casimir operator of $G$ attached to the form
$B$.

\begin{proposition} \label{trace1}
The trace of $f_t^p$ under $\pi_{\xi ,\nu}$ vanishes if $\dim \a
> 1 $. If ${\rm  dim\ } \a =1$ it equals
$$
e^{t\pi_{\xi ,\nu}(C)} \sum_{q=0}^{\dim(\p_-)-1} (-1)^q\
\dim\left( W_\xi \otimes \wedge^p\p_+ \otimes \wedge^q(\a^\perp
\cap \p_-)\right)^{K\cap M},
$$
where $\a^\perp$ is the orthocomplement of $\a$ in $\p$.
\end{proposition}

Proof: As before let $\a'$ be the span of all $X-X^c$, where
$X+X^c\in \a,\ X\in \p_+$ then ${\cal V} = \a \oplus \a' = \CV_+
\oplus \CV_-$ where $\CV_\pm = \CV \cap \p_\pm$. The group
$K_M=K\cap M$ acts trivially on $\a$ so for $x\in K_M$ we have
$X+X^c = {\rm Ad}(x)(X+X^c) = {\rm Ad}(x)X + {\rm Ad}(x)X^c$.
Since $K$ respects the decomposition $\p=\p_-\oplus \p_+$, we
conclude that $K_M$ acts trivially on $\CV$, hence on $\CV_-$.
Let $r=\dim\a=\dim\CV_-$. As a $K_M$-module we have $$
\begin{array}{cll}
\wedge^p\p_-    &       \=       &       \displaystyle
                                        \sum_{a+b=q} \wedge^b \CV_- \otimes
                                        \wedge^a \CV_-^\perp\\
                &       \=       &       \displaystyle
                                        \sum_{a+b=q} \binom{r}{b}
                                        \wedge^a \CV_-^\perp,
\end{array}
$$
where $\CV_-^\perp = \a^\perp \cap \p_-$. By definition we get
$$
\begin{array}{cll}
\tr\pi_{\xi ,\nu}(f_t^p)    & \= &   \displaystyle
                                \tr \pi_{\xi ,\nu}\left(\sum_{q=0}^m q(-1)^{q+1}
                                        h_t^{p,q}\right)\\
                & \= &   \displaystyle
                        e^{t\pi_{\xi ,\nu}(C)} \sum_{q=0}^m q(-1)^{q+1}
                                \dim(V_{\pi_{\xi ,\nu}} \otimes
                                \wedge^p\p_+ \otimes \wedge^q\p_-)^K
\end{array} $$
By Frobenius reciprocity this equals
\begin{eqnarray*}
                & {} &   \displaystyle
                        e^{t\pi_{\xi ,\nu}(C)} \sum_{q=0}^m q(-1)^{q+1}
                                \dim(W_{\pi_{\xi}} \otimes
                                \wedge^p\p_+ \otimes \wedge^q\p_-)^{K\cap M}\\
                & \= &   \displaystyle
                        e^{t\pi_{\xi ,\nu}(C)} \sum_{q=0}^m \sum_{a=0}^q q(-1)^{q+1}
                                \binom{r}{q-a} \dim(W_{\pi_{\xi}} \otimes
                                \wedge^p\p_+ \otimes
                                \wedge^a\CV_-^\perp)^{K\cap M}\\
                & \= &   \displaystyle
                        e^{t\pi_{\xi ,\nu}(C)} \sum_{a=0}^m \sum_{q=a}^m q(-1)^{q+1}
                                \binom{r}{q-a} \dim(W_{\pi_{\xi}} \otimes
                                \wedge^p\p_+ \otimes
                                \wedge^a\CV_-^\perp)^{K\cap M}
\end{eqnarray*}
By taking into account $a\leq m-r$ we get
 $$
 \sum_{q=a}^m q(-1)^q \binom{r}{q-a} \=
 \left\{ \begin{array}{cl}(-1)^{a+1}&{\rm if}\ r=1,\\ 0 & {\rm if}\ r>1,\end{array}\right.
 $$
 and
the claim follows. \qed

We will now combine this result with the computation of Casimir
eigenvalues. This requires an analysis whether the homogeneous
space $M/K_M$ embedded into $X=G/K$ inherits the complex structure
or not.

Fix a $\theta$-stable Cartan subgroup $H=A B$ with $\dim(A) =1$
and a parabolic $P=MA N$. Fix a system of positive roots
$\Phi^+=\Phi^+(\g,\h)$ in $\Phi(\g,\h)$ such that for
$\alpha\in\Phi^+$ and $\alpha$ nonimaginary it follows
$\alpha^c\in\Phi^+$. Further assume that $\phi^+$ is compatible
with the choice of $P$, i.e., for any $\alpha\in\Phi^+$ the
restriction $\alpha|_{\a}$ is either zero or positive. Let $\rho$
denote the half sum of positive roots. For $\xi \in \hat{M}$ let
$\la_\xi \in \b^*$ denote the infinitesimal character of $\xi$.
Recall that we have
$$ \pi_{\xi ,\nu}(C) \= B(\nu)+B(\la_\xi)-B(\rho). $$

\begin{lemma}
There exists a unique real root $\alpha_r\in\Phi^+$.
\end{lemma}

\prf Recall that a root $\alpha\in\Phi(\g,\h)$ is real if and only
if it annihilates $\b=Lie_\C(B)$. Hence the real roots are
elements of $\a^*$ which is one dimensional, so, if there were
two positive real roots, one would be positive multiple of the
other which is absurd.

To the existence. We have that $\dim X=\dim\a +\dim\n$ is even,
hence $\dim\n$ is odd. Now $\n=\oplus_\alpha\g_\alpha$, where the
sum runs over all $\alpha\in\Phi^+$ with $\alpha|_\a\ne 0$. The
complex conjugation permutes the $\g_\alpha$ and for any nonreal
$\alpha$ we have $(\g_\alpha)^c\ne\g_\alpha$, hence the nonreal
roots pair up, thus $\n$ can only be odd dimensional if there is a
real root.\qed

Let $c=c(H)$ denote the number of positive restricted roots in
$\Phi(\g,\a)$. Let $\Phi^+(\g,\a)$ denote the subset of positive
restricted roots.

\begin{lemma}
There are three possibilities:
\begin{itemize}
\item
$c=1$ and $\Phi^+(\a,\g)=\{ \alpha_r\}$,
\item
$c=2$ and $\Phi^+(\a,\g)=\{ \rez{2}\alpha_r, \alpha_r\}$ and
\item
$c=3$ and $\Phi^+(\a,\g)=\{ \rez{2}\alpha_r, \alpha_r,
\frac{3}{2}\alpha_r\}$.
\end{itemize}
\end{lemma}

\prf For any root $\alpha\in\Phi(\h ,\g)$ we have that $2B(\alpha
,\alpha_r)/B(\alpha_r)$ can only take the values $0,1,2,3$, so
the only possible restricted roots in $\Phi^+(\a,\g)$ are
$\rez{2}\alpha_r, \alpha_r, \frac{3}{2}\alpha_r$.

If $\frac{3}{2}\alpha_r\in \Phi^+(\a,\g)$ then there is a root
$\beta$ in $\Phi(\h,\g)$ with $\beta |_\a
=\frac{3}{2}\alpha_r|_\a$. Then $B(\alpha_r ,\beta)>0$ hence
$\eta =\beta-\alpha_r$ is a root. Since $\eta |_\a
=\rez{2}\alpha_r$ we get that then $\rez{2}\alpha_r\in\Phi(\a
,\g)$. From this the claim follows.
\qed

Recall that a {\it central isogeny}\index{central isogeny} $\ph
:L_1\ra L_2$ of Lie groups is a surjective homomorphism with
finite kernel which lies in the center of $L_1$.

\begin{lemma}
There is a central isogeny $M_1\times M_2\ra M$ such that the
inverse image of $K_M$ is of the form $K_{M_1}\times K_{M_2}$
where $K_{M_j}$ is a maximal compact subgroup of $M_j$ for
$j=1,2$ and such that $\p_M = \p_{M_1}\oplus\p_{M_2}$ as
$K_{M_1}\times K_{M_2} $-module and
\begin{itemize}
\item
the map $X\mapsto [X,X_{\alpha_r}]$ induces a $K_M$-isomorphism
$$
\p_{M_1}\ \cong\ [\p_M ,\g_{\alpha_r}],
$$
\item
with $\p_{M_2 ,\pm}:=\p_{M_2}\cap\p_\pm$ we have
$$
\p_{M_2} =\p_{M_2,+}\oplus\p_{M_2,-}.
$$
\end{itemize}
The latter point implies that the symmetric space $M_2/K_{M_2}$,
naturally embedded into $G/K$, inherits the complex structure.
This in particular implies that $M_2$ is orientation preserving.
Further we have $\p_{M_1}\cap\p_+ =\p_{M_1}\cap\p_- = 0$, which
in turn implies that the symmetric space $M_1/K_{M_1}$ does not
inherit the complex structure of $G/K$.
\end{lemma}

\prf At first we reduce the proof to the case that the center $Z$
of $G$ is trivial. So assume the proposition proved for $G/Z$ then
the covering $M_1\times M_2\ra M$ is gotten by pullback from that
of $M/Z$. Thus we may assume that $G$ has trivial center.

Let $H$ be a generator of $\a_0$. Write $H=Y+Y^c$ for some $Y\in
\p_+$. According to the root space decomposition $\g = \a \oplus
\k_M \oplus \p_M \oplus \n \oplus \theta(\n)$ we write
$Y=Y_a+Y_k+Y_p+Y_n+Y_{\theta(n)}$. Because of $\theta(Y)=-Y$ it
follows $Y_k=0$ and $Y_{\theta(n)}=-\theta(Y_n)$. For arbitrary
$k\in K_M$ we have $\Ad(k)Y=Y$ since $\Ad(k)H=H$ and the
projection $Pr_+$ is $K_M$-equivariant. Since the root space
decomposition is stable under $K_M$ it follows $\Ad(K)Y_*=Y_*$ for
$*=a,p,n$. The group $M$ has a compact Cartan, hence $K_M$ cannot
act trivially on any nontrivial element of $\p_M$, so $Y_p=0$ and
$Y_n\in \g_{\ar}$. Since $Y=Y_a+Y_n-\theta(Y_n)$ and $Y\notin\a$
it follows that $Y_n\ne 0$, so $Y_n$ generates $\g_{\alpha_r}$ and
so $K_M$ acts trivially on $\g_{\ar}$.

Let $\p_{M_2}\subset\p_M$ by definition be the kernel of the map
$X\mapsto [X,Y_n]$. Let $\p_{M_1}$ be its orthocomplement in
$\p_M$. The group $K_M$ stabilizes $Y_n$, so it follows that $K_M$
leaves $\p_{M_2}$ stable hence the orthogonal decomposition $\p_M
=\p_{M_1}\oplus\p_{M_2}$ is $K_M$-stable. Thus the symmetric space
$M/K_M$ decomposes into a product accordingly and so does the
image of $M$ in the group of isometries of $M/K_M$. It follows
that the Lie algebra $\m$ of $M$ splits as a direct sum of ideals
$\m=\m_0\oplus\m_1\oplus\m_2$, where $\m_0$ is the Lie algebra of
the kernel $M_0$ of the map $M\ra{\rm
Iso}(M/K_M)=\overline{M_1}\times\overline{M_2}$. Let $L=AM$ be the
Levi component and let $L_c$ denote the centralizer of $\a$ in the
group ${\rm Int}(\g)$. Then $L_c$ is a connected complex group
and, since $G$ has trivial center, $L$ injects into $L_c$. The Lie
algebra $\l$ of $L$ resp. $L_c$ decomposes as
$l=\a\oplus\m_0\oplus\m_1\oplus\m_2$ and so there is an isogeny
$A_c\times M_{0,c}\times M_{1,c}\times M_{2,c}\ra L_c$. The group
$A$ injects into $A_c$ and we have $A_c=A\times U$ where $U$ is a
compact one dimensional torus. Let $\hat{M_1}\subset M$ be the
preimage of $\overline{M_1}$ and define $\hat{M_2}$ analogously.
Then $\hat{M_1}$ and $\hat{M_2}$ are subgroups of $L_c$ and we may
define $M_1$ to be the preimage of $\hat{M_1}$ in $U\times
M_{0,c}\times M_{1,c}$ and $M_2$ to be the preimage of $\hat{M_2}$
in $M_{2,c}$. Then the map $M_1\times M_2\ra M$ is a central
isogeny.

By definition we already have $\p_{M_1}\cong [\p_M
,\g_{\alpha_r}]$. Now let $X\in\p_{M_2}$, then $[X,Y]
=[X,Y_n]+\theta([X,Y_n]) =0$ since $Y_n\in\g_{\alpha_r}$. By
$[X,H]=0$ this implies $[X,Y^c]=0$ and so
$[Pr_\pm(X),H]=0=[Pr_\pm(X),Y]=[Pr_\pm(X),Y^c]$ for $Pr_\pm$
denoting the projection $\p\ra\p_\pm$. Hence
$Pr_\pm(X)\in\p_{M_2}$ which implies
$$
\p_{M_2}\= \p_{M_2 ,+}\oplus\p_{M_2 ,-}.
$$
The lemma is proven.
\qed

Let
$$
\n_\ar := \bigoplus_{{\alpha\in\Phi(\h ,\g)}\ {\alpha |_\a =\ar
|_\a}}\g_\ar
$$
and define $\n_{\rez{2}\ar}$ and $\n_{\frac{3}{2}\ar}$
analogously. Let $\n^\ar := \n_{\rez{2}\ar}\oplus
\n_{\frac{3}{2}\ar}$, then $\n=\n_\ar\oplus\n^\ar$.

\begin{proposition}\label{zerlp}
There is a $K_M$-stable subspace $\n_-$ of $\n^{\ar}$ such that
as $K_M$-module $\n^\ar \cong \n_-\oplus\tilde{\n_-}$, where
$\tilde{.}$ denotes the contragredient, and as $K_M$-module:
$$
\p_-\ \cong\ \C \oplus \p_{M_1}\oplus\p_{M_2 ,-}\oplus\n_- .
$$
\end{proposition}

\prf The map $a+p+n\mapsto a+p+n-\theta(n)$ induces a
$K_M$-isomorphism $$ \a\oplus\p_M\oplus\n \ \cong\ \p . $$ We now
prove that there is a $K_M$-isomorphism $\n_\ar \cong \g_\ar\oplus
[\p_M ,\g_\ar]$. For this let $\alpha \ne\ar$ be any root in
$\Phi(\h,\g)$ such that $\alpha |_\a =\ar |_\a$. Then $B(\alpha
,\ar)>0$ and hence $\beta :=\alpha -\ar$ is a root. The root
$\beta$ is imaginary. Assume $\beta$ is compact, then
$\g_\beta\subset\k_M$ and we have $[\g_\ar ,\g_\beta ]=\g_\ar$
which contradicts $[\k_M ,\g_\ar ]=0$. It follows that $\beta$ is
noncompact and so $\g_\beta\subset\p_M$, which implies $\n_\ar
=\g_\ar\oplus [\p_M,\g_\ar ]$.

Now the last lemma implies the assertion.
\qed

We have $M\leftarrow M_1\times M_2$ and so any irreducible
representation $\xi$ of $M$ can be pulled back to $M_1\times M_2$
and be written as a tensor product $\xi =\xi_1\otimes\xi_2$. Let
$\xi$ be irreducible admissible then Proposition \ref{trace1} and
Proposition \ref{zerlp} imply
$$
\tr\pi_{\xi ,\nu}(f_t^0) \= e^{t\pi_{\xi ,\nu}(C)}\sum_{q\ge 0}
(-1)^q \dim\left(
W_\xi\otimes\wedge^q(\p_{M_1}\oplus\p_{M_2,-}\oplus\n_-)\right)^{K_M}.
$$
For $c\ge 0$ let
$$
\wedge^c\n_- \= \bigoplus_{i\in I_c}\sigma_{i}^c\otimes\tau_{i}^c
$$\index{$\sigma_i^c$}\index{$\tau_i^c$}
be a decomposition as $K_M\leftarrow K_{M_1}\times
K_{M_2}$-module of $\n_-$ where $\tau_i^c$ is an irreducible
$K_{M_2}$-module and $\sigma_i^c$ is the image of the projection
$\n^\ar =\n_-\oplus\tilde{\n_-}\ra\n_-$ of an irreducible
$M_1$-submodule of $\n^\ar$. Then we conclude $$ \tr\pi_{\xi
,\nu}(f_t^0) \= e^{t\pi_{\xi ,\nu}(C)}\sum_{a,b,c\ge 0}
\sum_{i\in I_c} (-1)^{a+b+c} \dim\left( W_{\xi_1}\otimes
\wedge^a\p_{M_1}\otimes\sigma_{i}^c\right)^{K_{M_1}} $$ $$
\hspace{190pt}\times\dim\left( W_{\xi_2}\otimes
\wedge^b\p_{M_2,-}\otimes\tau_{i}^c\right)^{K_{M_2}}. $$

Fix an irreducible representation $\xi =\xi_1\otimes\xi_2$ of $M$.
Let $\xi':=\tilde{\xi_1}\otimes\xi_2$, where $\tilde{\xi_1}$ is
the contragredient representation to $\xi_1$.

We will use the following notation: For an irreducible
representation $\pi$ we denote its infinitesimal character by
$\wedge_\pi$. We will identify $\wedge_\pi$ to a corresponding
element in the dual of a Cartan subalgebra (modulo the Weyl
group), so that it makes sense to write an expression like
$B(\wedge_\pi)$. In the case of $\sigma_i^c$ we have the
following situation: either $\sigma_i^c$ already extends to a
representation of $M_1$ or $\sigma_i^c\oplus\tilde{\sigma_i^c}$
does. In either case we write $\wedge_{\sigma_i^c}$ for the
corresponding infinitesimal character with respect to the group
$M_1$.

\begin{lemma}
If $\xi\cong\xi'$ and $\tr\pi_{\xi ,\nu}(f_t^0)\ne 0$ for some
$t>0$ then
$$
\pi_{\xi ,\nu}(C) \=
B(\wedge_{\sigma_i^c})+B(\wedge_{\tau_{i}^{c}\otimes\epsilon})+B(\nu)-B(\rho)
$$
for some $c\ge 0$ and some $i\in I_c$. Here $\epsilon$ is the one
dimensional representation of $K_{M_2}$ such that
$\epsilon\otimes\epsilon\cong \wedge^{top}\p_{M_2}$.
\end{lemma}

\prf Generally we have
\begin{eqnarray*}
\pi_{\xi ,\nu}(C) &=& B(\wedge_{\pi_{\xi ,\nu}})-B(\rho)\\
&=& B(\wedge_\xi)+B(\nu) -B(\rho)\\
&=& B(\wedge_{\xi_1})+B(\wedge_{\xi_2}) +B(\nu) -B(\rho).
\end{eqnarray*}

Now $\tr\pi_{\xi ,\nu}(f_t^0)\ne 0$ implies firstly
$$
\sum_{a\ge 0} (-1)^a \dim\left( W_{\xi_1}\otimes
\wedge^a\p_{M_1}\otimes\sigma_i^c\right)^{K_{M_1}}\ne 0
$$
for some $c,i$. Now either $\sigma_i^c$ already extends to an
irreducible representation of $M_1$ or
$\sigma_i^c\oplus\tilde{\sigma_i^c}$ does. In the second case we
get
$$
\sum_{a\ge 0} (-1)^a \dim\left( W_{\xi_1}\otimes
\wedge^a\p_{M_1}\otimes(\sigma_i^c\oplus\tilde{\sigma_i^c})\right)^{K_{M_1}}
$$ $$
= 2\sum_{a\ge 0} (-1)^a \dim\left( W_{\xi_1}\otimes
\wedge^a\p_{M_1}\otimes\sigma_i^c\right)^{K_{M_1}}\ne 0.
$$
So in either case Lemma \ref{Caseig1} implies
$$
B(\wedge_{\xi_1}) -B(\rho_{M_1}) \= \xi_1(C_M) \= \sigma_i^c(C_M)
\= B(\wedge_{\sigma_i^c})-B(\rho_M),
$$
hence $B(\wedge_{\xi_1})=B(\wedge_{\sigma_1^c})$.

Next $\tr\pi_{\xi ,\nu}(f_t^0)\ne 0$ implies
$$
\sum_{b\ge 0}(-1)^b\dim\left(
W_{\xi_2}\otimes\wedge^b\p_{M_2,-}\otimes\tau_{i}^{c}\right)^{K_{M_2}}\ne
0
$$
for the same $c,i$. In that case with $\tau =\tau_{i}^{c}$ Lemma
\ref{Caseig2} implies
$$
\xi_2(C_{M_2}) = \tau\otimes\epsilon(C_{K_{M_2}})-B(\rho_{M_2})
+B(\rho_{K_{M_2}}).
$$
Since $\xi_2(C_{M_2})=B(\wedge_{\xi_2})-B(\rho_{M_2})$ we conclude
$$
B(\wedge_{\xi_2}) \=
\tau\otimes\epsilon(C_{K_{M_2}})+B(\rho_{K_{M_2}}) \=
B(\wedge_{\tau\otimes\epsilon}).
$$
\qed

Along the same lines we get

\begin{lemma}
If $\xi$ is not isomorphic to $\xi'$ and $\tr\pi_{\xi ,\nu}(f_t^0)
+\tr\pi_{\xi' ,\nu}(f_t^0)\ne 0$ for some $t>0$ then
$$
\pi_{\xi ,\nu}(C) \=
B(\wedge_{\sigma_i^c})+B(\wedge_{\tau_{i}^{c}\otimes\epsilon})+B(\nu)-B(\rho)
$$
for some $c\ge 0$ and some $i\in I_c$.
\end{lemma}

We abbreviate $s_i^c := B(\wedge_{\sigma_i^c})+
B(\wedge_{\tau_i^c\otimes\epsilon})$.\index{$s_i^c$} We have shown
that if $\xi\cong\xi'$ then
 $$
 \tr\pi_{\xi ,\nu}(f_t^0)
 $$
equals
 $$ e^{t(B(\nu)-B(\rho))}\sum_{a,b,c\ge 0} \sum_{i\in I_c}
e^{t(s_i^c)}(-1)^{a+b+c} \dim\left( W_{\xi_1}\otimes
\wedge^a\p_{M_1}\otimes\sigma_{i}^c\right)^{K_{M_1}}
 $$ $$
\hspace{190pt}\times\dim\left( W_{\xi_2}\otimes
\wedge^b\p_{M_2,-}\otimes\tau_{i}^c\right)^{K_{M_2}}.
 $$
Define the Fourier transform of $f_t^0$ by
 $$
 \hat{f_t^0}_H(\nu,b^*)\= \tr \pi_{\xi_{b^*},\nu}(f_t^0).
 $$
Let $B^*$ be the character group of $B$. According to $M=M_1\times
M_2$ we can write $B=B_1\times B_2$ and so we see that any
character $b^*$ of $B$ decomposes as $b_1^*\times b_2^*$. Let
${b^*}':=\overline{b^*_1}\times b_2^*$. We get

\begin{lemma}
The sum of Fourier transforms
 $$
\hat{f_t^0}_H(\nu ,b^*)+\hat{f_t^0}_H(\nu ,{b^*}')
 $$
equals
 $$
 e^{t(B(\nu)-B(\rho))}\sum_{a,b,c\ge 0} \sum_{i\in I_c}
e^{ts_i^c}(-1)^{a+b+c} \dim\left( (V_{b_1^*}\oplus
V_{\overline{b_1^*}})\otimes
\wedge^a\p_{M_1}\otimes\sigma_{i}^c\right)^{K_{M_1}}
 $$ $$
\hspace{190pt}\times\dim\left( V_{b_2^*}\otimes
\wedge^b\p_{M_2,-}\otimes\tau_{i}^c\right)^{K_{M_2}}.
 $$
\end{lemma}

\section{The global trace of the heat kernel}

Fix a discrete cocompact torsion-free subgroup $\Ga$ of $G$ then
the quotient $X_\Ga :=\Ga\bs X = \Ga\bs G/K$ is a non-euclidean
Hermitian locally symmetric space.

Since $X$ is contractible and $\Ga$ acts freely on $X$ the group
$\Ga$ equals the fundamental group of $X_\Ga$. Let $(\omega
,V_\omega)$ denote a finite dimensional unitary representation of
$\Ga$.

D. Barbasch and H. Moscovici have shown in \cite{BM} that the
function $f_t^0$ satisfies the conditions to be plugged into the
trace formula. So in order to compute $\tr R_{\Ga ,\omega}(f_t^0)$
we have to compute the orbital integrals of $f_t^0$. At first let
$h\in G$ be a nonelliptic semisimple element. Since the trace of
$f_t^0$ vanishes on principal series representations which do not
come from splitrank one Cartan subgroups, we see that
$\CO_h(f_t^0)=0$ unless $h\in H$, where $H$ is a splitrank one
Cartan. Write $H=AB$ as before. We have $h=a_hb_h$ and since $h$
is nonelliptic it follows that $a_h$ is regular in $A$. We say
that $h$ is \emph{split regular} in $H$. We can choose a parabolic
$P=MAN$ such that $a_h$ lies in the negative Weyl chamber $A^-$.

Let $V$ denote a finite dimensional complex vector space and let
$A$ be an endomorphism of $V$. Let $\det(A)$ denote the
determinant of $A$, which is the product over all eigenvalues of
$A$ with algebraic multiplicities. Let $\det'(A)$ be the product
of all nonzero eigenvalues with algebraic multiplicities.

In \cite{HC-HA1}, sec 17 Harish-Chandra has shown that for $h_0\in
H$
 $$
\O_{h_0}(f_t^0) \= \frac{\varpi_{h_0} (\
'F_{f_t^0}^H(h))|_{h=h_0}}
                {c_{h_0} h_0^{\rho_P} \det'(1-h_0^{-1}|(\g /\h)^+)},
 $$
where $(\g /\h)^+$ is the positive part of the root space
decomposition of a compatible ordering. Further $\varpi_h$ is the
differential operator attached to $h$ as follows. Let $\g_h$
denote the centralizer of $h$ in $\g$ and let $\Phi^+(\g_h,\h)$
the positive roots then
 $$
 \varpi_h\= \prod_{\alpha \in\Phi^+(\g_h,\h)} H_\alpha,
 $$
where $H_\alpha$ is the element in $\h$ dual to $\alpha$ via the
bilinear form $B$. In comparison to other sources the formula
above for the orbital integral lacks a factor $[G_h:G_h^0]$ which
doesn't occur because of the choices of Haar measures made.  We
assume the ordering to come from an ordering of $\Phi(\b,\m)$
which is such that for the root space decomposition $\p_M
=\p_M^+\oplus \p_M^-$ it holds $\p_{M_2}\cap\p_M^+ =\p_{M_2,-}$.
For short we will henceforth write $\varpi_{h_0}F(h_0)$ instead of
$\varpi_{h_0}F(h)|_{h=h_0}$ for any function $F$.

Our results on the Fourier transform of $f_t^0$ together with a
computation as in the proof of Lemma \ref{trace-tau} imply that
$'F_{f_t^0}^H(h)$ equals
$$
\frac{e^{-l_h^2/4t}}{\sqrt{4\pi t}}
 e^{-tB(\rho)}
\det(1-h^{-1} | \k_M^+\oplus \p_{M_1}^+) \sum_{c\ge 0,\ i\in I_c}
(-1)^c e^{ts_i^c} \tr(b_h^{-1} | \sigma_i^c\oplus \tau_i^c).
$$
Note that this computation involves a summation over $B^*$ and
thus we may replace $\sigma_i^c$ by its dual $\breve{\sigma_i^c}$
without changing the result.

Now let $(\tau ,V_\tau)$ be a finite dimensional unitary
representation of $K_M$ and define for $b\in B$ the monodromy
factor: $$ L^M(b,\tau) \= \frac{\varpi_b( \det(1-b|(\k_M/|\b)^+)
\tr(\tau(b)))}
                {\varpi_b( \det(1-b|(\m/\b)^+))}.
$$

Note that the expression $\varpi_b( \det (1-b|(\m /\b)^+)) $
equals
$$
|W(\m_b ,\b)|\prod_{\alpha \in \Phi_b^+(\m ,\b)}(\rho_b,\alpha)
\det'(1-b|(\m /\b)^+),
$$
and that for $\ga =a_\ga b_\ga$ split-regular we have
$$|W(\m_b,\b)|\prod_{\alpha \in \Phi_b^+(\m ,\b)}(\rho_b,\alpha)
\= |W(\g_\alpha ,\h)|\prod_{\alpha \in \Phi_\ga^+}(\rho_\ga
,\alpha),
$$
so that writing $L^M(\ga ,\tau) = L^M(b_\ga
  ,\tau)$ we get
$$
L^M(\ga ,\tau) \= \frac{\varpi_\ga( \det (1-\ga |
  (\k_M/\b)^+) \tr(\tau(b_\ga)))}
       {|W(\g_\ga ,\h)|\prod_{\alpha \in \Phi_\ga^+}(\rho_\ga
                        ,\alpha) \det'(1-\ga | (\m /\b)^+)}.
$$

From the above it follows that $\vol(\Ga_\ga \bs G_\ga) \O_\ga
(f_t^0)$ equals:
$$
\frac{\chi_{_1}(X_\ga)l_{\ga_0}}
     {\det(1-\ga | \n)}\
\frac{e^{-l_\ga^2/4t}}
     {\sqrt{4\pi t}}\
a_\ga^{\rho_P} \sum_{c\ge 0,\ i\in I_c} (-1)^c\ e^{ts_i^c}\
\tr(b_\ga|\sigma_i^c)L^{M_2}(\ga ,\tau_i^c).
$$
Note that here we can replace $\sigma_i^c$ by
$\tilde{\sigma_i^c}$.

For a splitrank one Cartan $H=AB$ and a parabolic $P=MAN$ let
$\CE_P(\Ga)$ denote the set of $\Ga$-conjugacy classes
$[\ga]\subset\Ga$ such that $\ga$ is $G$-conjugate to an element
of $A^-B$, where $A^-$ is the negative Weyl chamber given by $P$.
We have proven:

\begin{theorem} \label{higher_heat_trace} Let $X_\Ga$ be a compact locally
Hermitian space with fundamental group $\Ga$ and such that the
universal covering is globally symmetric without compact factors.
Assume $\Ga$ is neat and write $\triangle_{p,q,\omega}$ for the
Hodge Laplacian on $(p,q)$-forms with values in a the flat
Hermitian bundle $E_\omega$, then it theta series defined by
 $$
 \Theta (t) \= \sum_{q=0}^{\dim_\C X_\Ga}
q(-1)^{q+1} \tr\ e^{-t\triangle_{0,q,\omega}}
 $$
equals
 $$
\sum_{P/conj.} \sum_{[\ga] \in \CE_P(\Ga)} \chi_{_1}(X_\ga)
\frac{l_{\ga_0}a_\ga^{\rho_P}} {\det(1-\ga | \n)}\
\frac{e^{-l_\ga^2/4t}} {\sqrt{4\pi t}}\
 $$ $$
\hspace{110pt}\times\sum_{c\ge 0,\ i\in I_c} (-1)^c\ e^{ts_i^c}\
\tr(b_\ga |\sigma_i^c)L^{M_2}(\ga ,\tau_i^c).
 $$ $$
        + f^0_t(e)\ \dim\omega\ vol(X_\Ga).\hspace{100pt}
 $$
\end{theorem}

The reader should keep in mind that by its definition we have for
the term of the identity:
$$
f^0_t(e)\ \dim\omega\ vol(X_\Ga) \= \sum_{q=0}^{\dim_\C
X}q(-1)^{q+1} \tr_\Ga(e^{-t\lap_{0,q,\omega}}),
$$
where $\tr_\Ga$ is the $\Ga$-trace. Further note that by the
Plancherel theorem the Novikov-Shubin invariants of all operators
$\lap_{0,q}$ are positive.

\section{The holomorphic torsion zeta function}

Now let $X$ be Hermitian again and let $\tau =\tau_1\otimes\tau_2$
be an irreducible representation of $K_{M}\leftarrow K_{M_1}\times
K_{M_2}$.

\begin{theorem}\label{hol_zeta}
Let $\Ga$ be neat and $(\omega ,V_\omega)$ a finite dimensional
unitary representation of $\Ga$. Choose a $\theta$-stable Cartan
$H$ of splitrank one. For $\Re(s)>>0$ define the zeta function
$Z_{P,\tau,\omega}^0(s)$ to be
$$
\exp \left( -\sum_{[\ga] \in \CE_H(\ga)} \frac{\chi_{_1}(X_\ga)
\tr(\omega(\ga)) \tr\tau_1((b_\ga)L^{M_2}(\ga ,\tau_2)}
     {\det (1-\ga|\n)}\
\frac{e^{-sl_\ga}}{\mu_\ga})\right).
$$
Then $Z_{P,\tau,\omega}^0$ has a meromorphic continuation to the
entire plane. The vanishing order of $Z_{P,\tau,\omega}^0(s)$ at
$s=\la (H_1)$, $\la \in \a^*$ is $(-1)^{\dim \n}$ times
$$
\sum_{\pi \in \hat{G}}N_{\Ga ,\omega}(\pi)
    \sum_{p,q,r}(-1)^{p+q+r} \dim \Big( H^q({\n},\pi_K)^\la
    \otimes \wedge^p\p_{M_1}\otimes\wedge^r\p_{M_2,-}
    \otimes V_{\breve{\tau}}\Big)^{K_M}.
$$
Note that
in the special case $M_2={1}$ this function equals the Selberg
zeta function.
\end{theorem}

\prf By Lemma \ref{orient} the group $M$ is orientation preserving
an so are $M_1$ and $M_2$. Therefore the Euler-Poincar\'e function
$f_{\sigma}^{M_1}$ for the representation $\sigma$ exists. Further
for $M_2$ the function $g_\tau^{M_2}$ of Theorem \ref{exist-ep}
exists. We set $h_{\sigma\otimes\tau}(m_1,m_2) :=
f_{\sigma}^{M_1}(m_1)g_\tau^{M_2}(m_2)$ and this function factors
over $M$. Then for any function $\eta$ on $M$ which is a product
$\eta =\eta_1\otimes\eta_2$ on $M_1\times M_2$ we have for the
orbital integrals:
$$
\CO_m^M(\eta) \= \CO_{m_1}^{M_1}(\eta_1)\CO_{m_2}^{{M_2}}(\eta_2).
$$
With this in mind it is straightforward to see that the proof of
Theorem \ref{hol_zeta} proceeds as the proof of Theorem
\ref{genSelberg} with the Euler-Poincar\'e function $f_\tau$
replaced by the function $h_{\sigma\otimes\tau}$.
\qed

Extend the definition of $Z_{P,\tau,\omega}^0(s)$ to arbitrary
virtual representations in the following way. Consider a finite
dimensional virtual representation $\xi = \oplus_i a_i \tau_i$
with $a_i\in \Z$ and $\tau_i\in \hat{K_M}$. Then let
$Z_{P,\xi,\omega}^0(s) = \prod_i Z_{P,\tau_i,\omega}^0(s)^{a_i}.$

\begin{theorem} \label{detformel}
Assume $\Ga$ is neat, then for $\la >>0$ we have the identity
$$
\prod_{q=0}^{\dim_\C X} \left(\frac{\det
(\lap_{0,q,\omega}+\la)}{{\det}^{(2)}
(\lap_{0,q,\omega}+\la)}\right)^{q(-1)^{q+1}}
            \=
$$ $$
\prod_{P/{\rm conj.}} \prod_{{c\ge 0}\ {i\in I_c}}
Z_{P,\sigma_i^c\otimes\tau_i^c ,\omega}(|\rho_P|+\sqrt{\la +
s_i^c})^{(-1)^c}
$$
\end{theorem}

\prf Consider Theorem \ref{higher_heat_trace}. For any
semipositive elliptic differential operator $D_\Ga$ the heat
trace $\tr e^{-tD_\Ga}$ has the same asymptotic as $t\ra 0$ as
the $L^2$-heat trace $\tr_\Ga e^{-tD}$. Thus it follows that the
function
$$
h(t):= \sum_{q=0}^{\dim_\C X_\Ga} q(-1)^{q+1} (\tr
e^{-t\lap_{0,q,\omega}}-\tr_\Ga e^{-t\lap_{0,q,\omega}})
$$
is rapidly decreasing at $t=0$. Therefore, for $\la >0$ the Mellin
transform of $h(t)e^{-t\la}$ converges for any value of $s$ and
gives an entire function. Let
$$
\zeta_\la (s) := \rez{\Ga (s)} \int_0^\infty t^{s-1} h(t)
e^{-t\la} dt.
$$
We get that
$$
\exp(-\zeta_\la'(0)) \= \prod_{q=0}^{\dim_\C X} \left(\frac{\det
(\lap_{0,q,\omega}+\la)}{{\det}^{(2)}
(\lap_{0,q,\omega}+\la)}\right)^{q(-1)^{q+1}}.
$$
On the other hand, Theorem \ref{higher_heat_trace} gives a second
expression for $\zeta_\la(s)$. In this second expression we are
allowed to interchange integration and summation for $\la >>0$
since we already know the convergence of the Euler products
giving the right hand side of our claim.
\qed

Let $n_0$ be the order at $\la =0$ of the left hand side of the
last proposition. Then
$$
n_0 \= \sum_{q=0}^{\dim_\C(X)} q(-1)^q(h_{0,q,\omega}
-h_{0,q,\omega}^{(2)}),
$$
where $h_{0,q,\omega}$ is the $(0,q)$-th Hodge number of $X_\Ga$
with respect to $\omega$ and $h_{0,q,\omega}^{(2)}$ is the
$L^2$-analogue. Conjecturally we have $h_{0,q,\omega}^{(2)} =
h_{0,q,\omega}$, so $n_0=0$. For a splitrank one Cartan $H$, for
$c\ge 0$ and $i\in I_c$ we let
$$
n_{P,c,i,\omega} \ := \ {\rm
ord}_{s=|\rho_P|+\sqrt{s_i^c}}Z_{P,\sigma_i^c\otimes\tau_i^c
,\omega}(s)
$$
so $n_{P,c,i,\omega}$ equals $(-1)^{\dim \n}$ times
$$
\sum_{\pi \in \hat{G}}N_{\Ga ,\omega}(\pi)
    \sum_{p,q,r}(-1)^{p+q+r} \dim \Big( H^q({\n},\pi_K)^\la\otimes
\wedge^p\p_{M_1}\otimes\wedge^r\p_{M_2,-} \otimes
V_{\breve{\sigma}\otimes\breve{\tau}}\Big)^{K_M},
$$
for $\la(H)=|\rho_P|+\sqrt{s_i^c}$.
 We then consider
$$
c(X_\Ga ,\omega)\= \prod_{P}\prod_{c\ge 0,\ i\in I_c} \left(
2\sqrt{s_i^c}\right)^{(-1)^cn_{P,c,i,\omega}}.
$$

We assemble the results of this section to

\begin{theorem}
Let
$$
Z_\omega(s) \= \prod_{H/{\rm conj.}} \prod_{{c\ge 0}\ {i\in
I_c}}Z_{H,\sigma_i^c\otimes\tau_i^c,\omega}\left(
s+|\rho_P|+\sqrt{s_i^c}\right),
$$
then $Z_\omega$ extends to a meromorphic function on the plane.
Let $n_0$ be the order of $Z_\omega$ at zero then
$$
n_0 \= \sum_{q=0}^{\dim_\C(X)} q(-1)^q \left( h_{0,q}(X_\Ga)
-h_{0,q}^{(2)}(X_\Ga)\right),
$$
where $h_{p,q}(X_\Ga)$ is the $(p,q)$-th Hodge number of $X_\Ga$
and $h_{p,q}^{(2)}(X_\Ga)$ is the $(p,q)$-th $L^2$-Hodge number of
$X_\Ga$. Let $R_\omega(s) = Z_\omega(s)s^{-n_0}/c(X_\Ga ,\omega)$
then
$$
R_\omega(0) \= \frac{T_{hol}(X_\Ga
,\omega)}{T_{hol}^{(2)}(X_\Ga)^{\dim\omega}}.
$$
\end{theorem}

\prf This follows from Theorem \ref{detformel}. \qed

\chapter{The Lefschetz formula for a $p$-adic group}

\section{The trace formula}
Let $F$ be a nonarchimedean local field with valuation
ring $\CO$ and uniformizer $\varpi$. Denote by $G$ be a
semisimple linear algebraic group over $F$. Let $K\subset G$ be a
good maximal compact subgroup. Choose a parabolic
subgroup $P=LN$ of $G$ with Levi component $L$. Let $A$
denote the greatest split torus in the center of $L$. Then
$A$ is called the {\it split component} of $P$. Let
$\Phi=\Phi(G,A)$ be the root system of the pair $(G,A)$,
i.e. $\Phi$ consists of all homomorphisms $\alpha :A\ra
\GL_1$ such that there is $X$ in the Lie algebra of $G$
with $\Ad(a)X= a^\alpha X$ for every $a\in A$. Given
$\alpha$, let $\n_\alpha$ be the Lie algebra generated by
all such $X$ and let $N_\alpha$ be the closed subgroup of
$N$ corresponding to $\n_\alpha$. Let $\Phi^+=\Phi(P,A)$ be
the subset of $\Phi$ consisting of all positive roots with
respect to $P$. Let $\Delta\subset\Phi^+$ be the subset of
simple roots. 
 Let $A^-\subset A$ be the set of all
$a\in A$ such that $|a^\alpha| <1$ for any $\alpha\in
\Delta$.

There is a reductive subgroup $M$ of $L$ with compact center such
that $MA$ has finite index in $L$. We can assume that $K_M=M\cap K$ is a good
maximal subgroup of $M$. An element $g$ of $G$ is called \emph{elliptic} if it is
contained in a compact subgroup. Let $M_{ell}$ denote the set of elliptic elements
in $M$.

Let $X^*(A)=\Hom(A,\GL_1)$ be the group of all
homomorphisms as algebraic groups from $A$ to $\GL_1$. This
group is isomorphic to $\Z^r$ with $r=\dim A$. Likewise let
$X_*(A)=\Hom(\GL_1,A)$. There is a natural $\Z$-valued
pairing
\begin{eqnarray*}
X^*(A)\times X_*(A) &\ra& \Hom(\GL_1,\GL_1)\cong\Z\\
(\alpha,\eta) &\mapsto& \alpha\circ\eta.
\end{eqnarray*}
For every root $\alpha\in\Phi(A,G)\subset X^*(A)$ let
$\breve{\alpha}\in X_*(A)$ be its coroot. Then
$(\alpha,\breve{\alpha})=2$. The valuation $v$ of $F$
gives a group homomorphism $\GL_1(F)\ra\Z$. Let $A_c$ be
the unique maximal compact subgroup of $A$. Let $\Sigma
=A/A_c$; then $\Sigma$ is a $\Z$-lattice of rank $r=\dim
A$. By composing with the valuation $v$ the group $X^*(A)$
can be identified with
$$
\Sigma^*\=\Hom(\Sigma,\Z).
$$
Let
$$
\a_0^*\=\Hom(\Sigma,\R)\ \cong\ X^*(A)\otimes\R
$$
be the real vector space of all group homomorphisms from
$\Sigma$ to $\R$ and let
$\a^*=\a_0^*\otimes\C=\Hom(\Sigma,\C)\cong
X^*(A)\otimes\C$. For $a\in A$ and $\la\in\a^*$ let
$$
a^\la\=q^{-\la(a)},
$$
where $q$ is the number of elements in the residue class field of
$F$. In this way we get an identification
$$
{\a^*}/\mbox{\small$\frac{2\pi i}{\log q}$}\Sigma^* \
\cong\ \Hom(\Sigma,\C^\times).
$$
A quasicharacter $\nu : A\ra\C^\times$ is called {\it
unramified} if $\nu$ is trivial on $A_c$. The set
$\Hom(\Sigma,\C^\times)$ can be identified with the set of
unramified quasicharacters on $A$. Any unramified
quasicharacter $\nu$ can thus be given a unique real part
$\Re(\nu)\in \a_0^*$. This definition extends to not
necessarily unramified quasicharacters $\chi:A\ra\C^\times$
as follows. Choose a splitting $s:\Sigma\ra A$ of the exact
sequence
$$
1\ra A_c\ra A\ra\Sigma\ra 1.
$$
Then $\nu=\chi\circ s$ is an unramified character of $A$. Set
$$
\Re(\chi)\=\Re(\nu).
$$
This definition does not depend on the choice of the splitting
$s$. For quasicharacters $\chi$, $\chi'$ and $a\in A$ we will
frequently write $a^\chi$ instead of $\chi(a)$ and
$a^{\chi+\chi'}$ instead of $\chi(a)\chi'(a)$. Note that the
absolute value satisfies $|a^\chi|=a^{\Re(\chi)}$ and that a
quasicharacter $\chi$ actually is a character if and only if
$\Re(\chi)=0$.

Let $\Delta_P : P\ra\R_+$ be the modular function of the
group $P$. Then there is $\rho\in\a_0^*$ such that
$\Delta_P(a)=|a^{2\rho}|$. For $\nu\in\a^*$ and a root
$\alpha$ let
$$
\nu_\alpha\= (\nu,\breve{\alpha})\ \in\
X^*(\GL_1)\otimes\C\ \cong\ \C.
$$
Note that $\nu\in\a_0^*$ implies $\nu_\alpha\in\R$ for
every $\alpha$. For $\nu\in\a_0^*$ we say that $\nu$ is
positive, $\nu>0$, if $\nu_\alpha>0$ for every positive
root $\alpha$.

{\bf Example.} Let $G=\GL_n(F)$ and let $\varpi_j\in G$ be
the diagonal matrix
$\varpi_j=\diag(1,\dots,1,\varpi,1,\dots,1)$ with the
$\varpi$ on the $j$-th position. Let $\nu\in\a^*$ and let
$$
\nu_j\=\nu(\varpi_j A_c)\ \in\ \C.
$$
Let $\alpha$ be a root, say
$\alpha(\diag(a_1,\dots,a_n))=\frac{a_i}{a_j}$. Then
$$
\nu_\alpha\= \nu_i-\nu_j.
$$
Hence $\nu\in\a_0^*$ is positive if and only if
$\nu_1>\nu_2>\dots >\nu_n$.

We will fix Haar-measures of $G$ and its reductive subgroups as follows.
For $H\subset G$ being a torus there is a maximal compact subgroup $U_H$ which is open and unique up to conjugation. 
Then we fix a Haar measure on $H$ such that $\vol (U_H)=1$.
If $H$ is connected reductive with compact center then we choose the unique positive Haar-measure which up to sign coincides with the Euler-Poincar\'e measure
\cite{Kottwitz}.
So in the latter case our measure is determined by the following property:
For any discrete torsionfree cocompact subgroup $\Ga_H\subset H$ we have
$$
\vol (\Ga_H \bs H) \= (-1)^{q(H)}\chi(\Ga_H ,\Q),
$$
where $q(H)$ is the $k$-rank of the derived group $H_{der}$ and $\chi(\Ga_H ,\Q)$ the
Euler-Poincar\'e characteristic of $H^\bullet(\Ga_H ,\Q)$. For the applications
recall that centralizers of tori in connected groups are connected
\cite{Borel-lingroups}.

Assume we are given a discrete subgroup $\Ga$ of $G$ such that the quotient space
$\Ga \bs G$ is compact. Let $(\omega ,V_\omega)$ be a finite dimensional unitary
representation of $\Ga$ and let $L^2(\Ga \bs G,\omega)$ be the Hilbert space
consisting of all measurable functions $f: G \ra V_\omega$ such that $f(\ga x) =
\omega(\ga) f(x)$ and $|f|$ is square integrable over $\Ga \bs G$ (modulo null
functions). Let $R$ denote the unitary representation of $G$ on $L^2(\Ga \bs
G,\omega)$ defined by right shifts, i.e. $R(g) \ph (x) = \ph (xg)$ for $\ph \in
L^2(\Ga \bs G,\omega)$. It is known that as a $G$-representation this space splits
as a topological direct sum:
$$
L^2(\Ga \bs G,\omega) \= \bigoplus_{\pi \in \hat{G}} N_{\Ga ,\omega}(\pi) \pi
$$
with finite multiplicities $N_{\Ga ,\omega}(\pi)<\infty$.

Let $f$ be integrable over $G$, so $f$ is in $L^1(G)$.
The integral
$$
R(f) := \int_G f(x) R(x) \ dx
$$
defines an operator on the Hilbert space $L^2(\Ga \bs G,\omega)$.

For $g\in G$ and $f$ any function on $G$ we define the \emph{orbital
integral}\index{orbital integral}
$$
\CO_g(f) := \int_{G_g\bs G} f(x^{-1} gx)\ dx,
$$
whenever the integral exists.
Here $G_g$ is the centralizer of $g$ in $G$.
It is known that the group $G_g$ is unimodular, so we have an invariant measure on
$G_g\bs G$.

A function $f$ on $G$ or any of its closed subgroups is called
\emph{smooth} if it is locally constant. It is called \emph{uniformly
smooth}\index{uniformly smooth} if there is an open subgroup
$U$ of
$G$ such that $f$ factors over
$U\bs G/U$. This is in particular the case if $f$ is smooth and compactly supported.

\begin{proposition}
(Trace formula)\\
Let $f$ be integrable and uniformly smooth, then we have 
$$
\sum_{\pi \in\hat{G}} N_{\Ga ,\omega}(\pi)\ \tr \pi (f) \= \sum_{[\ga]} \tr \omega (\ga)\ \vol(\Ga_\ga \bs G_\ga)\ \CO_\ga (f),
$$
where the sum on the right hand side runs over the set of $\Ga$-conjugacy classes $[\ga]$ in $\Ga$ and $\Ga_\ga$ denotes the centralizer of $\ga$ in $\Ga$.
Both sides converge absolutely and the left hand side actually is a finite sum.
\end{proposition}

\prf
At first fix a fundamental domain $\CF$ for $\Ga \bs G$ and let $\ph\in L^2(\Ga \bs G,\omega)$, then
\begin{eqnarray*}
R(f) &\=& \int_Gf(y) \ph(xy)\ dy\\
	&\=& \int_G f(x^{-1}y) \ph(y)\ dy\\
	&\=& \sum_{\ga\in\Ga} \int_\CF f(x^{-1}\ga y)\ph(\ga y)\ dy\\
	&\=& \int_{\Ga\bs G} \left( \sum_{\ga\in\Ga} f(x^{-1}\ga y)\omega (\ga)\right) \ph (y) \ dy
\end{eqnarray*}

We want to show that the sum $\sum_{\ga\in\Ga} f(x^{-1}\ga y)\omega (\ga)$ converges in $\End(V_\omega)$ absolutely and uniformly in $x$ and $y$.
Since $y$ can be replaced by $\ga y$, $\ga\in\Ga$ and since $\omega$ is unitary, we
only have to show the convergence  of $\sum_{\ga\in\Ga}|f(x^{-1}\ga y)|$ locally
uniformly in $y$. Let $\ga$ and $\tau$ be in $\Ga$  and assume that $x^{-1}\ga y$
and $x^{-1}\tau y$ lie in the same class in $G/U$. Then it follows $\tau yU\cap\ga
yU \ne \emptyset$ so with $V=yUy^{-1}$ we have $\ga^{-1} \tau V \cap V \ne
\emptyset$. It is clear that $V$ depends on $y$ only up to $U$ so to show locally
uniform convergence in $y$ it suffices to fix $V$. Since $V$ is compact also $V^2 =
\{ vv' |v,v'\in V\}$ is compact and so $\Ga \cap V^2$ is finite. This implies that
there are only finitely many $\ga\in\Ga$ with $\ga V\cap V\ne \emptyset$. Hence the
map $\Ga \ra G/U$, $\ga \mapsto x^{-1}\ga yU$ is finite to one with fibers having
$\le n$ elements for some natural number $n$. For $y$ fixed modulo $U$ we get
\begin{eqnarray*}
\sum_{\ga\in\Ga} |f(x^{-1} \ga y)| &\ \le\ & n\int_{G/U} |f(x)|\ dx\\
	&\=& \frac{n}{\vol(U)} \parallel f\parallel^1.
\end{eqnarray*}
We have shown the uniform convergence of the sum
$$
k_f(x,y) \= \sum_{\ga\in\Ga} f(x^{-1} \ga y) \omega(\ga).
$$

Observe that $R(f)$ factors over $L^2(\Ga\bs G,\omega)^U = L^2(\Ga\bs G/U,\omega)$, which is finite dimensional since $\Ga\bs G/U$ is a finite set.
So $R(f)$ acts on a finite dimensional space and $k_f(x,y)$ is the matrix of this operator.
We infer that $R(f)$ is of trace class, its trace equals
$$
\sum_{\pi\in\hat{G}} N_{\Ga ,\omega}(\pi)\ \tr\pi(f),
$$
and the sum is finite.
Further, since $k_f(x,y)$ is the matrix of $R(f)$ this trace also equals
\begin{eqnarray*}
\int_{\Ga \bs G} \tr k_f(x,x)\ dx &\=& \sum_{\ga\in\Ga} \int_\CF f(x^{-1}\ga x) \ dx\ \tr\omega(\ga)\\
	&\=& \sum_{[\ga]} \sum_{\sigma \in \Ga_\ga \bs \Ga} \int_\CF f((\sigma x)^{-1} \ga (\sigma x))\ dx\ \tr\omega(\ga)\\
	&\=& \sum_{[\ga]}  \int_{\Ga_\ga \bs G} f(x^{-1} \ga x)\ dx\ \tr\omega(\ga)\\
	&\=& \sum_{[\ga]}  \vol(\Ga_\ga\bs G_\ga)\ \int_{G_\ga \bs G} f(x^{-1} \ga x)\ dx\
\tr\omega(\ga).
\end{eqnarray*}
\qed

\section{The covolume of a centralizer}
Suppose $\ga\in\Ga$ is $G$-conjugate to some $a_\ga m_\ga\in A^+M_{ell}$.
We want to compute the covolume
$$
\vol(\Ga_\ga \bs G_\ga) .
$$
An element $x$ of $G$ is called 
 \emph{neat}\index{neat} if for  every
representation $\rho : G\ra GL_n(F)$ of $G$ the matrix $\rho(x)$ has no eigenvalue
which is a root of unity. A subset $A$ of $G$ is called neat if each element of it
is. Every arithmetic
$\Ga$ has a finite index subgroup which is neat \cite{Bor}. 

\begin{lemma}
Let $x\in G$ be neat and semisimple. Let $G_x$ denote its centralizer in $G$. Then
for every
$k\in\N$ we have $G_x=G_{x^k}$.
\end{lemma}

\prf
Since $G$ is linear algebraic it is a subgroup of some $H=\GL_n(\bar F)$, where $\bar
F$ is an algebraic closure of $F$. If we can show the claim for $H$ then it follows
for $G$ as well since $G_x=H_x\cap G$. In $H$ we can assume $x$ to be a diagonal
matrix. Since
$x$ is neat this implies the claim.
\qed

We suppose that $\Ga$ is
neat. This implies that for any
$\ga\in\Ga$ the Zariski closure of the group generated by $\ga$ is a torus. It then
follows  that $G_\ga$ is a connected reductive group
\cite{Borel-lingroups}.

An element $\ga\in\Ga$ is called \emph{primitive}\index{primitive element}
if
$\ga =\sigma^n$ with
$\sigma\in\Ga$ and $n\in\N$ implies $n=1$. It is a property of discrete cocompact
torsion free subgroups $\Ga$ of $G$ that every $\ga\in\Ga$, $\ga\ne 1$ is a positive
power of a unique primitive element. In other words, given a nontrivial $\ga\in\Ga$
there exists a unique primitive $\ga_0$ and a unique $\mu(\ga)\in\N$ such that
$$
\ga =\ga_0^{\mu(\ga)}.
$$

Let $\Sigma$ be a group of finite cohomological dimension $cd(\Sigma)$ over $\Q$.
We write
$$
\chi(\Sigma) \= \chi(\Sigma ,\Q) \ :=\ \sum_{p=0}^{cd(\Sigma)} (-1)^p \dim H^p(\Sigma
,\Q),
$$
for the Euler-Poincar\'e characteristic.
We also define the higher Euler characteristic as
$$
\chi_{_r}(\Sigma) \= \chi_{_r}(\Sigma ,\Q) \ :=\ \sum_{p=0}^{cd(\Sigma)} 
(-1)^{p+r}\binom pr
\dim H^p(\Sigma ,\Q),
$$
for $r=1,2,3,\dots$
It is known that $\Ga$ has finite cohomological dimension over $\Q$.

 We denote by $\CE_P(\Ga)$
the set of all conjugacy classes $[\ga]$ in $\ga$ such that $\ga$ is in $G$
conjugate to an element $a_\ga m_\ga\in AM$, where $m_\ga$ is elliptic and
$a_\ga\in A^-$.

Let $\ga\in\CE_P(\Ga)$. To simplify the notation let's assume that $\ga=a_\ga
m_\ga\in A^- M_{ell}$. Let $C_\ga$ be the connected component of the center of
$G_\ga$ then
$C_\ga = AB_\ga$, where
$B_\ga$ is the connected center of $M_{m_\ga}$ the latter group will also be
written as $M_\ga$. Let $M_\ga^{der}$ be the derived group of $M_\ga$. Then
$M_\ga=M_\ga^{der} B_\ga$.

\begin{lemma}
$B_\ga$ is compact.
\end{lemma}

\prf
Since $m_\ga$ is elliptic there is a compact Cartan subgroup $T$ of $M$ containing
$m_\ga$. Since $M$ modulo its center is a connected semisimple linear algebraic
group it follows that $T$ is a torus and therefore abelian. Therefore $T\subset
M_{m_\ga}$. Let
$b\in B_\ga$. Then $b$ commutes with every $t\in T$, therefore $b$ 
lies in the centralizer of $T$ in $M$ which equals $T$. So we have shown
$B_\ga\subset T$.
\qed

Let $\Ga_{\ga,A}=A\cap \Ga_\ga B_\ga$ and $\Ga_{\ga,M}=M_\ga^{der}\cap \Ga_\ga
AB_\ga$. Similar to the proof of Lemma 3.3 of \cite{Wolf} one shows that
$\Ga_{\ga,A}$ and
$\Ga_{\ga,M}$ are discrete cocompact subgroups of $A$ and $M_\ga^{der}$ resp.
Let
$$
\la_{\ga}\df \vol(\Ga_{\ga,A}\bs A)\= \vol(\Ga_\ga\bs AB_\ga).
$$

\begin{proposition}\label{2.3}
Assume $\Ga$ neat and let $\ga\in\Ga$ be $G$-conjugate to an element of $A^+M_{ell}$.
Then we get
$$
\vol (\Ga_\ga\bs G_\ga) = \la_\ga\ (-1)^{q(G)+r}\ \chi_{_r}(\Ga_\ga),
$$
where $r=\dim A$. If $\ga$ is a regular element of $G$, then $(-1)^{q(G)+1}\,\chi_{_r}(\Ga_\ga)=1$.
\end{proposition}

\prf The last sentence follows from the fact that if $\ga$ is regular, then $M_\ga=B_\ga$.
We normalize the volume of $B_\ga$ to be $1$. Then
\begin{eqnarray*}
\vol(\Ga_\ga\bs G_\ga) &=& \vol(\Ga_\ga\bs AM_\ga)\\
&=& \vol(\Ga_\ga B_\ga \bs AM_\ga)
\end{eqnarray*}
The space $\Ga_\ga B_\ga \bs AM_\ga$ is the total space of a fibration with fibre
$\Ga_{\ga,A}\bs A$ and base space $\Ga_\ga AB_\ga\bs M_\ga A\cong \Ga_{\ga,M}\bs
M_\ga^{der}$. Hence
$$
\vol(\Ga_\ga B_\ga \bs AM_\ga) \= \vol (\Ga_{\ga,A}\bs A) \,\vol(\Ga_{\ga,M}\bs
M_\ga^{der}).
$$
Since $\la_\ga=\vol (\Ga_{\ga,A}\bs A)$ it remains to show
$$
\vol(\Ga_{\ga,M}\bs M_\ga^{der})\= (-1)^r\chi_r(\Ga_\ga).
$$
We know that
$$
\vol(\Ga_{\ga,M}\bs M_\ga^{der})\= (-1)^{q(M_\ga)}\chi(\Ga_{\ga,M})\=
(-1)^{q(G)+r}\chi(\Ga_{\ga,M}).
$$
So it remains to show that $\chi(\Ga_{\ga,M})=\chi_r(\Ga_\ga)$.
The group $\Ga_{\ga,M}$ is isomorphic to $\Ga_\ga/\Sigma$, where $\Sigma=\Ga\cap
AB_\ga$ is isomorphic to $\Z^r$. So the proporsition follows from the next Lemma.

\begin{lemma}\label{chichi1}
Let $\Ga,\Lambda$ be of finite cohomological dimension over $\Q$.
Let $C_r$ be a group isomorphic to $\Z^r$ and assume there is an exact sequence
$$
1\ra C_r\ra \Ga\ra \Lambda\ra 1.
$$
Assume that $C_r$ is central in $\Ga$. 
Then 
$$
\chi (\Lambda,\Q) \= \chi_{r}(\Ga,\Q).
$$
\end{lemma}

\prf
We first consider the case $r=1$. In this case we want to prove for every $r$,
$$
\chi_{r-1} (\Lambda,\Q) \= \chi_{r}(\Ga,\Q).
$$
For this consider the Hochschild-Serre spectral sequence:
$$
E_2^{p,q} \= H^p(\Lambda,H^q(C_1,\Q))
$$
which abuts to
$$
H^{p+q}(\Ga,\Q).
$$
Since $C_1\cong \Z$ it follows
$$
H^q(C_1,\Q) \= \left\{ \begin{array}{cl} \Q& {\rm if}\ q=0,1\\ 0&{\rm
else}.\end{array}\right.
$$
Since $C_1$ is infinite cyclic and central it is an exercise to see that the
spectral sequence degenerates at $E_2$. Therefore,
\begin{eqnarray*}
\chi_r(\Ga) &=& \sum_{j\ge 0} (-1)^{j+r}\binom jr \dim H^j(\Ga)\\
&=& \sum_{j\ge r} (-1)^{j+r}\binom jr (\dim H^j(\Lambda)+\dim H^{j-1}(\Lambda)\\
&=& \sum_{j\ge r} (-1)^{j+r}\binom jr \dim H^j(\Lambda)\\
&& \quad -\sum_{j\ge{r-1}} (-1)^{j+r} \binom{j+1}r \dim H^j(\Lambda).
\end{eqnarray*}
Now replace $\binom{j+1}r$ by $\binom jr +\binom j{r-1}$ to get the claim. For the
general case write $C_r=C_1\oplus C^1$, where $C_1$ is cyclic and $C^1\cong
\Z^{r-1}$. Apply the above to $C_1$ and iterate this to get the lemma and hence the
proposition.
\qed

\section{The Lefschetz formula}
For a representation $\pi$ of $G$ let $\pi^\infty$ denote the subrepresentation of
\emph{smooth vectors}, ie $\pi^\infty$ is the representation on the space
$\bigcup_{H\subset G} \pi^H$, where $H$ ranges over the set of all open subgroups
of $G$. Further let $\pi_N$ denote the \emph{Jacquet module}\index{Jacquet
module} of
$\pi$. By definition $\pi_N$ is the largest quotient $MAN$-module of
$\pi^\infty$ on which
$N$ acts trivially. One can achieve this by factoring out the vector subspace
consisting of all vectors of the form $v-\pi(n)v$ for $v\in\pi^\infty$, $n\in N$. It
is known that if $\pi$ is an irreducible admissible representation, then $\pi_N$ is a
admissible $MA$-module of finite length. For a smooth $M$-module $V$ let
$H_c^\bullet(M,V)$ denote the continuous cohomology with coefficients in $V$ as in
\cite{BorWall}.

\begin{theorem}(Lefschetz Formula)\\
Let $\Ga$ be a neat discrete cocompact subgroup of $G$.
Let $\ph$ be a uniformly smooth function on $A$ with support in $A^-$. Suppose that
the function $a\mapsto \ph(a)|a^{-2\rho}|$ is integrable on $A$. Let
$\sigma$ be a finite dimensional representation of $M$. Let $q$ be the $F$-splitrank
of $G$ and $r=\dim A$.Then
$$
\sum_{\pi\in\hat G} N_{\Ga,\omega}(\pi)\sum_{q=0}^{\dim M} (-1)^a \int_{A^-}
\ph(a)\,\tr(a| H_c^q(M,\pi_N\otimes\sigma))\, da
$$
equals
$$
(-1)^{q+r}\sum_{[\ga]\in\CE_P(\Ga)} \la_\ga\,
\chi_r(\Ga_\ga)\,\tr\omega(\ga)\,\tr\sigma(m_\ga)\,\ph(a_\ga)\,|a_\ga^{2\rho}|.
$$
Both outer sums converge absolutely and the sum over $\pi\in\hat G$ actually is a
finite sum, ie, the summand is zero for all but finitely many $\pi$. For a given
compact open subgroup $U$ of $A$ both sides represent a continuous linear functional
on the space of all functions $\ph$ as above which factor over $A/U$, where this
space is equipped with the norm $\norm \ph=\int_A|\ph(a)| |a^{-2\rho}|\, da$.
\end{theorem}

Let $A^*$ denote the set of all continuous group homomorphisms $\la\colon
A\ra\C^\times$. For $\la\in A^*$ and an $A$-module $V$ let $V^\la$ denote the
generalized $\la$-eigenspace, ie,
$$
V^\la\df \bigcup_{k=1}^\infty \{ v\in V\mid (a-a^\la)^k v=0\ \forall a\in A\}.
$$
Then
$$
\int_{A^-} \ph(a)\,\tr(a|H_c^q(M,\pi_N\otimes\sigma))\, da
$$
equals
$$
\sum_{\la\in A^*}\dim
H_c^q(M,\pi_N\otimes\sigma)^\la\,\int_{A^-}\ph(a)\,a^\la\, da.
$$
For $\la\in A^*$ define
$$
m_\la^{\sigma,\omega} \df \sum_{\pi\in\hat G}N_{\Ga,\omega}(\pi)\sum_{q=0}^{\dim
M}(-1)^q\,\dim H_c^q(M,\pi_N\otimes\sigma)^\la.
$$
The sum is always finite.

On the other hand, for $[\ga]\in\CE_P(\Ga)$ let
$$
c_\ga\df \la_\ga\,\chi_r(\Ga_\ga)\,|a_\ga^{2\rho}|.
$$
Then the Theorem is equivalent to the following Corollary.

\begin{corollary}
(Lefschetz Formula)\\
As an identity of distributions on $A^-$ we have
$$
\sum_{\la\in A^*} m_\la^{\sigma,\omega}\, \la\= \sum_{[\ga]\in\CE_P(\Ga)}
c_\ga\,\tr\omega(\ga)\,\tr\sigma(m_\ga)\,\delta_{a_\ga}.
$$
\end{corollary}

{\bf Proof of the Theorem:}
Let $f_{EP}$ be an Euler-Poincar\'e function on $M$ which
is $K_M$-central \cite{Kottwitz}. For $m\in M$ regular we have
$$
\CO_m^M(f_{EP})\=\begin{cases} 1 & m {\ \rm elliptic},\\ 0 & {\rm
otherwise}.\end{cases}
$$
For $g\in G$ and a finite dimensional $F$ vector space  $V$ on which $g$ acts
linearly let
$E(g|V)$ be the set of all absolute values $|\mu|$, where $\mu$ ranges over
the eigenvalues of $g$ in the algebraic closure $\bar F$ of $F$. 
Let $\la_{min}(g|V)$ denote the minimum and $\la_{max}(v|V)$ the maximum of
$E(g|V)$. For $am\in AM$ define
$$
\la(am)\df \frac{\la_{min}(a|\bar\n)}{\la_{max}(m|\g)^2}
$$
Note that $\la_{max}(m|\g)$ is always $\ge 1$ and that
$\la_{max}(m|\g)\la_{min}(m|\g)=1$. We will consider the set
$$
(AM)^\sim \ :=\ \{ am\in AM | \la(am)>1 \}.
$$
Let $M_{ell}$ denote the set of elliptic elements in $M$.

\begin{lemma} \label{MA-pad}
The set $(AM)^\sim$ has the following properties:

\begin{enumerate}
\item
$A^-M_{ell}\subset (AM)^{\sim}$
\item
$am\in (AM)^\sim \Rightarrow a\in A^-$
\item
$am, a'm' \in (AM)^\sim, g\in G\ {\rm with}\ a'm'=gamg^{-1}
\Rightarrow a=a', g\in AM$.
\end{enumerate}
\end{lemma}

\prf The first two are immediate. For the third let $am, a'm' \in
(AM)^\sim$ and $g\in G$ with $a'm'=gamg^{-1}$. Observe that by the
definition of $(AM)^\sim$ we have
\begin{eqnarray*}
\la_{min}(am | \bar{\n}) &\ge& \la_{min}(a| \bar{\n})\la_{min}(m |
\g)\\
 &>& \la_{max}(m|\g)^2\la_{min}(m|\g)\\
 &=& \la_{max}(m|\g)\\
 &\ge & \la_{max}(m | \a +\m +\n)\\ &\ge&\la_{max}(am | \a +\m +\n)
\end{eqnarray*}
that is, any eigenvalue of $am$ on $\bar{\n}$ is strictly bigger
than any eigenvalue on $\a +\m +\n$. Since $\g = \a +\m +\n
+\bar{\n}$ and the same holds for $a'm'$, which has the same
eigenvalues as $am$, we infer that $\Ad(g)\bar{\n} =\bar{\n}$. So
$g$ lies in the normalizer of $\bar{\n}$, which is
$\bar{P}=MA\bar{N} =\bar{N}AM$. Now suppose $g=nm_1a_1$ and
$\hat{m} =m_1mm_1^{-1}$ then
$$
gamg^{-1} \= na\hat{m}n^{-1} \= a\hat{m}\
(a\hat{m})^{-1}n(a\hat{m})\ n^{-1}.
$$
Since this lies in $AM$ we have $(a\hat{m})^{-1}n(a\hat{m}) =n$
which since $am\in (AM)^\sim$ implies $n=1$. The lemma is proven.
\qed

 Let $G$ act on itself
by conjugation, write $g.x = gxg^{-1}$, write $G.x$ for the orbit,
so $G.x = \{ gxg^{-1} | g\in G \}$ as well as $G.S = \{ gsg^{-1} |
s\in S , g\in G \}$ for any subset $S$ of $G$.

Fix a smooth function $\eta$ on $N$ which has compact support, is
positive, invariant under $K_M$ and satisfies $\int_N\eta(n) dn
=1$. Extend the function $\ph$ from $A^-$ to a conjugation
invariant smooth function $\tilde{\ph}$ on $AM$ such that
$\tilde{\ph}(am)=\ph(a)$ whenever $m$ is elliptic and such that
there is a compact subset $C\subset A^-$ such that $\tilde{\ph}$
is supported in $CM\cap(AM)^\sim$. This is clearly possible since
the support of $\ph$ is compact. It follows that the function
$$
am\ \mapsto\ f_{EP}(m)\,\tr\sigma(m)\,\tilde{\ph}(am)\,|a^{2\rho}|
$$
is smooth and compactly supported on $AM$. Given these data let
$f = f_{\eta ,\tau ,\ph} : H\ra \C$ be defined by
$$
f (kn ma (kn)^{-1}) := \eta (n) f_{EP}(m)
\,\tr\sigma(m)\,\tilde{\ph}(am)\,|a^{2\rho}|,
$$
for $k\in K, n\in N, m\in M, a\in\overline{A^-}$. Further $f(x)=0$
if $x$ is not in $G.(AM)^\sim$.

\begin{lemma} \label{welldefpad}
The function $f$ is well defined.
\end{lemma}

\prf By the decomposition $G=KP=KNMA$ every element $x\in
G.(AM)^\sim$ can be written in the form $kn ma (kn)^{-1}$. Now
suppose two such representations coincide, that is
$$
kn ma (kn)^{-1}\= k'n' m'a' (k'n')^{-1}
$$
then by Lemma \ref{MA-pad} we get $(n')^{-1} (k')^{-1}kn\in MA$, or
$(k')^{-1}k\in n'MAn^{-1}\subset MAN$, hence $(k')^{-1}k\in K\cap
MAN=K\cap M=K_M$. Write $(k')^{-1}k=k_M$ and $n''=k_Mnk_M^{-1}$,
then it follows
$$
n'' k_Mmk_M^{-1} a (n'')^{-1}\= n'm'a'(n')^{-1}.
$$
Again by Lemma \ref{MA-pad} we conclude $(n')^{-1}n''\in MA$, hence
$n'=n''$ and so
$$
k_Mmk_M^{-1} a \= m'a',
$$
which implies the well-definedness of $f$.
\qed

We will plug $f$ into the trace formula. For the geometric side
let $\ga \in \Ga$. We have to calculate the orbital integral:
$$
\CO_\ga (f) = \int_{G_\ga \bs G} f(x^{-1}\ga x) dx.
$$
by the definition of $f$ it follows that $\CO_\ga(f)=0$ if
$\ga\notin G.(AM)^\sim$. It remains to compute $\CO_{am}(f)$ for
$am\in(AM)^\sim$. Again by the definition of $f$ it follows
\begin{eqnarray*}
\CO_{am}(f) &=&
\CO_m^M(f_{EP})\,\tr\sigma(m)\,\tilde{\ph}(am)\,|a^{2\rho}|\\
&=& \begin{cases} \tr\sigma(m)\,\ph(a)\,|a^{2\rho}| & {\rm if}\ m\ {\rm is\
elliptic},\\ 0 & {\rm otherwise.}\end{cases}
\end{eqnarray*}
Here $\CO_m^M$ denotes the orbital integral in the group $M$.
 Recall that Proposition \ref{2.3} says
$$
\vol (\Ga_\ga\bs G_\ga) = (-1)^{q(G)+r}\,\la_\ga\,\chi_{_r}(\Ga_\ga),
$$
so that for $\ga\in\CE_P(\Ga)$,
$$
\vol(\Ga_\ga\bs G_\ga)\,\CO_\ga(f)\=
(-1)^{q(G+r}\,\la_\ga\,\chi_r(\Ga_\ga)\,\tr\sigma(m_\ga)\,\ph(a_\ga)\,|a_\ga^{2\rho}|.
$$

To compute the spectral side let $\pi\in\hat{G}$. We want to compute
$\tr\pi(f)$. Let $\Theta_\pi^G$ be the locally integrable
conjugation invariant function  on $G$ such that
$$
\tr\pi(f)\= \int_G f(x)\, \Theta_\pi^G(x) dx.
$$
This function $\Theta_\pi$ is called the \emph{character} of $\pi$.
It is known that the Jacquet module $\pi_{{N}}$ is a finitely generated admissible
module for the  group $MA$ and therefore it has a character
$\Theta_{\pi_{{N}}}^{MA}$. In \cite{Casselman} it is shown that
$$
\Theta_\pi(ma) \= \Theta_{\pi_{{N}}}^{MA}(ma)
$$
for $ma\in M_{ell}A^-$.

Let $h$ be a function in $L^1(G)$ which is supported in the set $G.MA$. Comparing
invariant differential forms as in the proof of the Weyl integration formula one
gets that the integral $\int_G h(x)\, dx$ equals
$$
\frac 1{|W(G,A)|}
\int_A\int_M\int_{G/AM} h(y amy^{-1})\,|\det(1-am | \n +\bar \n)| dy\, da\, dm,
$$
where $W(G,A)$ is the Weyl group of $A$ in $G$. 

For $a\in A^-$ and $m\in M_{ell}$ every eigenvalue of $am$ on $\n$ is of absolute
value $<1$ and $>1$ on $\bar \n$. By the ultrametric ptoperty this implies
\begin{eqnarray*}
|\det(1-am|\n +\bar \n)| &=& |\det(1-am|\bar \n)|\\
&=& |\det(am|\bar \n)|\\
&=& |\det(a|\bar \n)|\\
&=& |a^{-2\rho}|,
\end{eqnarray*}

We apply this to
$h(x)=f(x)\Theta_\pi^G(x)$ and use conjugation invariance of $\Theta_\pi^G$ to get
that $\tr\pi(f)$ equals
$$
\frac 1{|W(G,A)|}\int_{AM}
f_{EP}(m)\,\tr\sigma(m)\,\tilde\ph(am)\,\Theta_{\pi_N}^M(am)\,da\, dm,
$$
which is the same as
$$
\int_{A^-M}
f_{EP}(m)\,\tr\sigma(m)\,\tilde\ph(am)\,\Theta_{\pi_N}^M(am)\,da\, dm.
$$
We recall the Weyl integration formula for $M$. Let $(H_j)_j$ be a maximal family of
pairwise non-conjugate Cartan subgroups of $M$. Let $W_j$ be the Weyl group of $H_j$
in $M$. For $h\in H_j$ let $D_j(h)=\det(1-h|\m/\h_j)$, where $\m$ and $\h_j$ are the
Lie algebras of $M$ and $H_j$ resp. Then, for every $h\in L^1(M)$,
\begin{eqnarray*}
\int_M h(m)\, dm&=& \sum_j \frac 1{|W_j|}\int_{H_j^{reg}}\int_{M/H_j}h(mxm^{-1})\,
D_j(x)\, dm\, dx\\
&=& \sum_j \frac 1{|W_j|}\int_{H_j} \CO_x^M(h)\, dx,
\end{eqnarray*}
where $H_j^{reg}$ is the set of $x\in H_j$ which are regular in $M$. We fix
$a\in A^-$ and apply this to
$h(m)=f_{EP}(m)\tilde\ph(am)\tr\sigma(m)\Theta_{\pi_N}^M(am)$.
Since
$\tilde\ph$ is conjugation invariant we get for $x\in H_j^{reg}$,
$$
\CO_x^M(h)\=\CO_x^M(f_{EP})\,\tilde\ph(am)\,\tr\sigma(m)\,\Theta_{\pi_N}^{AM}(ax).
$$
This is non-zero only if $x$ is elliptic. If $x$ is elliptic, then $\tilde\ph(ax)$
equals $\ph(a)$. So we can replace $\tilde\ph(ax)$ by $\ph(a)$ throughout. Thus
$\tr\pi(f)$ equals
$$
\int_{A^-M}
f_{EP}(m)\,\ph(a)\,\tr\sigma(m)\,\Theta_{\pi_N}^M(am)\,da\, dm.
$$
 The trace
$\tr\pi(f)$ therefore equals
$$
\int_{A^-M} f_{EP}(m)\, \ph(a)\,\Theta_{\pi_N\otimes\sigma}(am)\, da\,
dm.
$$
We write $H_c^\bullet(M,V)$ for the continuous cohomology of $M$ with coefficients
in the $M$-module $V$. By Theorem 2 in \cite{Kottwitz},
$$
\tr(\pi_N\otimes\sigma)(f_{EP})\=\sum_{q=0}^{\dim M} (-1)^q\, \dim
H_c^q(M,\pi_N\otimes\sigma).
$$
The cohomology groups $H_c^q(M,\pi_N\otimes\sigma)$ are finite dimensional
$A$-modules and
$$
\tr\pi(f)\=\sum_{q=0}^{\dim M} (-1)^q\,
\int_{A^-}\tr(a|H_c^q(M,\pi_N\otimes\sigma))\,\ph(a)\, da.
$$
The Lefschetz Theorem follows.
\qed

\section{Rationality of the zeta function}

We keep fixed a neat cocompact discrete subgroup $\Ga$ of $G$ and a 
$F$-parabolic $P=LN$. But now we assume that the $F$-rank of $P$ is one, ie,
that $\dim A=1$. Write
$\CE_P^p(\Ga)$ for the set of all
$[\ga]\in\CE_P(\Ga)$ such that $\ga$ is primitive. Fix a unitary character $\chi$
of $A$. For a finite dimensional unitary representation
$(\omega ,V_\omega)$ of $\Ga$ consider the infinite product
$$
Z_{P,\sigma ,\chi ,\omega} (T) \ :=\ \prod_{[\ga ]\in\CE_P^p(\Ga)} \det\left( 1-T^{l_\ga}\chi(a_\ga)\omega (\ga)\otimes\sigma(m_\ga)\right)^{\chi_{_1}(\Ga_\ga)}.
$$

\begin{theorem}
The infinite product 
$$
Z_{P,\sigma ,\chi ,\omega}(s)
$$ 
converges for $|T|$ small enough. 
The limit extends to a rational function in $T$.
\end{theorem}

\prf
Formally at first we compute the logarithmic derivative as
$$
\frac{Z_{P,\sigma ,\chi ,\omega}'}{Z_{P,\sigma ,\chi ,\omega}}(T) \= \frac{d}{dT} \left(-\sum_{[\ga]\in\CE_P^p(\Ga)} \chi_{_1}(\Ga_\ga) \sum_{n=1}^\infty \frac{T^{nl_\ga}}{n} \tr \omega(\ga^n)\tr\sigma(m_\ga^n)\chi(a_\ga^n)\right)
$$ $$
= -\sum_{[\ga]\in\CE_P^p(\Ga)}\chi_{_1}(\Ga_\ga) \sum_{n=1}^\infty l_\ga \frac{T^{nl_\ga}}{T} \tr \omega(\ga^n)\tr\sigma(m_\ga^n)\chi(a_\ga^n).
$$
Since $\Ga$ is neat we have $\Ga_\ga =\Ga_{\ga^n}$ for all $n\in\N$, so
$$
\frac{Z_{P,\sigma ,\chi ,\omega}'}{Z_{P,\sigma ,\chi ,\omega}}(T) \= 
-\sum_{[\ga]\in\CE_P(\Ga)}\chi_{_1}(\Ga_\ga) l_{\ga_0} \frac{T^{l_\ga}}{T} \tr
\omega(\ga)\tr\sigma(m_\ga)\chi(a_\ga).
$$
This last expression equals the geometric side of the Lefschetz formula for the test
function
$$
\ph(a)\= \chi(a) T^{l_a}.
$$

The Lefschetz formula implies that
$$
\frac{Z_{P,\sigma ,\chi ,\omega}'}{Z_{P,\sigma ,\chi ,\omega}}(T)
$$
equals
$$ 
(-1)^{q(G)} \sum_{\pi\in\hat{G}} N_{\Ga ,\omega} (\pi)\
\frac{1}{T} \sum_{i=0}^{\dim M} (-1)^i \int_{A^-} \tr
(a|H_c^i(M,\pi_{{N}}\otimes\sigma_\chi))\  |a|^{2|\rho|} T^{v(a)} \ da.
$$
Decompose the $A$-module $H_c^i(M,\pi_{{N}}\otimes\sigma_\chi)$ into generalized
eigenspaces under $A$. Clearly only that part of
$H_c^i(M,\pi_{{N}}\otimes\sigma_\chi)$ survives on which $A$ acts by unramified
characters. Suppose this space splits as
$$
\bigoplus_{\la\in\C} E_\la^i,
$$
where for $\la\in\C$ the space $E_\la^i$ is the generalized eigenspace to the character $|a|^\la$.
Let $m_\la^i := \dim_\C E_\la^i$.
Then the contribution of $\pi\in\hat G$ equals
$$
 \sum_{i=1}^{\dim M} (-1)^i \sum_{\la\in\C} m_\la^i \int_{A^+}
|a|^{2|\rho|+\la} T^{v(a)-1}\ da
$$
\begin{eqnarray*}
	&\=& \sum_{i=1}^{\dim M} (-1)^i \sum_{\la\in\C} m_\la^i \sum_{n=1}^\infty (Tq^{2|\rho|+\la})^n\ T^{-1}\\
	&\=& -\sum_{i=1}^{\dim M} (-1)^i \sum_{\la\in\C} m_\la^i  \rez{T-q^{-2|\rho|-\la}}.
\end{eqnarray*}
Since only finitely many $\pi\in\hat{G}$ contribute, it follows that $Z_{P,\sigma ,\chi ,\omega}'/Z_{P,\sigma ,\chi ,\omega}$ is a finite sum of functions of the form $\rez{T-a}$, this implies that  $Z_{P,\sigma ,\chi ,\omega}$ is a rational function. 
The theorem follows.
\qed

Since $\dim A=1$ there is an isomorphism of topological groups $\psi\colon A\ra
F^\times$ such that $A^-$ is mapped to the set of all $x\in F^\times$ with $|x|<1$.
For $a\in A$ we will write $|a|$ for $|\psi(a)|$. This does not depend on the choice
of $\psi$. Every 
$\la\in\C$ thus gives a quasicharacter $a\mapsto |a|^\la$ and we can represent all
unramified quasicharacters in this way.  The proof of the theorem actually gives us:

\begin{proposition}
For $\la\in\C$ let $H_c^i(M,\pi_{{N}}\otimes\sigma_\chi)_\la$ denote the
generalized $|.|^\la$-eigenspace of the $A$-action on
$H_c^i(M,\pi_{{N}}\otimes\sigma_\chi)$, then the (vanishing-) order of
$Z_{P,\sigma ,\chi ,\omega}(q^{-s-2|\rho |})$ in $s=s_0$  equals
$$
(-1)^{q(G)+1}\sum_{\pi\in\hat{G}} N_{\Ga ,\omega}(\pi) \sum_{i=0}^{\dim M} (-1)^i \dim
H_d^i(M,\pi_{{N}}\otimes\sigma_\chi)(-s_0).
$$
\end{proposition}
\qed

\section{The divisor of the zeta function}

First we introduce some standard notation.
A representation $(\rho ,V)$ on a complex vector space will be called
\emph{algebraic} if every vector $v\in V$ is fixed by some compact open subgroup.
Let $Alg(G)$ be the category of algebraic representations (comp. \cite{BerZel}).

The representation $(\rho ,V)$ is called \emph{admissible} if for any compact open
subgroup $K\subset G$ the space of fixed vectors $V^K$ is finite dimensional. Let
$\CHC_{alg}(G)$ denote the full subcategory of $Alg(G)$ consisting of all admissible
representations of finite length in $Alg(G)$.

In the following we will mix the representation theoretic notion with the module-theoretic notion, so $V$ will be a $G$-module with $x\in G$ acting as $v\mapsto x.v$.

We introduce the setup of topological (continuous) representations.
Let $\CC(G)$ denote the category of continuous representations of $G$ on locally convex, complete, Hausdorff topological vector spaces.
The morphisms in $\CC(G)$ are continuous linear $G$-maps.

To any $V\in \CC(G)$ we can form $V^\infty$, the algebraic part, which is by definition the space of all vectors in $V$ which have an open stabilizer.
By the continuity of the representation it follows that $V^\infty$ is dense in $V$.
The association $V\mapsto V^\infty$ gives an exact functor from $\CC(G)$ to $Alg(G)$.
An object $V\in \CC(G)$ is called admissible if $V^\infty$ is and so we define $\CHC_{top}$ to be the full subcategory of $\CC(G)$ consisting of all admissible representations of finite length.
Let $F$ be the functor:
\begin{eqnarray*}
F : \CHC_{top} &\ra & \CHC_{alg}\\
W&\mapsto & W^\infty .
\end{eqnarray*}

Let $V$ be in $\CHC_{alg}$ then any $W\in\CHC_{top}$ with $FW=V$ will be called a
\emph {completion} of $V$.

To $V\in Alg(G)$ let $V^*=Hom_\C(V,\C)$ be the dual module and let $\tilde{V}$ be the algebraic part of $V^*$.
If $V$ is admissible we have that the natural map from $V$ to $\tilde{\tilde{V}}$ is an isomorphism.

On the space $C(G)$ of continuous maps from $G$ to the complex numbers we have two actions of $G$, the left and the right action given by
$$
L(g)\ph(x):=\ph(g^{-1}x)\ \ \  R(g)\ph(x) := \ph(xg),
$$
where $g,x\in G$ and $\ph\in C(G)$. 
For $V\in\CHC_{alg}(G)$ let
$$
V^{-\infty} := Hom_G(\tilde{V} ,C(G)),
$$
where we take $G$-homomorphisms with respect to the right action, so $V^{-\infty}$ is the space of all linear maps $f$ from $\tilde{V}$ to $C(G)$ such that $f(g.v^*)=R(g)f(g^*)$.
Then $V^{-\infty}$ becomes a $G$-module via the left translation: for $\alpha\in V^{-\infty}$ we define $g.\alpha (v^*) := L(g)\alpha(v^*)$.
We call $V^{-\infty}$ the \emph{maximal completion} of $V$.
The next lemma and the next Proposition will justify this terminology.
First observe that the topology of locally uniform convergence may be installed on $V^{-\infty}$ to make it an element of $\CC(G)$.
By $V\mapsto V^{-\infty}$ we then get a functor $R:\CHC_{alg}(G)\ra\CHC_{top}(G)$.

\begin{lemma}
We have $(V^{-\infty})^\infty \cong V$.
\end{lemma}

\prf
There is a natural map $\Phi : V\ra V^{-\infty}$ given by $\Phi(v)(v^*)(g) := v^*(g^{-1}.v)$, for $v^*\in \tilde{V}$ and $g\in G$.
This map is clearly injective.
To check surjectivity let $\alpha\in(V^{-\infty})^\infty$, then the map
$\tilde{V}\ra \C$, $v^*\mapsto \alpha(v^*)(1)$ lies in $\tilde{\tilde{V}}$.
Since the natural map $V\ra \tilde{\tilde{V}}$ is an isomorphism, there is a $v\in V$ such that $\alpha(v^*)(1)=v^*(v)$ for any $v^*\in \tilde{V}$, hence $\alpha(v^*)(g)= \alpha(g.v^*)(1)=g.v^*(v)=v^*(g^{-1}.v)=\Phi(v)(v^*)(g)$, hence $\alpha =\Phi(v)$.
\qed

It follows that the functor $R$ maps $\CHC_{alg}(G)$ to $\CHC_{top}(G)$ and that $FR=Id$.

\begin{proposition}
The functor $R: \CHC_{alg}(G)\ra \CHC_{top}(G)$ is right adjoint to $F: W\mapsto W^\infty$.
So for $W\in \CHC_{top}(G)$ and $V\in \CHC_{alg}(G)$ we have a functorial isomorphism:
$$
Hom_{alg}(FW,V)\cong Hom_{top}(W,RV).
$$
\end{proposition}

\prf
We have a natural map
\begin{eqnarray*}
Hom_{top}(W,V^{-\infty}) &\ra& Hom_{alg}(W^\infty ,V)\\
\alpha &\mapsto& \alpha |_{W^\infty}.
\end{eqnarray*}
Since $W^\infty$ is dense in $W$ this map is injective.
For surjectivity we first construct a map $\psi : W\ra (W^\infty)^{-\infty}$.
For this let $W'$ be the topological dual of $W$ and for $w'\in W'$ and $w\in W$ let
$$
\psi_{w',w}(x) := w'(x^{-1}.w),\ \ \ x\in G.
$$
The map $w\mapsto \psi_{.,w}$ gives an injection
$$
W\hookrightarrow Hom_G(W',C(G)).
$$
Let $\ph\in \tilde{FW}$, then $\ph$ factors over $(FW)^K$ for some compact open subgroup $K\subset G$.
But since $(FW)^K =W^K$ it follows that $\ph$ extends to $W$.
The admissibility implies that $\ph$ is continuous there.
So we get $\tilde{(FW)} \hookrightarrow W'$.
The restriction then gives
$$
Hom_G(W',C(G))\ra Hom_G(\tilde{(FW)},C(G)).
$$
From this we get an injection
$$
\psi : W\hookrightarrow (FW)^{-\infty},
$$
which is continuous.
Let $\zeta : W^\infty \ra V$ be a morphism in $\CHC_{alg}(G)$.
We get
$$
\alpha(\zeta) : W\hookrightarrow (FW)^{-\infty} \begin{array}{c} \zeta^{-\infty}\\ \longrightarrow\\ {}\end{array} V^{-\infty},
$$
with $\alpha(\zeta)|_{W^\infty}=\zeta$, i.e. the desired surjectivity.
\qed

The Proposition especially implies that any completion of $V\in\CHC_{alg}(G)$ injects into $V^{-\infty}$.
This assertion is immediate for irreducible modules and follows generally by induction on the length.
Thus the use of the term 'maximal completion' is justified.

Let now $\Ga$ denote a cocompact torsion-free discrete subgroup of $G$, let $(\omega ,V_\omega)$ be a finite dimensional unitary representation of $\Ga$ and let $C(\Ga \bs G,\omega)$ denote the space of all continuous functions $f$ from $G$ to $V_\omega$ satisfying $f(\ga x)=\omega(\ga)f(x)$.

For $U,V\in Alg(G)$ and $p,q\ge 0$ let $H_d^q(G,V)$ and $Ext_{G,d}^q(U,V)$ denote the differentiable cohomology and Ext-groups as in
\cite{BorWall}, chap. X.

The next theorem is a generalization of the classical duality theorem
\cite{Gelf}.

\begin{theorem}\label{higher_duality}
(Higher Duality Theorem) 
Let $\Ga$ be a cocompact torsion-free discrete subgroup of $G$ and $(\omega ,V_\omega)$ a finite dimensional unitary representation of $\Ga$, then for any $V\in \CHC_{alg}(G)$:
\begin{eqnarray*}
H^q(\Ga ,V^{-\infty}\otimes \omega) &\cong& H_d^q(G,Hom_\C(C(\Ga\bs G,\tilde{\omega})^\infty ,V)^\infty)\\
	&\cong& Ext_{G,d}^q(C(\Ga\bs G,\tilde{\omega})^\infty ,V).
\end{eqnarray*}
\end{theorem}

\prf
Let 
$$
\hat{M} := Hom_\C(\tilde{V},C(G)\otimes V_\omega).
$$
Consider the action of $G$ on $\tilde{M}$:
\begin{eqnarray*}
l(g)\alpha (v^*)(x) &:=& \alpha(g^{-1}.v^*)(xg),\\
\end{eqnarray*}

It is clear that
$$
\hat{M}^{l(G)} \cong V^{-\infty}\otimes V_\omega.
$$
Let $M$ denote the algebraic part of $\hat{M}$ with respect to $l$. i.e.
\begin{eqnarray*}
M &=& Hom_\C(\tilde{V},C(G)\otimes V_\omega)^\infty\\
&=& \lim_{\begin{array}{c}\ra\\ K\end{array}} Hom_K(\tilde{V},C(G)\otimes V_\omega)
\end{eqnarray*}
where $K$ runs over the set of all compact open subgroups of $G$.

The group $\Ga$ acts on $M$ by
$$
\ga \alpha(v^*)(x) := \omega(\ga) \alpha(v^*)(\ga x).
$$
This action commutes with $l(G)$.
Now let $\CC^\infty(G)$ be the subcategory of $\CC(G)$ consisting of \emph{smooth} representations.
This means that $V\in\CC(G)$ lies in $\CC^\infty(G)$, when $V^\infty =V$.
The inductive limit topology makes $M$ an element of $\CC^\infty(G\times \Ga)$.
Let $Ab$ denote the category of abelian groups.
The functor $\CC^\infty(G\times \Ga)\ra Ab$ defined by taking $G\times \Ga$-invariants can be written as the composition of the functors $H^0(\Ga,.)$ and $H^0(G,.)$ in two different ways:
$$
\begin{array}{ccc}
\CC^\infty(G\times \Ga) & \begin{array}{c}{H^0(\Ga ,.)}\\ \longrightarrow\\ {}\end{array}& \CC^\infty(G)\\
{ H^0(G,.)}\downarrow& {}& \downarrow H^0(G,.)\\
\CC^\infty(\Ga) &\begin{array}{c}{H^0(\Ga ,.)}\\ \longrightarrow\\ {}\end{array}& Ab,
\end{array}
$$
giving rise to two spectral sequences:
\begin{eqnarray*}
^1E_2^{p,q} &=& H^p_d(G,H^q(\Ga ,M)),\\
^2E_2^{p,q} &=& H^p(\Ga ,H^q_d(G,M)),
\end{eqnarray*}
both abutting to $H^*(G\times \Ga ,M)$.

\begin{lemma}
The module $M$ is $(G,d)$-acyclic and $\Ga$-acyclic.
\end{lemma}

When saying $(G,d)$-acyclic we mean differentiable or strong acyclicity
\cite{BorWall}. For the discrete group $\Ga$ this notion coincides
with usual acyclicity.

\prf
We show the $G$-acyclicity first.
We have
\begin{eqnarray*}
H_d^q(G,M) &=& H_d^q(G,Hom(\tilde{V},C(G)\otimes \omega)^\infty)\\
	&=& Ext_{G,d}^q(\tilde{V},C(G)^\infty\otimes\omega).
\end{eqnarray*}
By \cite{BorWall} X 1.5 we know that $C(G)^\infty$ is s-injective.
Therefore $H_d^q(G,M)=0$ for $q>0$ as desired.

For the $\Ga$-acyclicity we consider the standard resolution $C^q:=\{ f:\Ga^{q+1}\ra M\}$ with the differential $d: C^q\ra C^{q+1}$ given by
$$
df(\ga_0 ,\dots ,\ga_{q+1}) = \sum_{j=0}^{q+1} (-1)^j f(\ga_0 ,\dots ,\hat{\ga_j}, \dots ,\ga_{q+1}).
$$
The group $\Ga$ acts on $C^q$ by
$$
(\ga f) (\ga_0 ,\dots ,\ga_q) = \ga .f(\ga^{-1}\ga_0 ,\dots ,\ga^{-1}\ga_q).
$$
Then $H^*(\Ga ,M)$ is the cohomology of the complex of $\Ga$-invariants $(C^q)^\Ga$.

Choose an open set $U\subset G$ such that $(\ga U)_{\ga\in\Ga}$ is a locally finite open covering of $G$.
Then there exists a locally constant partition of unity $(\rho_\ga)_{\ga\in\Ga}$ such that $\supp \rho_\ga\subset \ga U$ and $L_{\ga'}\rho_\ga =\rho_{\ga'\ga}$.
Consider a cocycle $f\in(C^q)^\Ga$, $df=0$, $q\ge 1$.
Let $F\in C^{q-1}$ be defined by
$$
F(\ga_0,\dots ,\ga_{q-1})(x) := \sum_{\ga\in\Ga} f(\ga_0 ,\dots ,\ga_{q-1},\ga)(x)\rho_\ga(x),\ \ \ x\in G.
$$
The sum is locally finite in $x$.
A computation shows that $F$ is $\Ga$-invariant and that $dF=(-1)^qf$.
The lemma follows.
\qed

By the lemma both our spectral sequences only live on the $x$-axis, ie, they only have a nonzero entry  at $(p,q)$ if $q=0$.
Since they have the same abutment, they agree. Hence
\begin{eqnarray*}
H^p(\Ga ,H^0(G,M)) &=& ^1E_2^{p,0}\\
	&=& ^2E_2^{p,0}\\
	&=& H^p_d(G,H^0(\Ga ,M)).
\end{eqnarray*}
The first term is $H^p(\Ga,V^{-\infty}\otimes V_\omega)$.
The last is
$$
H_d^q(G ,Hom_\C (\tilde{V},C(\Ga\bs G,\omega))^\infty).
$$
Since $V$ is of finite length we get an isomorphism as $G$-modules:
$$
Hom_\C(\tilde{V} ,C(\Ga \bs G,\omega))^\infty \cong \bigoplus_{\pi\in\hat{G}} N_{\Ga ,\omega}(\pi) Hom_\C (\tilde{V},\pi^\infty).
$$

We now employ the following

\begin{lemma}
Let $V,W\in\CHC_{alg}$, then as  $G$-modules,
$$
Hom_\C(V,W)^\infty \cong Hom_\C(\tilde{W},\tilde{V})^\infty.
$$
\end{lemma}

\prf
We have a natural map $Hom_\C(V,W)^\infty \ra Hom_\C(\tilde{W},\tilde{V})^\infty$ given by dualizing: $\alpha\mapsto \tilde{\alpha}$.
This map is injective.
Iterating, we get a map $Hom_\C(\tilde{V},\tilde{W})^\infty \ra Hom_\C(\tilde{\tilde{W}},\tilde{\tilde{V}})^\infty\cong Hom_\C(V,W)^\infty$, and these maps are inverse to each other.
\qed

Recall that $C(\Ga\bs G,\omega)^\infty$ is a direct sum of irreducibles, each occurring with finite multiplicity.
So every $G$-morphism joining $C(\Ga\bs G,\omega)^\infty$ with a module of finite length will factor over a finite sum of irreducibles in $C(\Ga\bs G,\omega)^\infty$.
Therefore the above applies to give
$$
Hom_\C(\tilde{V} ,C(\Ga \bs G,\omega)^\infty) \cong Hom_\C(C(\Ga\bs G,\tilde{\omega})^\infty ,V),
$$
where we additionally have used that the invariant integral puts $C(\Ga\bs G,\omega)^\infty$ in duality with $C(\Ga\bs G,\tilde{\omega})^\infty$.
The theorem is proven.
\qed

Let by a slight abuse of notation $\pi_{\sigma ,\chi ,s_0}$ be the representation parabolically induced from $\sigma \chi |.|^{s_0}$, i.e. the space of this representation consists of measurable functions
$f:G\ra V_\sigma$ such that $f(manx)=\chi(a) |a|^{s_0 +|\rho |}\sigma(m) f(x)$, $f$ is square-integrable on $K$ modulo functions vanishing on a set of measure zero.

The next theorem is the Patterson conjecture as adapted to our situation.

\begin{theorem} 
Assume that the Weyl-group $W(G,A)$ is non-trivial.
The vanishing order of the function $Z_{P,\sigma ,\chi ,\omega}(q^{-s-|\rho|})$ at $s=s_0$ is
$$
(-1)^{q(G)+1}\chi_{_1}(\Ga ,\pi_{\sigma ,\chi ,s_0}^{-\infty} \otimes \omega)
$$
if $s_0\ne 0$ and
$$
(-1)^{q(G)+1}\chi_{_1}(\Ga ,\hat{\pi}_{\sigma ,\chi ,0}^{-\infty} \otimes \omega)
$$
if $s_0=0$, where $\hat{\pi}_{\sigma ,\chi ,0}^{-\infty}$ is a certain extension of $\pi_{\sigma ,\chi ,0}^{-\infty}$ with itself.

Further we have that the usual Euler-characteristic $\chi(\Ga ,\pi_{\sigma ,\chi ,-s_0}^{-\infty} \otimes \omega)$ always vanishes.
\end{theorem}

\prf
So far we know that the vanishing order of $Z_{P,\sigma ,\chi ,\omega}(q^{-s-2|\rho|})$ at $s=s_0$ equals
$$
(-1)^{q(G)+1}\sum_{\pi\in\hat{G}} N_{\Ga ,\omega}(\pi) \sum_{i=0}^{\dim M} (-1)^i \dim
H_d^i(M,\pi_{{N}}\otimes\sigma_\chi)(-s_0).
$$

Note that $H_d^i(M,\pi_{{N}}\otimes\sigma_\chi)(-s_0) =
H_d^i(M,(\pi_{{N}}\otimes\sigma_\chi)(-s_0))$.

\begin{lemma}
Let $\pi\in\hat{G}$.
If $s\ne -|\rho |$, then $A$ acts semisimply on $(\pi_{{N}})(s)$.
In any case the $L$-module $\pi_{{N}}$ has length at most $2$.
\end{lemma}

\prf
Suppose $\pi_{{N}}\ne 0$. Then there is an induced representation
$I=Ind_{{P}}^G(\xi\otimes 1)$ such that $\pi$ injects into $I$. Here $\xi$ is an
irreducible representation of $L$. Suppose $A$ acts on $\xi $ by the character $\mu$.
By the exactness of the Jacquet functor it follows that
$\pi_{{N}}$ injects into
$I_{{N}}$. It therefore suffices to show the assertion for induced
representations of the form $\pi =I$. These are given on spaces of locally constant
functions $f:G\ra V_{\xi ,\mu}$ such that $f(l{n}x)= \xi (l) f(x)$. The
map $\psi : I\ra V_{\xi}$, given by $f\mapsto f(1)$ is a
${P}$-homomorphism. The group $A$ acts on the image by $a^\mu$ and this
defines an irreducible quotient of $\pi_{{N}}$. We will show that $\ker \psi$
gives an irreducible subrepresentation of $\pi_{{N}}$ on which $A$ acts by
$a^{-\mu +2|\rho|}$. Bruhat-Tits theory tells us that $G$ splits into two
${P}$-double cosets $G={P} \cup {P}w{P} = {P} \cup
{P}w{N}$, where $w$ is the nontrivial element of the Weyl-group $W(G,A)$.
Let $f\in\ker\psi$. Then $f$ vanishes in a neighborhood of $1$ and thus the integral
$\int_{{N}}f(w{n})\,d{n}$ converges. (For this recall that by the
Bruhat-decomposition ${P}\bs G$ is the one point compactification  of ${N}$.
Therefore $f(1)=0$ implies that the function ${n}\mapsto f(w{n})$ has
compact support.)  This integral defines a map $(\ker\psi)_{{N}}\ra V_{\xi}$.
The group $A$ acts as
\begin{eqnarray*}
\int_{{N}}f(w{n}a)d{n} &=& a^{-2|\rho|}
\int_{\bar{N}}f(wa{n})d{n}\\
	&=& a^{-2|\rho|} \int_{{N}}f(a^{-1}w{n})d{n}\\
	&=& a^{-\mu -2|\rho|} \int_{{N}}f(w{n})d{n}.
\end{eqnarray*}
The claim is now clear.
\qed

To show the theorem consider the case $s\ne -|\rho|$.
Then by semisimplicity it follows
$$
H_d^q(M,(\pi_{{N}}\otimes\sigma_\chi)_{-s_0}) =
H^0(A,H_d^q(M,\pi_{{N}}\otimes\sigma_{\chi s_0})),
$$
where $\sigma_{\chi s_0}$ is the $AM$-module given by $\sigma_{\chi s_0}(am)=|a|^{s_0}\chi(a)\sigma(m)$.
Similar to \ref{chichi1} we get that
$$
\sum_{q=0}^{\dim M} (-1)^q \dim H^0(A,H_d^q(M,\pi_{{N}}\otimes\sigma_{\chi s_0}))
$$
equals
$$
 \sum_{q=0}^{\dim AM} q(-1)^{q+1} \dim H_d^q(AM,\pi_{{N}}\otimes\sigma_{\chi
s_0}).
$$
By \cite{BorWall} p.262 we know
\begin{eqnarray*}
H_d^q(AM ,\pi_{{N}}\otimes \sigma_{\chi s_0})&\cong& H_d^q(AM
,Hom_\C(\tilde{(\pi_{{N}})}, \sigma_{\chi s_0}))\\
	&\cong& Ext_{AM,d}^q((\tilde{\pi})_N,\sigma_{\chi s_0})
\end{eqnarray*}

Now recall that the functors $(.)_N$ and $Ind_{MAN}^G(.\otimes 1)$ are
adjoint to each other, i.e.
$$
Hom_{AM,d}((.)_N,.)\cong Hom_{G,d}(.,Ind_{MAN}^G(.\otimes 1)).
$$
This especially implies that $Ind_{MAN}^G(.\otimes 1)$ maps injectives to injectives
and that $(.)_N$ maps projectives to projectives. Further they are both exact, hence
preserve resolutions. It follows that for any $q$:
\begin{eqnarray*}
Ext_{AM,d}^q((\tilde{\pi})_N,\sigma_{\chi s_0})
	&\cong& Ext_{G,d}^q(\tilde{\pi},Ind_{MAN}^G(\sigma_{\chi s_0}))\\
	&\cong& Ext_{G,d}^q(\tilde{\pi},\pi_{\sigma ,\chi ,{s_0}-|\rho|}).
\end{eqnarray*}

Therefore we have that the order $ord_{s=s_0} Z_{P,\sigma ,\chi
,\omega}(q^{-s-2|\rho|})$ equals
\begin{eqnarray*}
&=& (-1)^{q(G)+1}\sum_{\pi\in\hat{G}} N_{\ga ,\omega}(\pi) \sum_{q=0}^{\dim AM} q(-1)^{q+1} \dim Ext_{G,d}^q(\tilde{\pi},\pi_{\sigma ,\chi ,s_0-|\rho|})\\
	&=& (-1)^{q(G)+1}\sum_{q=0}^{\dim AM} q(-1)^{q+1} \dim Ext_{G,d}^q(C(\Ga\bs G,\tilde{\omega})^\infty,\pi_{\sigma ,\chi ,s_0-|\rho|})\\
&=& (-1)^{q(G)+1}\sum_{q=0}^{\dim AM} q(-1)^{q+1} \dim H^q(\Ga,\pi_{\sigma ,\chi ,s_0-|\rho|}^{-\infty}\otimes\omega),
\end{eqnarray*}
the last equality comes by \ref{higher_duality}.

In the case $s=-|\rho|$ finally, one replaces $A$ by $A^2=v^{-1}(2\Z)$, the set of all $a\in A$ with even valuation.
Then $\pi_{\sigma ,\chi ,s_0-|\rho |}$ is replaced by $Ind_{A^2MN}^G(\sigma_{\chi s_0})$, which is the desired extension.
\qed

\input{hr.ind}
\end{document}